\documentclass[final,5p,times,twocolumn,authoryear]{elsarticle}

\usepackage{microtype}
\usepackage{graphicx}
\usepackage{booktabs} 

\usepackage{url}

\usepackage{amsmath}
\usepackage{amssymb}
\usepackage{mathtools}
\usepackage{amsthm}

\usepackage{gensymb}

\usepackage{caption}
\usepackage{subcaption}

\usepackage{algorithm}
\usepackage{algorithmic}

\theoremstyle{plain}

\theoremstyle{definition}

\theoremstyle{remark}

\usepackage[textsize=tiny]{todonotes}

\journal{Energy and Buildings}

\begin{document}

\begin{frontmatter}



\title{Employing Federated Learning for Training Autonomous HVAC Systems}


\author[sys]{Fredrik Hagström}


\affiliation[sys]{organization={Department of Mathematics and Systems Analysis, Aalto University},
            city={Espoo},
            country={Finland}}

\author[comp,yaiyai]{Vikas Garg}
\author[sys]{Fabricio Oliveira}

\affiliation[comp]{organization={Department of Computer Science, Aalto University},
            city={Espoo},
            country={Finland}}

\affiliation[yaiyai]{organization={YaiYai Ltd},
            city={Espoo},
            country={Finland}}

\begin{abstract}
Buildings account for 40 \% of global energy consumption. A considerable portion of building energy consumption stems from heating, ventilation, and air conditioning (HVAC), and thus implementing smart, energy-efficient HVAC systems has the potential to significantly impact the course of climate change. In recent years, model-free reinforcement learning algorithms have been increasingly assessed for this purpose due to their ability to learn and adapt purely from experience. They have been shown to outperform classical controllers in terms of energy cost and consumption, as well as thermal comfort. However, their weakness lies in their relatively poor data efficiency, requiring long periods of training to reach acceptable policies, making them inapplicable to real-world controllers directly. 

In this paper, we demonstrate that using federated learning to train the reinforcement learning controller of HVAC systems can improve the learning speed, as well as improve their ability to generalize, which in turn facilitates transfer learning to unseen building environments. In our setting, a global control policy is learned by aggregating local policies trained on multiple data centers located in different climate zones. The goal of the policy is to simultaneously minimize energy consumption and maximize thermal comfort. We perform a thorough set of experiments, evaluating three different optimizers for local policy training, as well as three different federated learning algorithms against two alternative baselines. We demonstrate through experimental evaluation that these effects lead to a faster learning speed, as well as greater generalization capabilities in the federated policy compared to any individually trained policy. Furthermore, the learning stability is significantly improved, with the learning process and performance of the federated policy being less sensitive to the choice of parameters and the inherent randomness of reinforcement learning. 
\end{abstract}




\end{frontmatter}


\section{Introduction}
\label{sec:Intro}
One of the greater challenges of modern society is that of climate change. Efforts to mitigate climate change must focus not only on the supply side of energy, e.g., renewable and nuclear energy, but also on the demand side, considering factors such as energy consumption and efficiency \citep{fawzy2020strategies}. Of the global energy consumption, buildings alone are responsible for roughly 40 \% of the total consumption \citep{biemann2021experimental}. Heating, ventilation and air conditioning (HVAC) are major factors in building energy consumption \citep{fawzy2020strategies}, and hence, developing smart and energy-efficient HVAC control systems can play an important role in mitigating climate change.  

Most of the current HVAC systems in residential buildings are managed by classical algorithms, such as rule-based controllers and proportional, integral and derivative controllers \citep{biemann2021experimental}. These controllers not only lack knowledge of the thermal dynamics of the building environment but are also unable to take weather predictions into account. Hence, they are unable to react and adapt to changes in the environment, leading to sub-optimal energy performance \citep{wang2020reinforcement}. To utilize predictive data and knowledge of the building environment for improved building control performance, one can rely on Model Predictive Control (MPC) techniques \citep{wang2020reinforcement}. MPC can anticipate when to, e.g., preheat a building based on weather and occupant forecasts, in order to improve energy efficiency. MPC has been shown to be effective at reducing energy consumption on both simulated and real building environments \citep{wang2020reinforcement}. However, a serious drawback to MPC is that it requires accurate models of the environment in which the controller operates. On the other hand, every building is unique and, as such, developing a general MPC-based energy management system that can be deployed to various buildings is extremely difficult, and MPC is yet to be adopted by the building industry on a wider scale \citep{wang2020reinforcement}.

In recent years, through the emergence and rapid development of deep learning, it has become increasingly popular to apply machine learning techniques in multiple different research fields \citep{perera2021applications}. Reinforcement learning, a sub-field of machine learning concerned with control problems, has also started to gain considerable interest in research on energy system applications, including HVAC control systems. In particular, model-free reinforcement learning algorithms provide a promising direction for building control. As the name suggests, these algorithms do not require any model of the building environment or of its dynamics within. Instead, they learn purely from data collected while interacting with the environment. This eliminates the need for expert domain knowledge to develop models of the environment and allows the algorithms to be applied to any building in general, which are the main challenges of developing and deploying MPC-based controllers. Reinforcement learning algorithms also have greater adaptability to changes in the environment, as they can learn from the environment indefinitely, and as such, they can take into account long-term changes, such as changes in climate and occupant behaviour. 

Reinforcement learning has been successfully applied for building-related control tasks, though mostly in simulated environments \citep{wang2020reinforcement}. A major hurdle for the deployment of reinforcement learning algorithms to real buildings is their poor data efficiency. They need to collect large amounts of experience data to learn decision policies that take reasonable actions, and thus, it takes a long time to train them. For example, soft actor-critic, a state-of-the-art algorithm with best-performant data efficiency and learning speed requires more than a year of training to produce an acceptable policy in terms of thermal comfort \citep{biemann2021experimental}. Currently, this data inefficiency makes it infeasible to train reinforcement learning algorithms directly in physical building environments. A promising approach to overcome this data efficiency is to use transfer learning, i.e., to pre-train a controller on a simulated environment and then move it to a real environment for fine tuning \citep{wang2020reinforcement}. Still, it is not known how to generalize a controller trained on a small set of buildings for use in another building not seen during training \citep{wang2020reinforcement}.

Our main contribution is to address this gap by investigating how federated learning can improve data efficiency and generalization in reinforcement learning-based HVAC control. To the best of our knowledge, this is the first study to systematically evaluate the impact of federated optimization on the learning dynamics and performance of reinforcement learning agents in a distributed HVAC control setting. 

Federated learning is a paradigm for decentralized distributed machine learning \citep{mcmahan2017communication}. A shared global model is trained on data distributed locally over a network of participating nodes by sending copies of the global model to the nodes, training the copies on the local data, and sending the local updates back to a central server for model aggregation. The local data of each node is never explicitly shared with other nodes, nor with the central server. This reduces the communication costs associated with transmitting data and eliminates the need for large storage capacity at the coordinating central server, while simultaneously ensuring a higher degree of data privacy at the nodes. Federated learning also makes no assumptions about the distribution of the data, and thus it can be applied to systems with heterogeneous components. These features make federated learning an ideal distributed learning scheme for smart HVAC system controllers, since every building will have its own unique data distribution, and sensitive information, e.g., occupancy behaviour, will be kept private. By training a controller on multiple buildings simultaneously, we effectively collect the total experience data at a higher rate than any single building, which counteracts the low data efficiency of reinforcement learning algorithms. Also, since the data distribution is heterogeneous, the collected experience data varies from building to building, leading the total experience to be more diverse, thereby facilitating greater generalization capabilities in the shared controller.

In this paper, we demonstrate the effectiveness of federated learning for training reinforcement learning-based HVAC controllers. In a real-world deployment, the proposed system would consist of a network of HVAC controllers operating across multiple buildings, each equipped with local reinforcement learning agents. These agents would collect sensor and forecasted data (e.g., temperature, air relative humidity, and energy consumption; for a complete list of those used in our experiments, please see Table \ref{tab:obsVar}) and update their control policies accordingly. Instead of training in isolation, the agents would participate in a federated learning framework, where local updates are periodically aggregated at a central server to refine a global control policy. This global policy is then redistributed to each building, aiming to improve learning efficiency and enable generalization across different environments.

We perform an experimental evaluation of a federated controller trained in multiple simulated data center environments using the \textit{Federated Averaging} algorithm \citep{mcmahan2017communication}. The objective of the controller is to minimize energy consumption while maintaining thermal comfort, i.e., keeping the temperature within a user-specified (deemed acceptable) range of values. We evaluate and compare the performance of three different optimizers on the local nodes: \textit{stochastic gradient descent}, \textit{stochastic gradient descent with momentum}, and \textit{Adam} \citep{kingma2014adam}. The performance of the federated controller with the best local optimizer is then compared to that of individual controllers trained exclusively on each respective data center. Furthermore, we evaluate two additional federated learning algorithms: \textit{Federated Averaging with server momentum} \citep{hsu2019measuring} and \textit{FedAdam} \citep{reddi2020adaptive}. We apply a \textit{gradient masking} technique \citep{tenison2022gradient} to each federated algorithm to improve learning stability. The reinforcement learning algorithm used is the Soft Actor-Critic (SAC) algorithm \citep{haarnoja2018softfirst, haarnoja2018softapp, haarnoja2018softlearning}, which has previously shown to outperform other alternatives in HVAC control tasks \citep{biemann2021experimental, hagstrom2023using}. Our main findings from applying federated learning to train a reinforcement learning HVAC control agent are:

\begin{itemize}
    \item Improved generalization: The federated control agent outperforms all individual agents when applied to an unseen environment.
    \item Increased learning speed: The federated control agent is able to converge to the best policy at a faster rate than an individually trained agent.
    \item Improved learning stability: There is an inherent randomness to the training process of reinforcement learning agents. Federated learning reduces the variance across different training runs, leading to more consistent results.
    \item Benefits of adaptivity: The federated training process can benefit from adaptivity on the local optimizers, as Adam outperforms stochastic gradient descent with and without momentum.
\end{itemize}

The rest of this paper consists of the following parts. Section \ref{sec:RelWork} presents an overview of the related literature. Section \ref{sec:Methodology} discusses the methodology used in this paper by presenting the key technical aspects related to the SAC algorithm, as well as those related to federated learning. Section \ref{sec:Experiments} describes the simulated environment and the setup of our experiments, as well as the results obtained. Finally, we draw conclusions in section \ref{sec:Conclusion}.

\section{Related Work}
\label{sec:RelWork}

In this section, we provide a non-exhaustive survey of the literature on the employment of reinforcement learning for HVAC control. The purpose is to provide a general overview of different techniques available, whilst delineating our contribution to the literature. For more extensive surveys, we refer the reader to the works of \cite{vazquez2019reinforcement, wang2020reinforcement, perera2021applications} and \cite{weinberg2022review}.

\subsection{Early approaches}

Some of the first applications of reinforcement learning to the control of HVAC systems were made around the turn of the millennium. \cite{anderson1997synthesis} combined a proportional plus integral (PI) controller with a reinforcement learning component to control a heating coil. They evaluated it in a simulated environment, showing improved performance compared to the PI controller alone. \cite{mozer1998neural} utilized reinforcement learning in the control of the HVAC, Domestic Hot Water (DHW) and lighting systems of a real house, with the objective of minimizing both electricity cost and occupant discomfort. In \cite{henze2003evaluation}, the authors investigated a reinforcement learning solution for the operation of a simulated thermal storage system to reduce energy costs, showing favorable results when compared to conventional controllers. \cite{liu2005evaluation} used Q-learning to train both passive and active thermal storage controllers for reduced energy costs. They found the performance to be sensitive to the learning parameters and the sizes of the state and action spaces. The training time was also observed to be unacceptably long for real-world applications. They followed up their research with a hybrid learning approach in \cite{liu2006experimental-1, liu2006experimental-2}, where the agent is first pre-trained in a simulation of the environment, after which it is applied to and further trained on the true environment, making it an early example of transfer learning. They found the approach to significantly reduce the training time needed in the true environment. However, this approach requires an accurate model of the environment for the simulation phase, therefore having the same drawback as MPC. 

\subsection{Value-based approaches}

The last decade has seen an increase in research on reinforcement learning in the energy domain, 
\citep{vazquez2019reinforcement, perera2021applications}. \cite{sun2013event} minimized the day-ahead energy costs using an event-based approach, where the reinforcement learning agent takes actions only ``as needed'', instead of in regular time intervals. This reduces computational requirements while maintaining similar performance in cost savings and human comfort compared to time-based approaches. \cite{barrett2015autonomous} reduced the cost of energy while meeting the temperature set-point specified by the user during periods of occupancy. They employ a Bayesian learning approach to predict occupancy and a Q-learning agent to control the thermostat unit. \cite{li2015multi} trained a Q-learning agent to simultaneously minimize energy consumption and maximize thermal comfort. They improve upon the learning speed of standard Q-learning by utilizing a multi-grid approach, where the discretization of the state and action spaces are highly coarse at the beginning for early convergence, after which both spaces are iteratively refined during training for more fine control of the HVAC system. \cite{ruelens2016residential} minimized the energy costs of thermostatically controlled loads in both a dynamic pricing and day-ahead scheduling scenario using a Fitted Q-iteration controller equipped with a backup controller to ensure comfort. The controller converges much faster than standard Q-learning, and yields significant cost savings compared to the default controller, though increasing the energy consumption. A similar approach was taken by \cite{costanzo2016experimental}. \cite{wei2017deep} minimized the energy costs and thermal comfort violations of a multi-zone building using a Deep Q-network (DQN), which achieves comparable levels of comfort violations while yielding greater cost savings than standard Q-learning.

The papers reviewed thus far focus on value-based reinforcement learning, in most cases Q-learning. Their limitations lie in that they must discretize the state and actions spaces, and scale poorly in terms of computation and memory to both the increase in dimension of the space and the granularity of the features \citep{wiering2012reinforcement, kochenderfer2022algorithms}. Hence, in practice, the discretization is often coarse. HVAC control tasks are often naturally formulated as continuous control problems. Modeling the control task with value-based methods can therefore lead to oversimplification, since the rough discretization of state and action spaces sacrifices finer control. In contrast, policy gradient and actor-critic algorithms learn continuous policy functions and, as such, can provide more suitable alternatives.

\subsection{Policy-based approaches}

Policy gradient and actor-critic algorithms, while applicable to continuous control, have not seen nearly as much interest as value-based methods in the HVAC control literature, as well as 
building control in general \citep{vazquez2019reinforcement, wang2020reinforcement}. This is likely due to earlier algorithms being either difficult to train due to high hyperparameter sensitivity or having poor data efficiency, making them unfeasible for any potential real-world application \citep{biemann2021experimental}. 

Still, policy-based and actor-critic methods have been the algorithm of choice in some applications. \cite{gao2019energy} combined a deep neural network for thermal comfort prediction with Deep Deterministic Policy Gradient (DDPG) to control an HVAC system. DDPG is shown to outperform the value-based Q-learning, SARSA and DQN algorithms in terms of energy consumption and thermal comfort. A similar comparison and conclusion was made between DDPG and DQN in \cite{du2021intelligent}. \cite{biemann2021experimental} evaluated and compared the performance of four actor-critic algorithms; Trust Region Policy Optimisation (TRPO), Proximal Policy Optimisation (PPO), Twin Delayed DDPG (TD3) and SAC, which have received little attention in the energy domain, despite their success in other domains \citep{perera2021applications}. \cite{biemann2021experimental} concluded that while all four algorithms reduce energy consumption compared to their model-based baseline controller, SAC provides the best trade-off between energy savings and thermal comfort, while simultaneously displaying significantly greater learning speed and stability. In \cite{chen2020gnu}, PPO was used to reduce energy consumption while maintaining thermal comfort. The control policy was pre-trained on historical data of the existing controller using imitation learning. Thus, the policy learns to emulate the existing controller, performing reasonably well already at deployment, and quickly improving through fine-tuning with the PPO algorithm. The performance was evaluated in both simulated environments and a real conference room. The pre-trained PPO controller managed to reduce the cooling demand in the real environment, making the approach reasonable for real-world deployment, assuming the existence of historical controller data.

\subsection{Model-based approaches}

Model-based reinforcement learning approaches have also been explored, albeit the role the model plays varies. For example, in \cite{gao2023comparative}, a model of the environment was learned through function approximation. The learned model is used to generate additional simulated experience in conjunction with the real experience, leading to faster convergence of the reinforcement learning algorithm. \cite{nagy2018deep} also learned a model of the environment, but instead used the model to plan the actions multiple steps ahead. While model-based approaches demonstrate greater sample efficiency than model-free algorithms, leading them to learn significantly faster, their success depends on how accurate the model is. In \cite{nagy2018deep}, their model-based algorithm converges in only about 20 days, while simultaneously outperforming the model-free approach in terms of both consumption and comfort. However, they showed that if the learned model is incorrect or if the dynamics of the environment change, the algorithm fails to adapt and is in turn outperformed by the model-free algorithm. As with MPC-based approaches, the main drawback of model-based approaches is that they require accurate models to achieve successful performance. As the dynamics of different buildings vary greatly and are difficult to model, developing model-based control systems that can be deployed generally is a challenging task.

\subsection{Federated learning in the building domain}

Federated learning has seen some application in the building energy domain. \cite{khalil2021federated} used Federated Averaging to train a thermal comfort predictive model, which is used as input for a rule-based temperature set-point controller. They follow up in \cite{khalil2022federated} with a modified implementation of Federated Averaging for reduced overhead in communication. \cite{guo2020towards} used federated learning to train machine learning models to predict the coefficient of performance of a chiller. \cite{gao2021decentralized} trained a federated model for forecasting the energy demand of buildings. \cite{lu2023residential} also take a federated approach to residential energy consumption forecasting, incorporating a reinforcement learning agent to assign weights to each local model when performing model aggregation. In \cite{wang2022privacy}, federated learning was used to train a model for regulation capacity evaluation of an HVAC system. \cite{lee2021privacy} used a federated reinforcement learning model to schedule the energy consumption of the HVAC systems of three buildings with solar photovoltaic systems and a shared controllable energy storage system. In \cite{fujita2022federated}, a similar approach to ours was taken, training a SAC agent for HVAC control using Federated Averaging, though in a notably different setting. They evaluate two different scenarios. In their power-saving scenario, the temperature setting of the AC is fixed, and the task of the agent is to turn the AC on when people are present in a room and off when the room is empty. In the second, normal operation scenario, the agent also aims to control the charging and discharging of a storage battery, with the goal of maintaining the temperature below a threshold. The agent is able to perform in the power-saving scenario, and \cite{fujita2022federated} observe an increase in the rate of convergence when using federated learning, but the agent is unable to achieve ideal control in the normal operation scenario.

Our survey suggests that, while there has been effort dedicated to the employment of reinforcement learning for controlling HVAC systems with a degree of success, there is a lack of focus on investigating whether federated learning can be used to address some of the challenges faced by these studies, such as data efficiency and generalization. This is precisely where our contribution lies. To the best of our knowledge, our work is the first to thoroughly evaluate the effects federated optimization has on the learning and performance of reinforcement learning agents for direct control of HVAC systems.

\section{Methodology}
\label{sec:Methodology}

\subsection{Reinforcement Learning}
\label{sec:RL}
Reinforcement learning is, in its essence, a computational paradigm where how to optimize a decision-making problem is ``learned by doing'' \citep{sutton2018reinforcement}. The two main components of reinforcement learning are the agent and the environment. The agent aims to learn how to optimally interact with the environment in which it exists through trial and error. The agent-environment interaction follows a \textit{Markov Decision Process} (MDP), which is a stochastic control process that evolves in a sequence of discrete time steps $t \in \mathbb{N}$. An MDP can be formally represented as a tuple ($\mathcal{S}, \mathcal{A}, R, P$), where
\begin{itemize}
    \item $\mathcal{S}$ is the \textit{state space}, i.e., the set of possible states $s$,
    \item $\mathcal{A}$ is the \textit{action space}, i.e., the set of possible actions $a$,
    \item $R : \mathcal{S} \times \mathcal{A} \rightarrow \mathbb{R}$ is the \textit{reward function} $R(s_t, a_t)$,
    \item $P : \mathcal{S} \times \mathcal{S} \times \mathcal{A} \rightarrow [0, 1]$ is the \textit{transition probability function} $P(s' | s, a)$.
\end{itemize}
At time step $t$, the agent chooses and performs an action $a_t$ based on the current state $s_t$. The environment then transitions to state $s_{t+1}$ following the dynamics of the environment described by the transition probability function $P(s_{t+1} | s_t, a_t)$. As a consequence of its actions, the agent receives a reward $r_t = R(s_t, a_t)$, which measures the quality of the chosen action. This process continues in the same way, resulting in a sequence of states and actions:
\begin{align*}
    (s_0, a_0, s_1, a_1, s_2, a_2, ...).
\end{align*}
This sequence is known as a \textit{trajectory}. In the literature, it is also commonly referred to as an episode or a rollout.

To decide what action to take in state $s_t$, the agent follows a so-called \textit{policy} $\pi$. The policy can be either deterministic or stochastic. A deterministic policy is defined as a mapping $\pi : \mathcal{S} \rightarrow \mathcal{A}$, such that $a_t = \pi (s_t)$. A stochastic policy is a probabilistic function $\pi : \mathcal{A} \times \mathcal{S} \rightarrow [0, 1]$, where $a_t \sim \pi(\cdot | s_t)$ and $\sum_{a_t \in \mathcal{A}} \pi (a_t | s_t) = 1$.

\subsubsection{Soft Actor-Critic}
\label{sec:SAC}

\textit{Soft Actor-Critic} (SAC) \citep{haarnoja2018softfirst, haarnoja2018softapp, haarnoja2018softlearning} is a state-of-the-art deep reinforcement learning algorithm that learns a continuous stochastic policy. It is \textit{model-free}, meaning that the policy is learned without knowledge of the transition dynamics $P$. It is also \textit{off-policy}, meaning that it can learn from experience samples generated by any arbitrary policy, making it more sample-efficient than \textit{on-policy} algorithms, which can only utilize samples collected from the current policy. These factors make SAC a suitable option for HVAC control, where environment dynamics are difficult to model, collecting experience is time expensive, and continuous actions allow for finer control.
 
The objective in classical reinforcement learning is to find the policy $\pi$ that maximizes the expected \textit{return}, i.e., the expected sum of rewards $\sum_{t} \mathbb{E}_{(s_t, a_t) \sim p_\pi} \Big[ R(s_t, a_t) \Big]$, where $p_\pi$ refers to the state-action marginal of the trajectory distribution induced by $\pi$. The SAC algorithm considers instead an alternative maximum-entropy objective by adding an entropy term to the expectation as follows
\begin{align}
    J(\pi) &= \max_{\pi} \sum_{t} \mathbb{E}_{(s_t, a_t) \sim p_\pi} \Big[ \big(R(s_t, a_t) - \alpha \log \pi (\cdot | s_t)) \Big],
\end{align}
where $\alpha \in [0, \infty)$ is the temperature variable that controls the trade-off between exploration (entropy) and exploitation (reward maximization).

\cite{haarnoja2018softfirst, haarnoja2018softapp, haarnoja2018softlearning} derived the SAC algorithm from an algorithm called \textit{Soft Policy Iteration} (SPI). SPI learns a policy by repeating two main steps: \textit{policy evaluation} and \textit{policy improvement}. The policy evaluation step evaluates the soft \textit{action-value function} $Q: \mathcal{S} \times \mathcal{A} \rightarrow \mathbb{R}$ of the current policy $\pi$, i.e., the expected return of starting in state $s$, taking action $a$, and adhering to the policy thereafter. The soft Q-value is evaluated by iteratively updating the soft Q-function until convergence according to the soft \textit{Bellman equation}
\begin{align}
    Q_{k+1} (s_t, a_t) = R(s_t, a_t) + \gamma \mathbb{E}_{s_{t+1} \sim p_s} \big[ V_k (s_{t+1}) \big], \label{eq:softQ}
\end{align}
where $\gamma \in [0, 1]$ is the \textit{discounting factor}, $p_s$ is the state marginal of the trajectory distribution induced by $\pi$, and $V: \mathcal{S} \rightarrow \mathbb{R}$ is the soft \textit{state-value function}, i.e., the expected return starting from state $s$ and following policy $\pi$ thereafter. The state-value function $V$ is given by
\begin{align}
    V_k (s_t) = \mathbb{E}_{a_t \sim \pi} \big[ Q_k (s_t, a_t) - \alpha \log \pi (a_t | s_t) \big]. \label{eq:softV}
\end{align}
In the policy improvement step, the policy is updated towards the exponential of the soft Q-function. In practice, it is preferable to have tractable policies, so the policy is restricted to a set of policies $\Pi$, which can be, e.g., a family of parameterized distributions. In the update, the new policy must therefore be projected onto the set $\Pi$. \cite{haarnoja2018softapp} use information projection, and so the new policy is computed, for all states $s \in \mathcal{S}$, according to
\begin{align}
    \pi_{new} = \arg \min_{\pi \in \Pi} D_{KL} \Bigg( \pi(\cdot | s_t)~ \Bigg|\Bigg| ~ \frac{\exp \big(\frac{1}{\alpha} Q^{\pi_{old}} (s_t, \cdot) \big)}{Z^{\pi_{old}} (s_t)} \Bigg), \label{eq:softPolicyImprov}
\end{align}
where $Z^{\pi_{old}} (s)$ is a partition function that normalizes the distribution and $D_{KL}$ is the Kullback-Leibler divergence. 

SPI is only applicable to discrete state and action spaces. To extend SPI to continuous spaces, \cite{haarnoja2018softapp} introduce function approximators for the soft Q-function $Q_\theta$ and policy $\pi_\phi$, and alternate between optimizing their parameterization via gradient descent (instead of performing their evaluations) and policy improvement steps, yielding the SAC algorithm. The SAC algorithm models the soft Q-function using a neural network. The policy $\pi_\phi$ is typically modeled as a Gaussian distribution, where the mean $\mu_\phi$ and standard deviation $\sigma_\phi$ vectors are given by a neural network. The Q-function is updated via gradient descent, by minimizing a loss function based on the Bellman equations: 
\begin{align}
    \mathcal{L}(\theta) = \frac{1}{|\mathcal{D}|} \sum_{(s_t, a_t, r_t, s_{t+1}) \in \mathcal{D}} \frac{1}{2} \big(Q_\theta (s_t, a_t) - y \big)^2, \label{eq:loss}
\end{align}
where $\mathcal{D}$ is a mini-batch of experience examples $(s_t, a_t, r_t, s_{t+1})$, and $y$ is the target of the Q-network, derived from combining equations \eqref{eq:softQ} and \eqref{eq:softV}:
\begin{align}
    y = r_t + \gamma \Big(Q_{\bar{\theta}} (s_{t+1}, \Tilde{a}_{t+1}) - \alpha \log \pi (\Tilde{a}_{t+1} | s_{t+1}) \Big). \label{eq:SACtarget}
\end{align}
Here, the next action $\Tilde{a}_{t+1}$ is sampled from the current policy $\Tilde{a}' \sim \pi_\phi(\cdot | s_{t+1})$. The update utilizes a \textit{target} Q-network parameterized by $\bar{\theta}$ to stabilise training. The target Q-network is obtained by Polyak averaging the Q-network weights with smoothing constant $\rho$ over the course of training as
\begin{align}
    \Bar{\theta} \leftarrow \rho \theta + (1 - \rho) \Bar{\theta}.
\end{align}
The policy update can be computed by minimizing the expected KL-divergence in equation \eqref{eq:softPolicyImprov} via gradient descent
\begin{align}
    J(\phi) = \mathbb{E}_{s_t \sim \mathcal{D}} \Big[ \mathbb{E}_{a_t \sim \pi_\phi (\cdot | s_t)} \big[ \alpha \log \pi(a_t | s_t) - Q_\theta (s_t, a_t)  \big]  \Big]. \label{eq:PolicyUpdate}
\end{align}
Notice that the expression has been multiplied by $\alpha$ and the constant partition function $Z$ is ignored since it does not affect the gradient. The performance $J(\phi)$ is an expectation over actions, which are dependent on the policy parameters $\phi$, and so it is not possible to get an estimate of the gradient based on equation \eqref{eq:PolicyUpdate} directly. To get an expression for the gradient of the performance that can be estimated with samples, \cite{haarnoja2018softapp} use the \textit{reparameterization trick}. The policy is reparameterized by the transformation
\begin{align}
    \hat{a} = f_\phi(\epsilon_t, s_t)
\end{align}
where $\epsilon$ is some noise sampled from a fixed distribution. The transformation depends on the policy distribution used. For example, \cite{haarnoja2018softapp} use a squashed Gaussian in practice to ensure that the action values are bounded, in which case the appropriate transformation is
\begin{align}
    f_\phi(\epsilon_t, s_t) = \tanh \big( \mu_\phi (s_t) + \sigma_\phi (s_t) \odot \epsilon_t \big), \epsilon_t \sim \mathcal{N} (0, 1).
\end{align} 
With the transformation, the performance is then rewritten as
\begin{align}
     J(\phi) = \mathbb{E}_{s_t \sim \mathcal{D}, \epsilon_t \sim \mathcal{N}} \Big[  \alpha \log \pi(f_\phi(\epsilon_t, s_t) | s_t) - Q_\theta (s_t, f_\phi(\epsilon_t, s_t)) \Big]. \label{eq:ReParamPolicyUpdate}
\end{align}
One can notice that the expectation is no longer dependent on the policy parameters and so the gradient can be moved into the expectation and approximated. The full SAC algorithm is presented in algorithm \ref{alg:SAC} in \ref{appendix:SAC}.

\subsection{Federated learning}
\label{sec:FL}

Federated learning is a framework for learning a shared global model on decentralized data across multiple nodes, without the nodes sharing their private data. Unlike typical distributed learning, federated learning makes no assumptions about the data distribution across nodes being \textit{independent and identically distributed} (IID), and so it can be applied to non-IID settings as well. Furthermore, federated learning can also handle unbalanced data, i.e., some nodes having significantly larger local data sets than others. These characteristics allow federated learning to take advantage of massive amounts of data spread out over a large, heterogeneous network, e.g., pictures taken and stored on mobile phones, to learn a global model that generalizes well, while never communicating the local data itself. This maintains a higher degree of privacy across nodes while simultaneously eliminating the need for a central data center capable of storing the entire global data set. 

Federated learning is well-suited for smart HVAC system controllers due to its ability to accommodate the unique data distribution of each building and maintain the privacy of potentially sensitive information such as occupancy behavior. By training a controller across multiple buildings at once, we indirectly gather experience data more efficiently compared to training on a single building, which helps overcome the data efficiency limitations of reinforcement learning algorithms. Additionally, the heterogeneous data distribution results in more diverse experience data from different buildings, enhancing the generalization capabilities of the controller agent.

\subsubsection{Federated Averaging}

The federated learning setting consists of two main components. Firstly, we have a set of $K$ nodes, referred to as \textit{clients}, which compute updates to a shared global model independently of each other by training on their local data. Secondly, we have a central server, which coordinates the clients and updates the global model. One round of communication between the server and clients consists of the server sending the current global model parameters to a fraction $C \in (0, 1]$ of clients, chosen at random, the chosen clients computing their local updates, and finally sending their respective locally updated parameters to the server for model aggregation.

The federated optimization algorithm presented by \cite{mcmahan2017communication} can be applied to any problem with a finite-sum objective of the form
\begin{align}
    \min_{w \in \mathbb{R}^d} f(w) \quad \textrm{where} \quad f(w) \equiv \frac{1}{N} \sum_{i=1}^N f_i(w).
\end{align}
When applying federated optimization to, e.g., an actor-critic algorithm, we are optimizing two different objectives, where $f(w)$ corresponds to both $\mathcal{L}(\theta)$ and $J(\phi)$. Assuming the global data set is partitioned over $K$ clients, where $\mathcal{P}_k$ denotes the set of indexes of data points at client $k$, with $n_k = | \mathcal{P}_k |$, the objective can be rewritten as
\begin{align}
    f(w) \equiv \sum_{i=1}^N \frac{n_k}{N} F_k(w) \quad \textrm{where} \quad F_k(w) = \frac{1}{n_k} \sum_{i \in \mathcal{P}_k} f_i(w).
\end{align}

\cite{mcmahan2017communication} focus on the application of federated optimization to deep learning models, which are typically trained using some variant of stochastic gradient descent (SGD) to optimize their objective, the loss function. Hence, they use a federated version of SGD, called \textit{FedSGD}, as a starting point for their developed federated optimization algorithm. For one round of FedSGD, with fixed learning rate $\eta$ and fraction $C = 1$, each client $k$ computes the average gradient on their local data $g_k = \nabla F_k(w_t)$, where $w_t$ is the current global model. The local gradients are then aggregated at the central server and used to update the model according to
\begin{align}
    w_{t+1} \leftarrow w_t - \eta \nabla f(w_t) \quad \textrm{where} \quad \nabla f(w_t) = \sum_{k=1}^K \frac{n_k}{N} g_k. \label{FedSGDUpdate}
\end{align}
An equivalent update to \eqref{FedSGDUpdate} can be performed by taking one step of gradient descent on each local model $w_{t+1}^k \leftarrow  w_t - \eta \nabla g_k, \forall k$, and then aggregating the local model parameters via the following weighted average
\begin{align}
    w_{t+1} \leftarrow \sum_{k=1}^K \frac{n_k}{N} w_{t+1}^k. \label{FedAVGUpdate}
\end{align}
Since the update \eqref{FedAVGUpdate} is just an average over the parameters of each local model, it is possible to perform multiple local steps of gradient descent $w^k \leftarrow  w^k - \eta \nabla F_k(w^k)$ before averaging in order to increase the amount of computation per communication round. This is the core of the \textit{Federated Averaging} (FedAvg) algorithm.

\subsubsection{FedOpt}

In FedAvg, the updated global model parameters $w_{t+1}$ are computed by averaging the updated local parameters $w_{t+1}^k$ according to equation \eqref{FedAVGUpdate}. Alternatively, this update can be performed by computing the ``pseudo-gradient'' $\Delta_{t+1}$, which is the average of differences between the local parameters and the current global model, $\Delta_{t+1}^k = w_{t+1}^k - w_t$,
and adding it to the current parameters according to 
\begin{align}
    w_{t+1} \leftarrow w_t + \Delta_{t+1} \quad \textrm{where} \quad \Delta_{t+1} = \sum_{k=1}^K \frac{n_k}{N} \Delta_{t+1}^k. \label{PseudoGradUpdate}
\end{align}
Through this formulation, the server update in FedAvg can be viewed as taking one gradient ascent step using the pseudo-gradient and a \textit{global learning rate} $\eta_g$ = 1. \cite{reddi2020adaptive} recognize the possibility of choosing other values of $\eta_g$. They also suggest the possible use of alternative server update rules based on the pseudo-gradient, as well as utilizing other optimizers than SGD on the client side. Combining these ideas, \cite{reddi2020adaptive} generalize FedAvg into a framework called \textit{FedOpt}, presented in algorithm \ref{alg:FedOpt}.

\newcommand{\alglinelabel}{%
  \addtocounter{ALC@line}{-1}
  \refstepcounter{ALC@line}
  \label
}

\begin{algorithm}
\caption{FedOpt} \label{alg:FedOpt}
    \begin{algorithmic}[1]
        \STATE Initialise global model $w_0$
        \FOR{each communication round $t = 0,1,...,T$}
            \STATE $m \leftarrow \max \{ C \cdot K, 1 \}$
            \STATE $S_t \leftarrow$ random set of $m$ clients 
            \STATE $w_{k, 0}^t = w_t, \forall k \in S_t$
            \FOR {each client $k \in S_t$ \textbf{in parallel}}
                \FOR{$u = 0, 1,..., U - 1$}
                    \STATE Compute estimate $g_{k, u}^t$ of $\nabla F_k (w_{k, u}^t)$
                    \STATE $w_{k, u+1}^t$ = \text{ClientOpt}$(w_{k, u}^t, g_{k, u}^t, \eta_l, t)$
                \ENDFOR
                \STATE $\Delta_{t}^k = w_{k, U}^t - w_t$
            \ENDFOR
        \STATE $n_{tot} = \sum_{k \in S_t} n_k$
        \STATE $\Delta_{t} = \sum_{k \in S_k} \frac{n_k}{n_{tot}} \Delta_{t}^k$
        \STATE $w_{t+1} =$ \text{ServerOpt} $(w_t, \Delta_{t}, \eta_g, t)$ \alglinelabel{line:ServerOpt}
        \ENDFOR
    \end{algorithmic}
\end{algorithm}

ClientOpt and ServerOpt in algorithm \ref{alg:FedOpt} refer to the optimizers used at the clients and server, respectively. Any gradient-based optimizer can be applied. The hyperparameter $\eta_l$ sets the \textit{local learning rate} at the clients. The hyperparameter $U$ determines how many local updates to perform in each communication round. \cite{reddi2020adaptive} also allow the optimizers to depend on the communication round $t$ to facilitate the potential use of learning rate schedulers. 

\textit{FedAvgM}, which stands for Federated Averaging with Server Momentum \citep{hsu2019measuring}, slightly modifies the FedAvg algorithm by adding a momentum term $v$. During a server update (line \ref{line:ServerOpt} in algorithm \ref{alg:FedOpt}), the momentum is updated according to
\begin{align}
    v_t \leftarrow \mu v_{t-1} + \eta_g \Delta_t,
\end{align}
where $\mu \in [0, 1)$ determines the level of momentum. The global weight parameters are then updated using the momentum as
\begin{align}
    w_{t+1} \leftarrow w_t + v_t. \label{eq:FedAvgMUpdate}
\end{align}

\textit{FedAdam} is an adaptation of the Adam optimizer \citep{kingma2014adam} to ServerOpt, presented by \cite{reddi2020adaptive}. FedAdam uses two momentum terms $m$ and $v$ in the server update. The first momentum $m$ is computed as the exponential moving average
\begin{align}
    m_t \leftarrow \beta_1 m_{t-1} + (1 - \beta_1) \Delta_t
\end{align}
and the second momentum as the squared exponential moving average
\begin{align}
    v_t \leftarrow \beta_2 v_{t-1} + (1 - \beta_2) \Delta_t^2,
\end{align}
where $\beta_1, \beta_2 \in [0, 1)$ are hyperparameters. The global model update is then computed according to
\begin{align}
    w_{t+1} \leftarrow w_t + \eta_g \frac{m_t}{\sqrt{v_t} + \epsilon}. \label{eq:FedAdamUpdate}
\end{align}
Here, $\epsilon > 0$ controls the \textit{degree of adaptivity}.

\subsubsection{Gradient masking} 

\textit{Gradient masking} \citep{tenison2022gradient} can improve the performance of FL algorithms in heterogeneous settings. The idea of gradient masking is to apply a soft mask to the server update, which assigns higher importance to the components of the pseudo-gradients which are in agreement with the dominant direction, thereby better capturing the invariances across clients. In \cite{hagstrom2023using} gradient masking was found to improve the stability of the learning process by reducing the randomness across different seeds. The importance is determined by the sign agreement across parameters over the client updates $\Delta_{t}^k$. \cite{tenison2022gradient} define the agreement score $A \in [0, 1]$, which is given by
\begin{align}
    A \equiv \Big| \frac{1}{K} \sum_{k = 1}^K \textrm{sign} (\Delta^k) \Big|.
\end{align}
The agreement score is then used compute the mask $\Tilde{m}_\tau$ element-wise according to 
\begin{align}
    [\Tilde{m}_\tau]_j = 1 \quad \textbf{if} \quad A_j \geq \tau \quad \textbf{else} \quad A_j,
\end{align}
where $\tau \in (0, 1]$ is a hyperparameter determining the desired level of agreement. The mask $\Tilde{m}_\tau$ is then applied to the final computed update in ServerOpt before addition to the current model parameters via the element-wise product. The updates of FedAvg \eqref{PseudoGradUpdate}, FedAvgM \eqref{eq:FedAvgMUpdate} and FedAdam \eqref{eq:FedAdamUpdate} with gradient masking are thus
\begin{align}
    \text{FedAvg:} \quad w_{t+1} &\leftarrow w_t + \Tilde{m}_\tau \odot \Delta_{t+1} \\
    \text{FedAvgM:} \quad w_{t+1} &\leftarrow w_t + \Tilde{m}_\tau \odot v_t \\
    \text{FedAdam:} \quad w_{t+1} &\leftarrow w_t + \eta_g \Tilde{m}_\tau \odot \frac{m_t}{\sqrt{v_t} + \epsilon}.
\end{align}

\section{Experiments}
\label{sec:Experiments}

\subsection{Simulation Environment} \label{sec:SimEnv}

In our experiments, we use the open-source building simulation and control framework Sinergym (v.2.0.0) \citep{2021sinergym}. Sinergym provides an interface for interacting with the building energy model simulation tool EnergyPlus via the OpenAI Gym API \citep{brockman2016openai}, a popular API for implementing and evaluating reinforcement learning algorithms. Sinergym provides a handful of different building environments as well as several weather profiles. We conduct our experiments on the available data center environment\footnote{The name of the environment file is \textit{2ZoneDataCenterHVAC\_wEconomizer.idf}.}. The data center has a total area of 491.3 $m^2$. It is split into two asymmetrical zones; the west and east zone, equipped with their own respective HVAC systems. The HVAC systems are composed of air economizers, evaporative coolers, a direct expansion cooling coil, a chilled water coil and a variable air volume fan. The heating and cooling setpoints of each zone are controllable, and one episode of simulation runs for one year.

\subsubsection{Markov Decision Process Formulation} \label{sec:MDPExperiment}

To apply reinforcement learning algorithms to the control of the HVAC systems, we must provide an MDP formulation of the building environment. We define a state space $\mathcal{S}$, an action space $\mathcal{A}$ and a reward function $R$. One environment step, or control action, is taken every 15 minutes within the simulation, leading to a total number of 35 040 steps for one simulation episode.

The agent observes a state vector $s \in \mathcal{S} \subset \mathbb{R}^{18}$ of $18$ features. The complete list of features is presented in table \ref{tab:obsVar} in \ref{appendix:ObsSpace}. The features consist of the factors that we aim to control, namely the temperature of the zones and indirectly the energy consumption of the IT equipment and HVAC system, as well as other factors that relate to the temperature in the zones, e.g., outside air temperature. We also include ``forecasted'' outside temperature and air relative humidity values. This allows the agent to anticipate large changes in temperature and potentially counteract them by pre-heating or pre-cooling the zones. How the forecasted values are observed is described further in \ref{appendix:weather}.

The control variables of the data center model are the heating and cooling setpoint temperatures of each zone, and so the action $a \in \mathcal{A} \subset \mathbb{R}^4$ taken by the agent is a vector of 4 features, which determines these setpoint temperatures. The action space is described in table \ref{tab:actVar}. The actions are bounded by a range of possible values, which also include ``bad'' values that can lead the temperature in the zones to lie outside the comfortable range of values. The notion of good values should instead be encoded into the reward function and learned by the agent, irrespective of the possible range of values of the HVAC equipment available, as argued by \cite{biemann2021experimental}.
\begin{table}[htb]
    \caption{Description of the action space.} \label{tab:actVar}
    \centering
    \begin{tabular}{l c c} \\ \toprule 
         Feature & Range & Unit  \\  \midrule
         West zone cooling setpoint & [15.0, 22.5] & $\degree \textrm{C}$ \\ 
         West zone heating setpoint & [22.5, 30.0] & $\degree \textrm{C}$ \\
         East zone cooling setpoint & [15.0, 22.5] & $\degree \textrm{C}$ \\  
         East zone heating setpoint & [22.5, 30.0] & $\degree \textrm{C}$ \\  
         \bottomrule
    \end{tabular}
\end{table}

The goal is to train an agent that minimizes the total energy consumption of the data center. At the same time, the temperature inside the building must remain within the target range. Hence, we need to encode information about the energy consumption and the \textit{thermal comfort} into the reward signal. The reward function defined by \cite{biemann2021experimental} does precisely this, and so, we use it in our MDP formulation. They define the following reward function 
\begin{align}
    R(s) = r_{west} + r_{east} - \lambda_p (P_{it} + P_{hvac}), \label{eq:RewExperiments}
\end{align}
where $r_i$ is computed based on the thermal comfort in zone $i$, and $P_{it}$ and $P_{hvac}$ are the power demands of the IT and HVAC equipment, respectively. The term $\lambda_p \geq 0$ is a scaling factor for the energy component of the reward. Given the observed temperature $T_i$ in zone $i$, the thermal comfort component is computed as 
\begin{align}
    r_i &= \exp \big( -\lambda_g(T_i - T_{tgt})^2 \big) \notag \\
    &- \lambda_t \big( \max (T_{min} - T_i, 0) + \max (T_i - T_{max}, 0) \big) \label{eq:RewComfort}, 
\end{align}
where $T_{tgt}$ is the desired target temperature, and $T_{min}$ and $T_{max}$ are the lower and upper bounds of the comfortable temperature range. Scalars $\lambda_g, \lambda_t \geq 0$ are hyperparameters that determine the shape of the reward function. The first term in equation \eqref{eq:RewComfort} gives the function a Gaussian shape, with the purpose of motivating the agent to stay close to the target temperature, providing a more robust reward than a simple trapezoidal reward function. The second term, the trapezoid penalty, is added to extend the function to yield negative rewards far away from the center, helping the agent to better distinguish moderately bad actions from very bad ones than it would with the zero rewards of a simple Gaussian.

\begin{figure}[htb]
    \centering
     \includegraphics[width=0.49\textwidth]{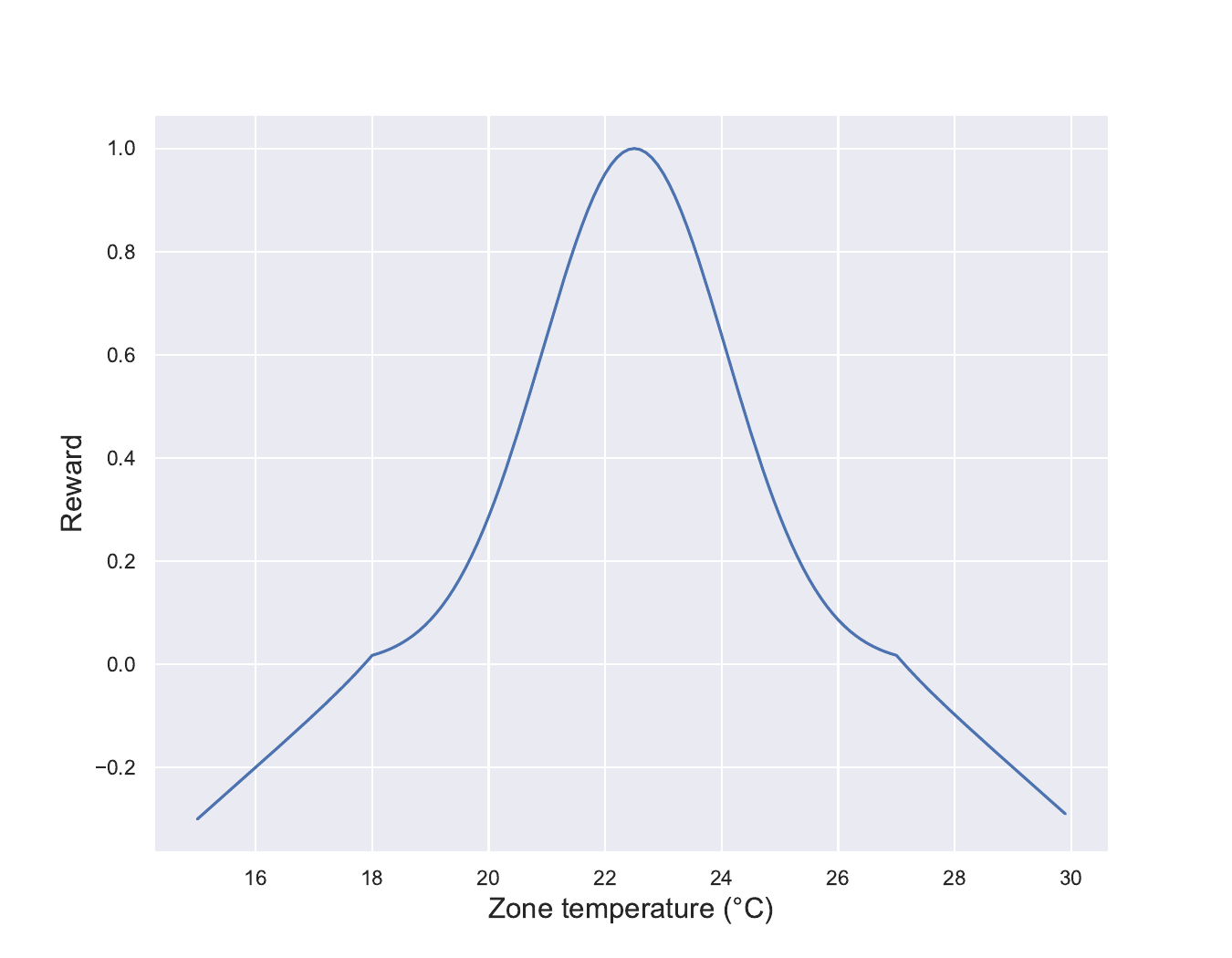}
     \caption{Graph of the zone thermal comfort reward $r_i$. The hyperparameters are set to $\lambda_g = 0.2$ and $\lambda_t = 0.1$, and the comfort range is bounded to $T_{min} = 18 \degree \textrm{C}$ and $T_{max} = 27 \degree \textrm{C}$. The target temperature is set to the midpoint of the comfort range, $T_{tgt} = 22.5 \degree \textrm{C}$.}
     \label{fig:RewComfort}
\end{figure}

The thermal comfort reward $r_i$ is close to $1$ when the temperature of zone $i$ is close to the target, and small or negative when close to or outside the comfort bounds. The total power demand $P_{tot} = P_{it} + P_{hvac}$ of the data center is in the order of 100 kW, and so to bring the energy penalty component in the reward function \eqref{eq:RewExperiments} to the same scale as the comfort component, we use the scaling factor $\lambda_p = 10^{-5}$ in our experiments. We set the comfort range bounds to $T_{min} = 18 \degree \textrm{C}$ and $T_{max} = 27 \degree \textrm{C}$ according to the recommended temperature range by the ASHRAE guidelines for data center power equipment \citep{tc2016data}. The target temperature is set to the midpoint of the comfort range $T_{tgt} = (T_{min} + T_{max}) / 2$, so as to motivate the agent to stay as far away from the edges of comfort as possible. Finally, we set the hyperparameters $\lambda_g = 0.2$ and $\lambda_t = 0.1$ as in \cite{biemann2021experimental}. Figure \ref{fig:RewComfort} displays the shape of the thermal comfort reward $r_i$ with the chosen parameters.

Sinergym provides 12 different weather profiles from significantly different climates. Each profile is fixed and provides hourly weather observations over a one-year period. The training of our agents spans multiple years, and so we do not wish to use the same weather profile for every year of training since we cannot know if the agent learns a useful policy for variable weather or if it simply overfits the weather profile. Thankfully, Sinergym allows us to add stochasticity to the weather from year to year. In \ref{appendix:weather} we provide further details and the full list of the weather profiles considered in table \ref{tab:weatherFiles}.

\subsection{Experiment configurations} \label{sec:Configurations}

We perform two main sets of experiments. In the first set, we train a federated HVAC control agent using FedAvg as the server optimizer. We evaluate the performance of three different client optimizers: SGD, SGD with momentum (SGDM), and Adam. We have 12 available weather profiles, and so we train on 11 client data centers, each with its own unique weather conditions. The Helsinki weather profile is reserved for evaluating the performance of the global agent in unseen environments. We consider two performance comparison baselines. The first is the employment of a proportional-integral-derivative (PID) controller, using temperature as its process variable and defining its error according to the setpoints described in figure \ref{fig:RewComfort} and with hyperparameters set as described by \cite{biemann2021experimental}. This choice is justified by its widespread use in HVAC control applications. We also train individual agents for each client and include their performance as a baseline. 

Lastly, in the second set of experiments, we evaluate two alternative federated algorithms, FedAvgM and FedAdam, using the best-performing client-side configuration from the first set.

Since our set of training clients is relatively small, we choose to include all clients in every global communication round, i.e., we set the fraction $C = 1$ for all our experiments. We also set the masking threshold to $\tau = 0.4$ in all experiments since it was found to generally perform well in \cite{tenison2022gradient} and \cite{hagstrom2023using}. We evaluate FedAvg, and so the global learning rate is set to $\eta_g = 1$. For the client optimizers, we only vary the learning rate, and use the default values for other hyperparameters. See table \ref{tab:cliopt:hyperparams} in \ref{appendix:ImpDetails} for the complete list of client optimizer hyperparameters. For the first set of experiments, we have two controllable hyperparameters, the client learning rate $\eta_l$ and the total of local updates per round $U$. For each client optimizer, we perform a search over the following grid of values
\begin{align*}
    &\eta_l \in \{ 0.0003, 0.001, 0.01, 0.1 \} \\
    &U \in \{4, 12, 24 \} . 
\end{align*}
For FedAvgM, the controllable hyperparameters are the global learning rate $\eta_g$, the number of local updates per round $U$, and the server momentum $\beta$. We perform a search over the following grid of values
\begin{align*}
    &\eta_g \in \{ 0.001, 0.01, 0.1, 1.0 \}  \\
    &U \in \{4, 12, 24 \}   \\
    &\beta \in \{0.8, 0.9, 0.99\}. 
\end{align*}
For FedAdam we set the degree of adaptivity to $\epsilon = 10^{-3}$, as \cite{reddi2020adaptive} find it to perform well across multiple different tasks. The controllable hyperparameters then are the global learning rate $\eta_g$, the number of local updates per round $U$, and the moment parameters $\beta_1$ and $\beta_2$. We perform a search over the following grid of values
\begin{align*}
    &\eta_g \in \{ 0.001, 0.01, 0.1, 1.0\} \\
    &U \in \{4, 12, 24 \}  \\
    &\beta_1 \in \{0.8, 0.9, 0.99\} \\
    &\beta_2 \in \{0.9, 0.99, 0.999\}. 
\end{align*}
In all experiments, each configuration is repeated 3 times with different random seeds to evaluate the robustness of each configuration. The training runs over a period of 15 years. The simulator takes a step in the environment, i.e., sends observations to the agent and executes the actions chosen by the agent, every 15 minutes, and so a full training run consists of a total of 525 600 environment interactions. For further implementation details, see \ref{appendix:ImpDetails}. The source code is available at https://github.com/hagstromf/FedHVAC.

\subsection{Results} \label{sec:Results}

In evaluating the performance of the models, we focus on the energy consumption and thermal comfort of the data center. The total energy consumption $E_{tot}$ is the cumulative total power consumption $P_{tot}$ over a year. The thermal comfort of the data center is evaluated in terms of thermal comfort violations. A thermal comfort violation takes place when the temperature in either or both zones of the building is outside the specified comfort range. The comfort violations 
are reported as the percentage of comfort-violating environment steps over a year.

\subsubsection{Evaluation results} \label{sec:Results:eval}

First, we consider the performance of FedAvg using different client optimizers. At the end of training, each model is run for three episodes on the Helsinki evaluation environment. The results are presented in table \ref{tab:cliopt:performance}, where the performance values are the means over the three evaluation episodes over all three random seed iterations. We report the values of the configuration that yielded the highest mean return. Further discussion on the performance of different hyperparameter configurations is provided in \ref{appendix:sensAnal}. 
\begin{table}
    \caption{Performance of the federated agent for different client optimizers on the evaluation environment (Helsinki) after 15 episodes of training. We choose the configuration that yields the highest return for reporting the performance of the federated agent, which are $\eta_l = 0.001, U = 24$ for Adam, $\eta_l = 0.1, U = 24$ for SGD, and $\eta_l = 0.1, U = 12$ for SGDM. The reported values are the means over three episodes of evaluation. $E_{tot}$ is the cumulative power consumption of the data center over one year, and Viol. is the comfort violation rate.} \label{tab:cliopt:performance}
    \centering
    \begin{tabular}{clrr} \\ \toprule  
        &&$E_{tot}$ (GWh)&Viol. (\%)\\ \midrule
        & Adam  &  $\mathbf{0.9189}$ & 0.0016 \\
        & SGDM  & 0.9220 & 0.0092 \\
        & SGD & 0.9266 & 0.0438 \\
        & PID-Baseline & 0.9311 & $\mathbf{0.0}$ \\
        \bottomrule
    \end{tabular}
\end{table} 

\begin{figure}[htb]
    \centering
    \begin{subfigure}[b]{0.49\textwidth}
         \centering
         \includegraphics[width=\textwidth]{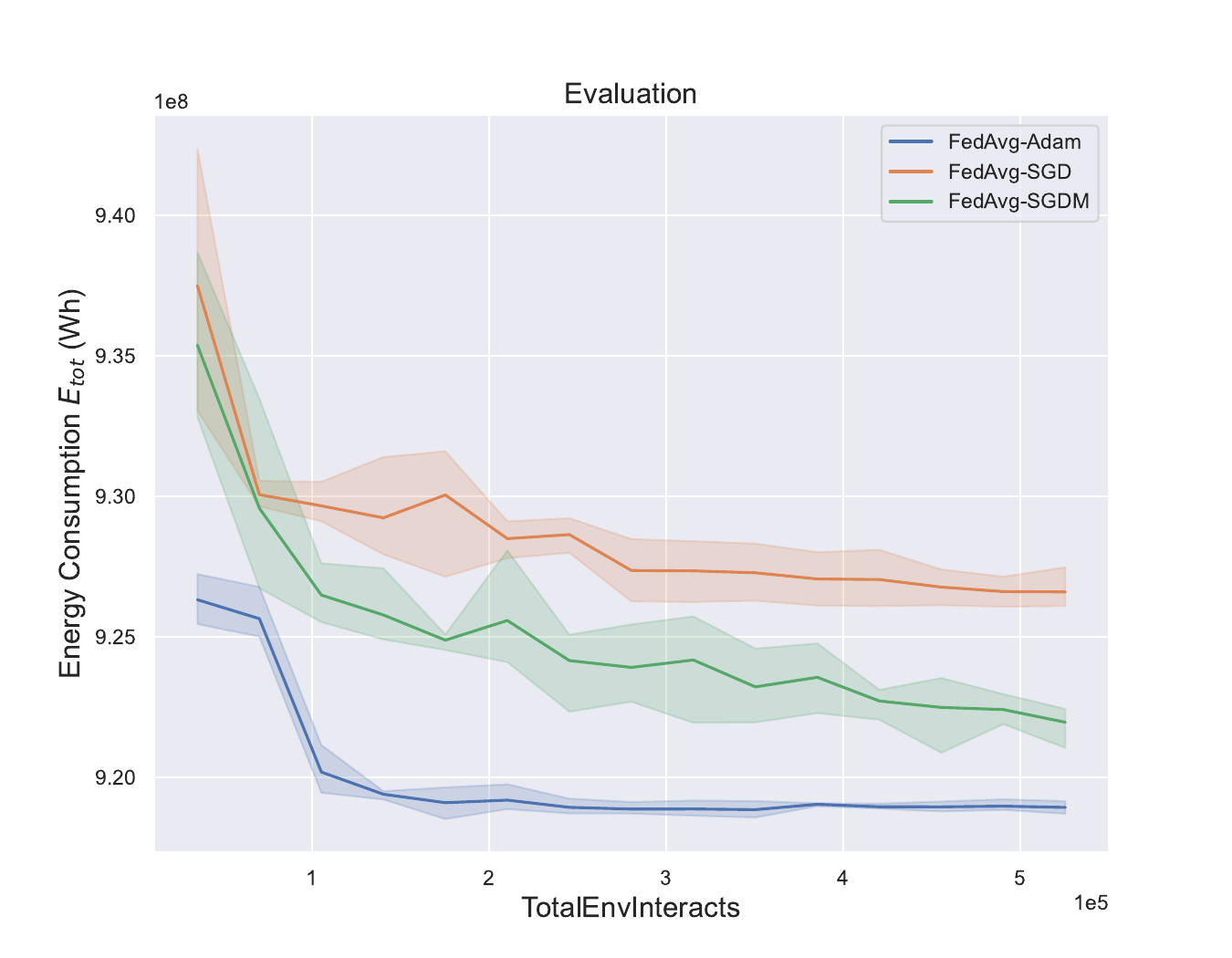}
         \caption{Energy consumption}
         \label{eval:fedavg:all:energy}
    \end{subfigure}
    \hfill
    \begin{subfigure}[b]{0.49\textwidth}
         \centering
         \includegraphics[width=\textwidth]{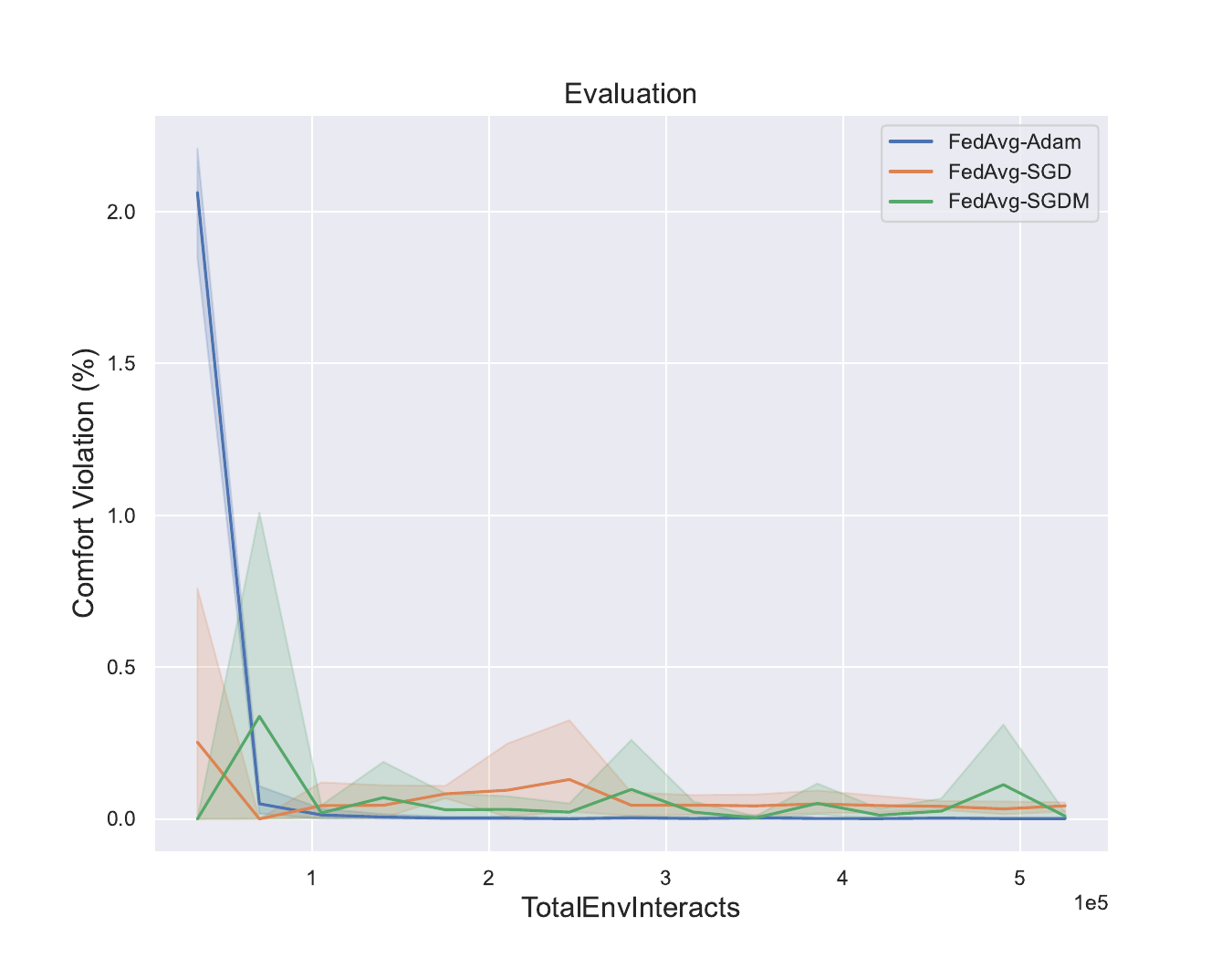}
         \caption{Comfort violation}
         \label{eval:fedavg:all:comfort}
    \end{subfigure}
    \caption{Progression of the energy consumption and comfort violation on the Helsinki evaluation environment of the FedAvg agent with different client optimizers.}
    \label{eval:fedavg:all}
\end{figure}

From table \ref{tab:cliopt:performance}, we notice that FedAvg with Adam outperforms SGD and SGDM in terms of both energy consumption and comfort violations when deployed on an unseen environment, indicating the best generalisation capabilities of the three. In figure \ref{eval:fedavg:all}, we show the progression of the energy consumption and comfort violation of the FedAvg agents on the evaluation environment for all client optimizers. The agents are evaluated for three episodes at the end of each episode of training. We plot the mean values over the three episodes over all random seeds with their bootstrapped 95 \% confidence intervals. Based on the progression plots, we notice that FedAvg with Adam does not only offer improved generalisation compared to the others. It also displays faster learning speeds, with the energy consumption converging after about five episodes and the comfort violations converging after just three episodes, while SGD converges in roughly eight episodes, and SGDM has not converged yet at the end of training. FedAvg with Adam has better learning stability as well. The tighter confidence intervals regarding both energy consumption and comfort violation indicate that the learning is more robust to randomness in the model initialization and training process, and thus its performance is more reliable. Lastly, it can be noticed that all federated learning agents overperform the PID controller in terms of average power consumption. No comfort violation is observed for the PID controller. 

\begin{figure}[htb]
    \centering
    \begin{subfigure}[b]{0.49\textwidth}
         \centering
         \includegraphics[width=\textwidth]{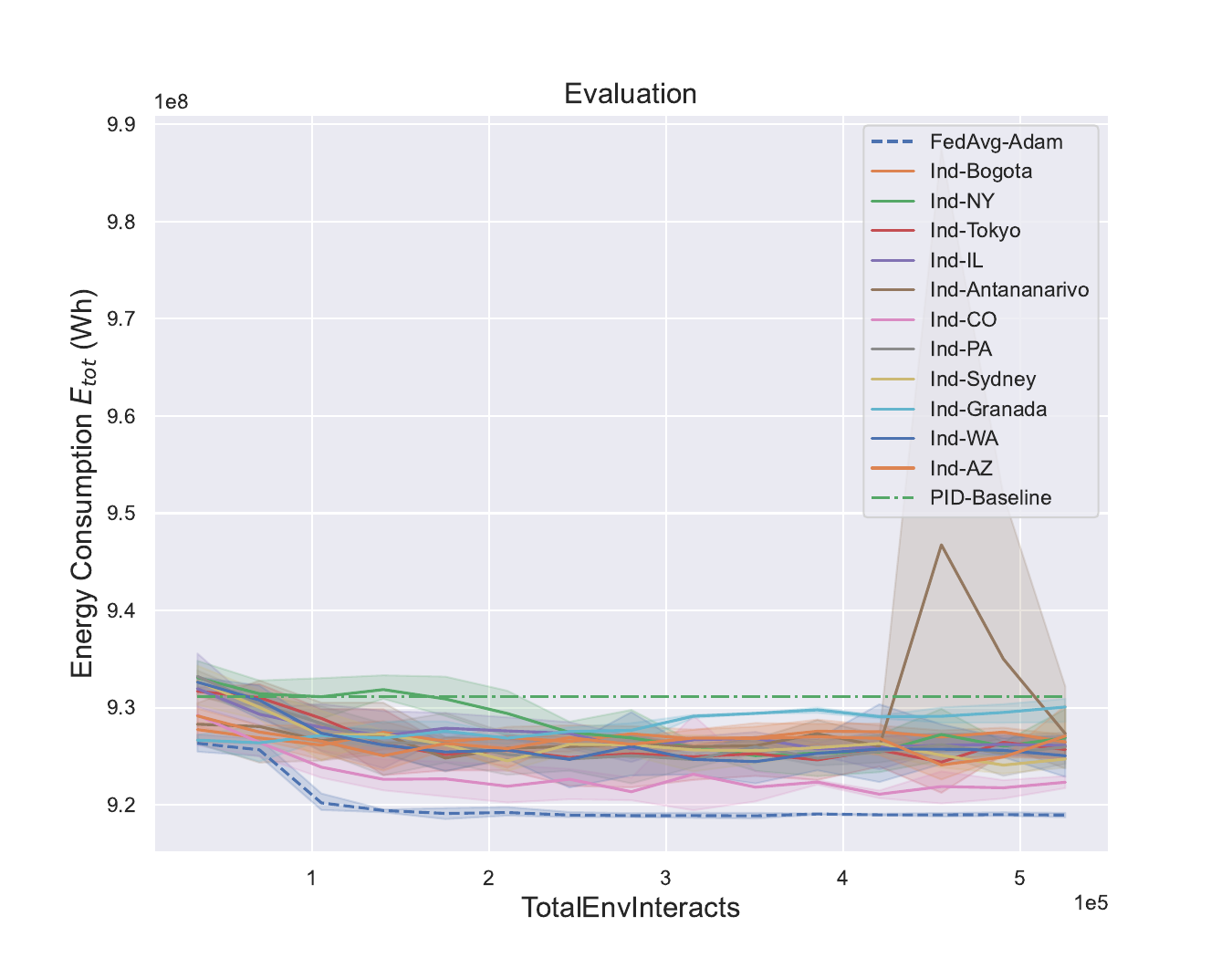}
         \caption{Energy consumption}
         \label{eval:fedavg:ind:energy}
    \end{subfigure}
    \hfill
    \begin{subfigure}[b]{0.49\textwidth}
         \centering
         \includegraphics[width=\textwidth]{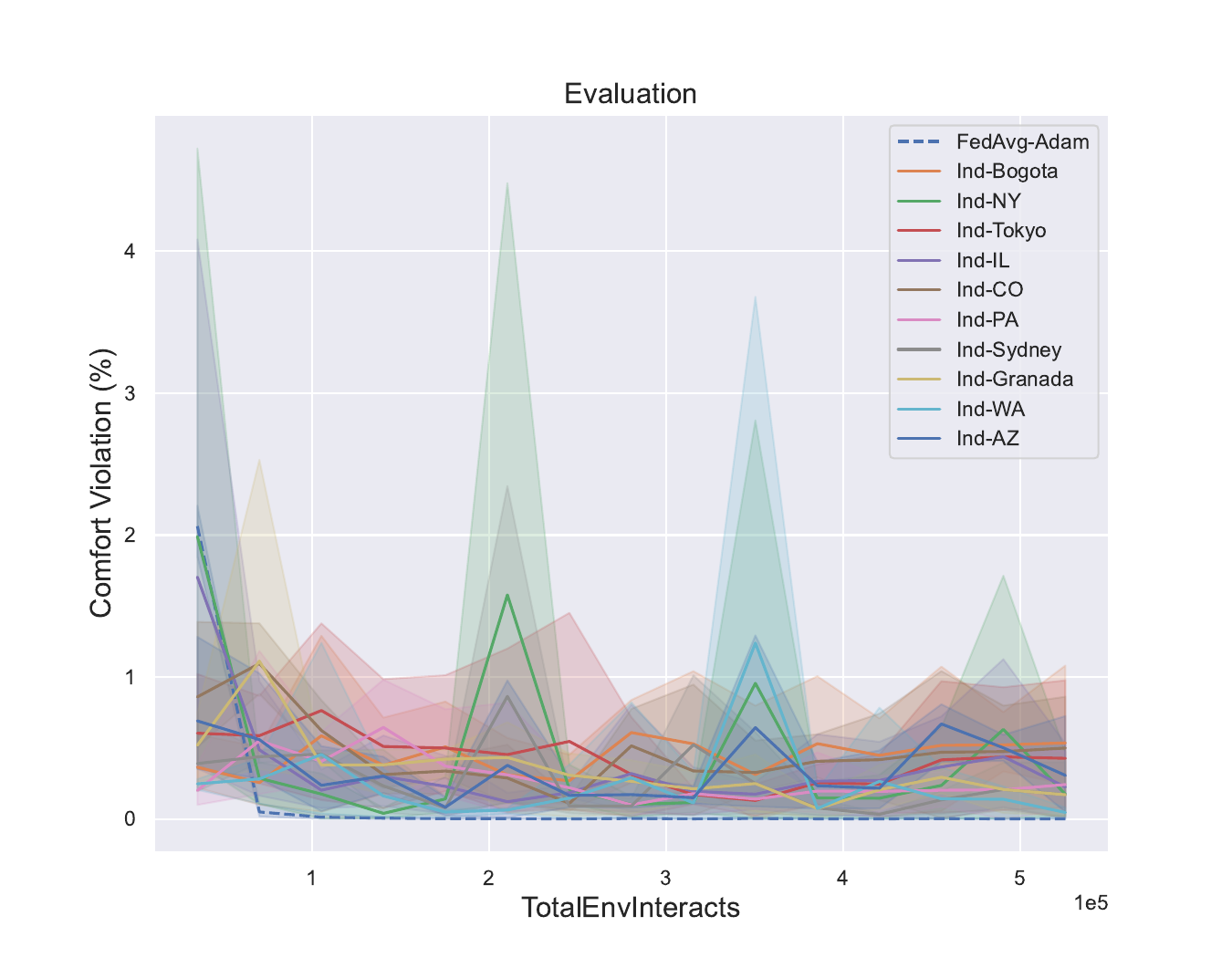}
         \caption{Comfort violation}
         \label{eval:fedavg:ind:comfort}
    \end{subfigure}
    \caption{Progression of the energy consumption and comfort violation on the Helsinki evaluation environment of FedAvg and independent agents with Adam as client optimizer. In the comfort violation plot \ref{eval:fedavg:ind:comfort} we omit the outlier Antananarivo for the sake of legibility.}
    \label{eval:fedavg:ind}
\end{figure}

Next, we analyze how the evaluation performance of a federated agent compares to agents trained independently on the clients. We focus on the best performing federated agent, i.e., with Adam as the client optimizer and $\eta_l = 0.001, U=24$, and compare it to the best performing individual agents, with Adam and $\eta_l = 0.01$. The progression of energy consumption and comfort violations of the federated and independent agents in the evaluation environment are presented in figure \ref{eval:fedavg:ind}. From figure \ref{eval:fedavg:ind:energy}, we see that the federated agent outperforms every independent agent in terms of energy consumption, converging to a lower value, and at a faster rate. We also notice a high variance in the energy consumption, both across different clients as well as across different runs for each client, with a significant outlier in the agent trained in the Antananarivo environment. Remarkably, the variance for the federated agent is significantly lower. Similar observations are made regarding the comfort violation in figure \ref{eval:fedavg:ind:comfort}, though we note that the independent agents tend to outperform the federated agent in the first episode. 

These observations support our conclusion that using federated optimization to train an HVAC control agent can significantly improve generalization, with a better performance in an unseen environment than any independently trained agent. Federated training can also improve the learning speed, generally converging faster, as well as learning stability, displaying a significant reduction in the variance in performance over different random seeds. In any real-world application, this consistency is a highly desirable trait, since we are not able to train the agent multiple times and therefore need a model that can reliably learn a good policy despite the inherent randomness of the real environment and training.

\subsubsection{Training results} \label{sec:TrainPerf}
\begin{figure}[htb]
    \centering
    \begin{subfigure}[b]{0.49\textwidth}
         \centering
         \includegraphics[width=\textwidth]{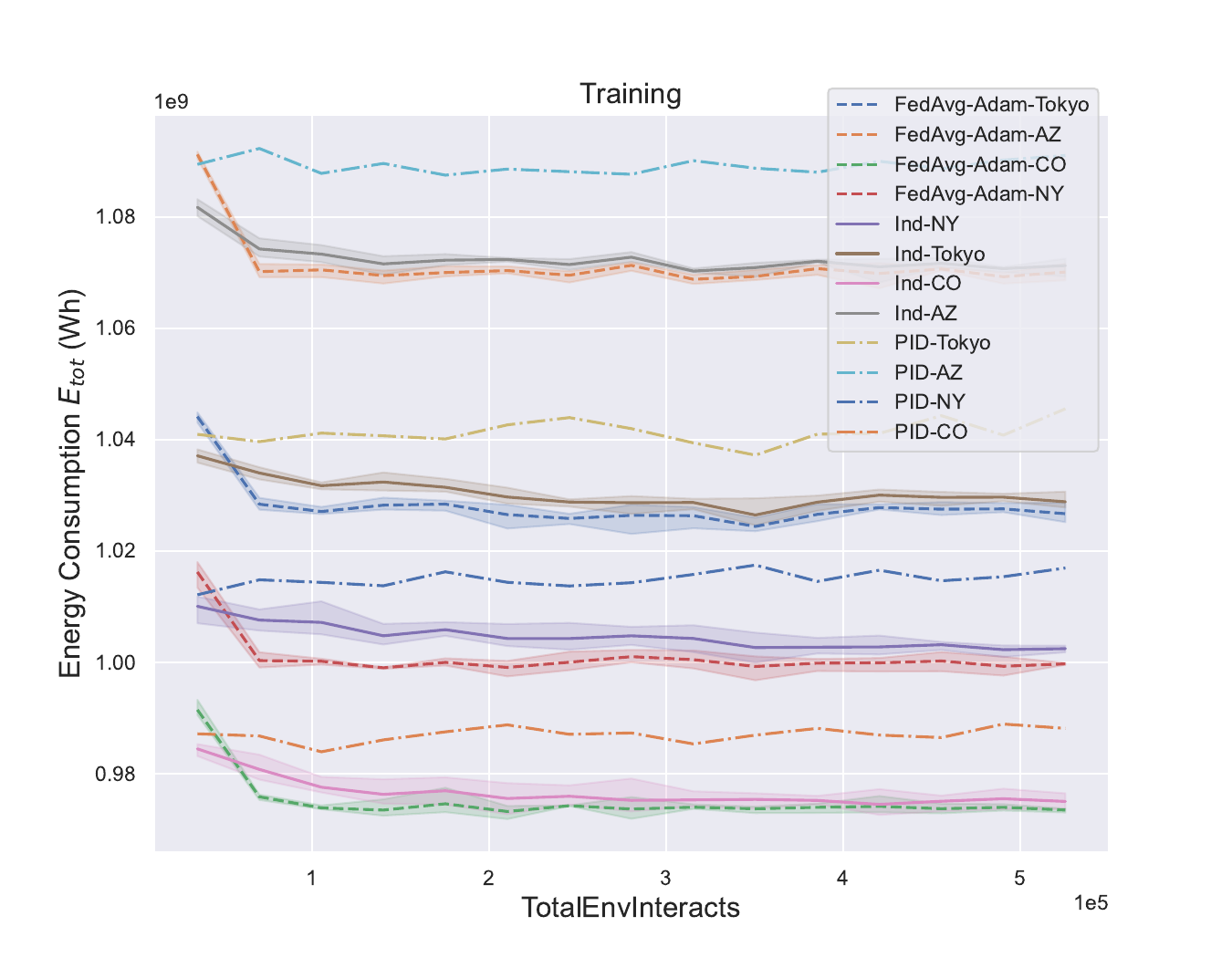}
         \caption{Energy consumption}
         \label{train:fedavg:ind:energy}
    \end{subfigure}
    \hfill
    \begin{subfigure}[b]{0.49\textwidth}
         \centering
         \includegraphics[width=\textwidth]{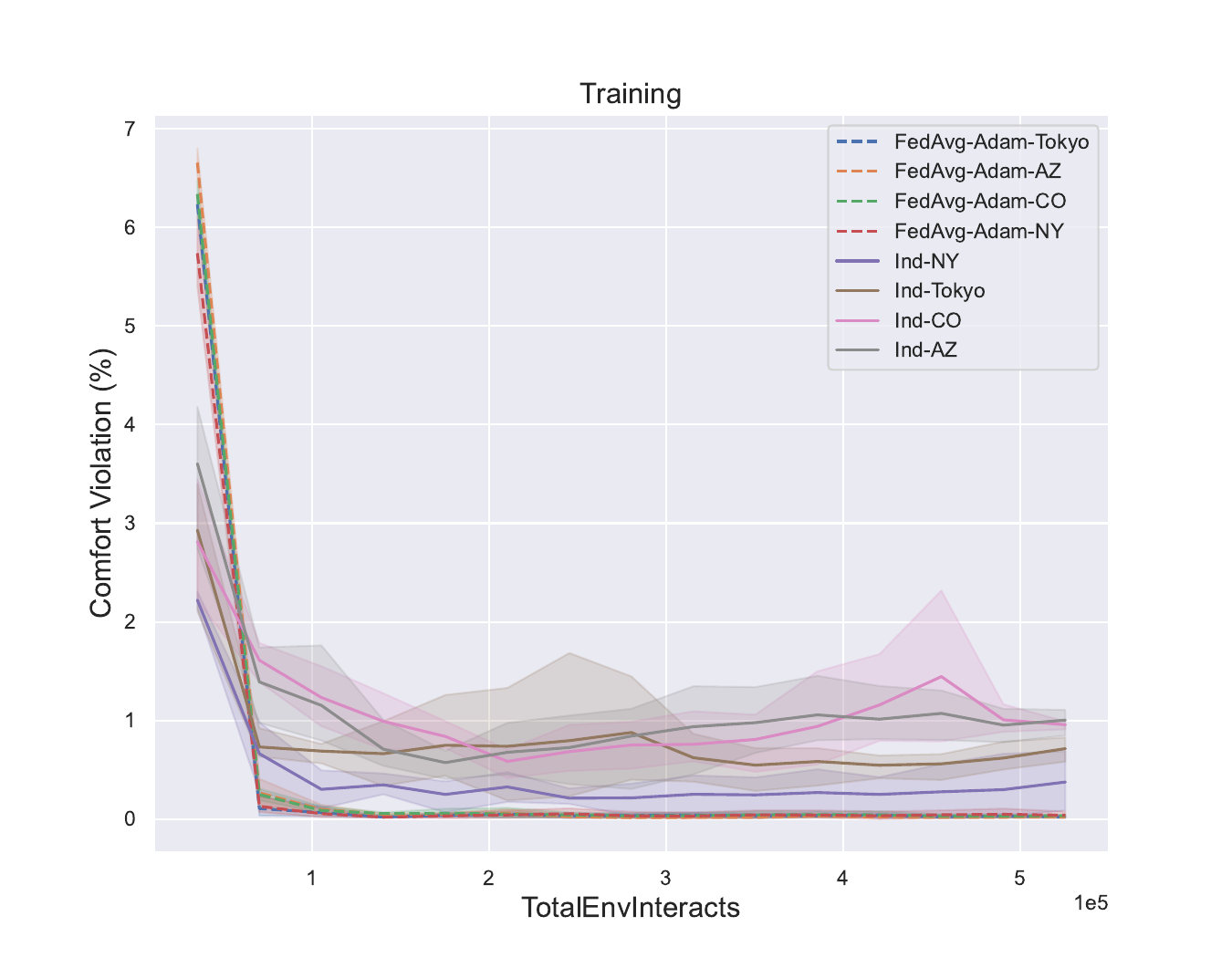}
         \caption{Comfort violation}
         \label{train:fedavg:ind:comfort}
    \end{subfigure}
    \caption{Progression of the energy consumption and comfort violation of FedAvg and independent agents on training environments Tokyo, AZ, CO and NY.}
    \label{train:fedavg:ind}
\end{figure}

We have seen that applying federated optimization can improve the performance of a reinforcement learning HVAC control agent in an unseen environment. Thankfully, this does not come at the expense of poorer performance in the training environments. In figure \ref{train:fedavg:ind}, we present the evolution of energy consumption and comfort violation of the federated agent and independent agents in the training environments. We only plot a subset of the environments for the sake of legibility. The behaviour of the omitted environments is consistent with the ones shown and analysed in this section. For additional figures of the remaining environments, please refer to \ref{appendix:AddPlots}.

The energy consumption in figure \ref{train:fedavg:ind:energy} displays similar improvements from using FedAvg in the evaluation environment. FedAvg generally converges faster and manages to reach a lower level of energy consumption. We also see improved learning stability, with slightly less variance across training runs. These improvements are even more pronounced when analysing the progression of comfort violation in figure \ref{train:fedavg:ind:comfort}. While the federated agent converges to near-zero comfort violations after two or three episodes, the independent agents never achieve near-zero violations. They also exhibit significantly more variance across both agents and different training runs. This shows that the model does not only benefit from the improved generalization, learning speed and learning stability of federated learning when applied to an unseen environment but also during training itself.

While the federated agent outperforms the independent agents in the long run, we notice that the independent agents tend to perform better during the first episode, both in terms of energy consumption and comfort violation. This, however, seems to be an effect of the larger client learning rate $\eta_l$ used for the independent agents. All federated agents perform better than the PID controller regarding energy consumption throughout the episodes.

\begin{figure}[htb]
    \centering
    \begin{subfigure}[b]{0.49\textwidth}
         \centering
         \includegraphics[width=\textwidth]{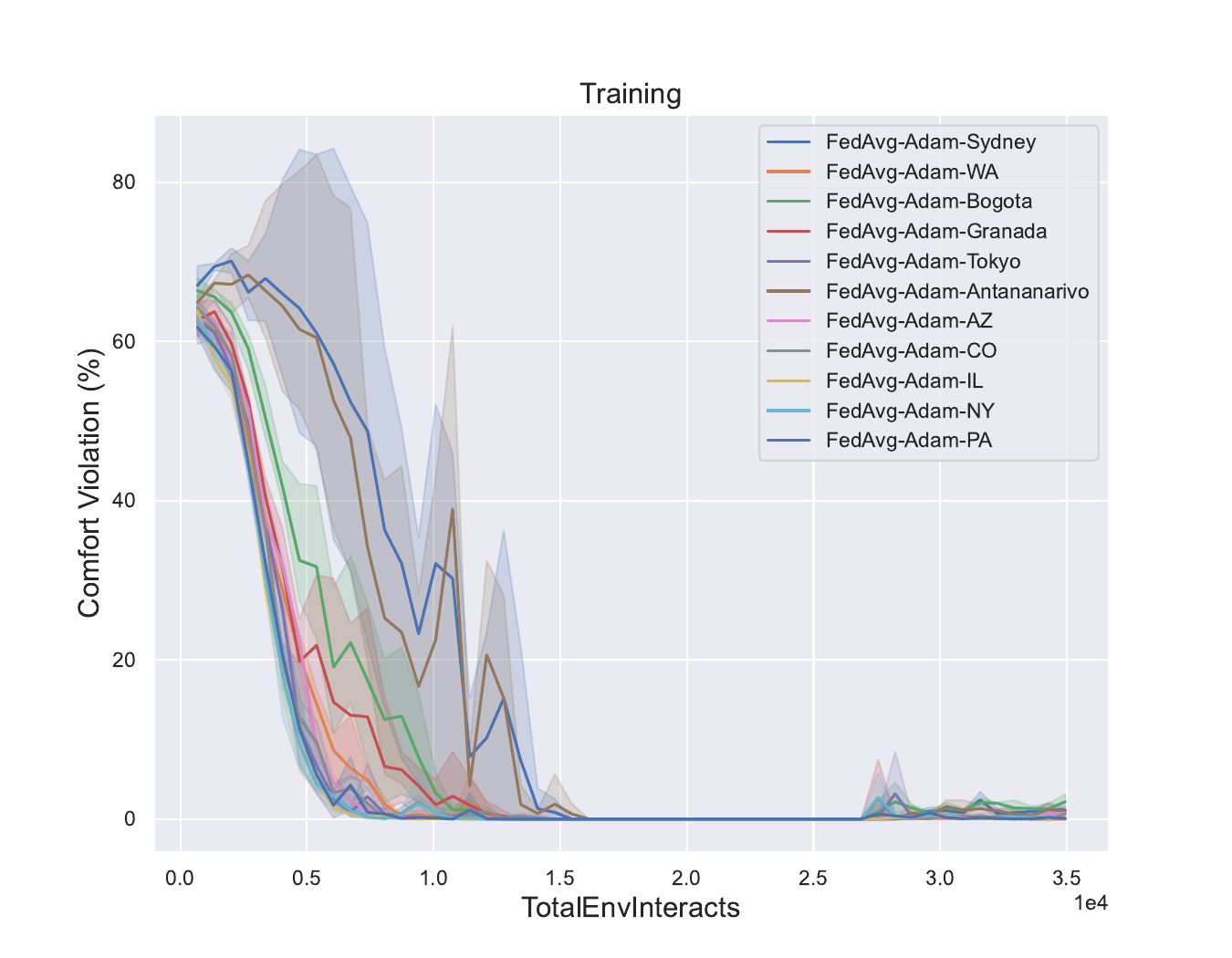}
         \caption{$\eta_l = 0.001$}
         \label{train:fedavg:weekly:comfort:0.001}
    \end{subfigure}
    \hfill
    \begin{subfigure}[b]{0.49\textwidth}
         \centering
         \includegraphics[width=\textwidth]{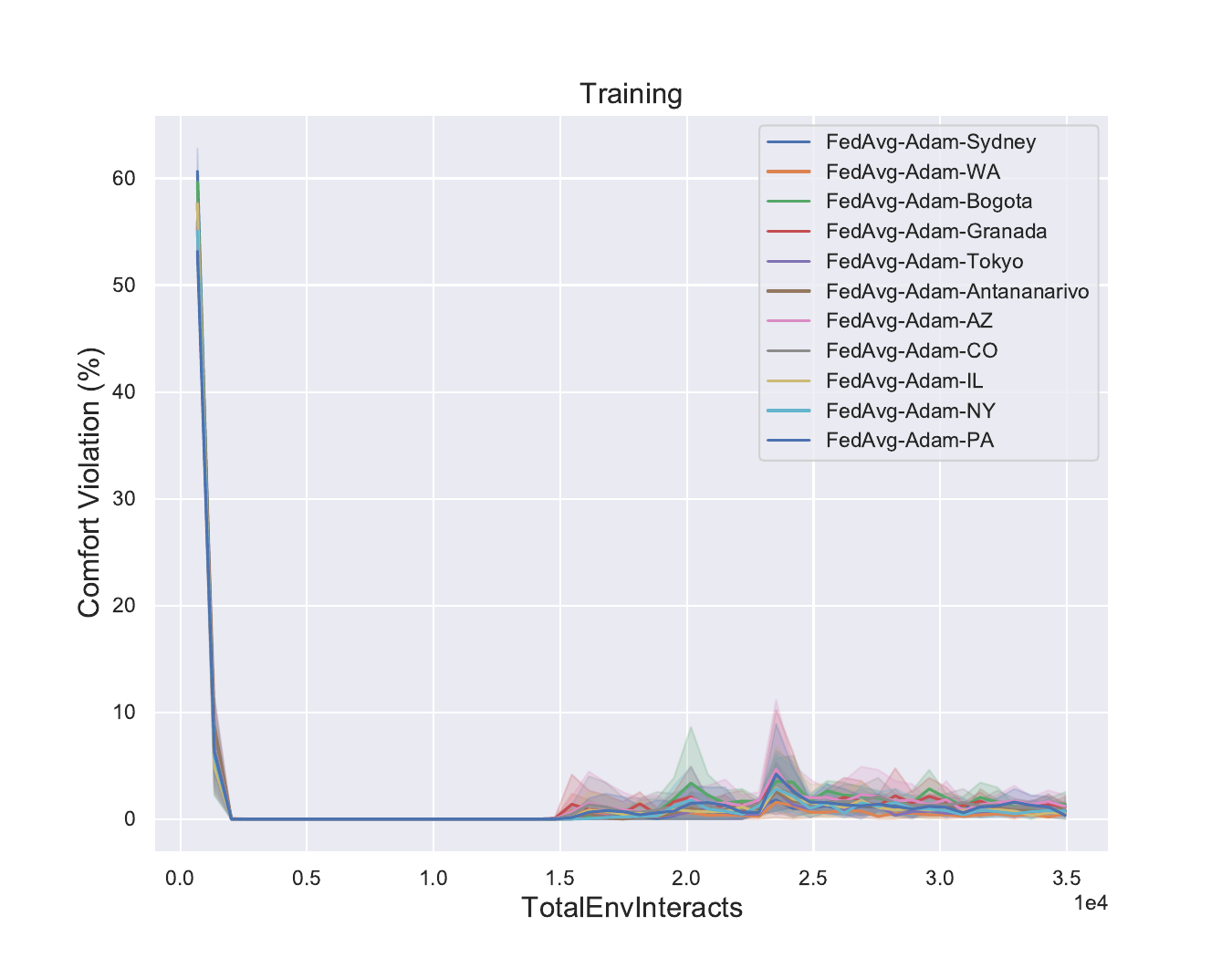}
         \caption{$\eta_l = 0.01$}
         \label{train:fedavg:weekly:comfort:0.01}
    \end{subfigure}
    \caption{Progression of the energy consumption and comfort violation of FedAvg and independent agents on training environments Tokyo, AZ, CO and NY.}
    \label{train:fedavg:weekly:comfort}
\end{figure}

In figure \ref{train:fedavg:weekly:comfort} we present the weekly comfort violations of the federated agent on the training environments over the first year of training for client learning rates $\eta_l = 0.001$ and $\eta_l = 0.01$. In this setting, the federated agent requires less than a full year of training to reach near-zero comfort violation. Depending on the environment, the agent with the lower client learning rate $\eta_l = 0.001$ requires between around 8000 to 17000 steps to reach near-zero violations, corresponding to roughly 12 to 25 weeks. Some of the environments experience a small increase towards the end of the year. However, if we increase the client learning rate to $\eta_l = 0.01$, we can achieve near-zero comfort violation significantly faster, in just three weeks, though there is an increase in violations for the second half of the year. While increasing the client learning rate can lead to significantly faster comfort violation reduction, it comes at the cost of significantly worse performance in the long run (see figure \ref{sens:Adam:clr} in \ref{appendix:sensAnal:Client}). On the other hand, while training with a lower client learning rate leads to great performance, the violations during the first few months are inadmissible, and so this federated agent would not be suitable for a real-world setting.  

\subsubsection{Server optimizers}
\begin{figure}[htb]
    \centering
    \begin{subfigure}[b]{0.49\textwidth}
         \centering
         \includegraphics[width=\textwidth]{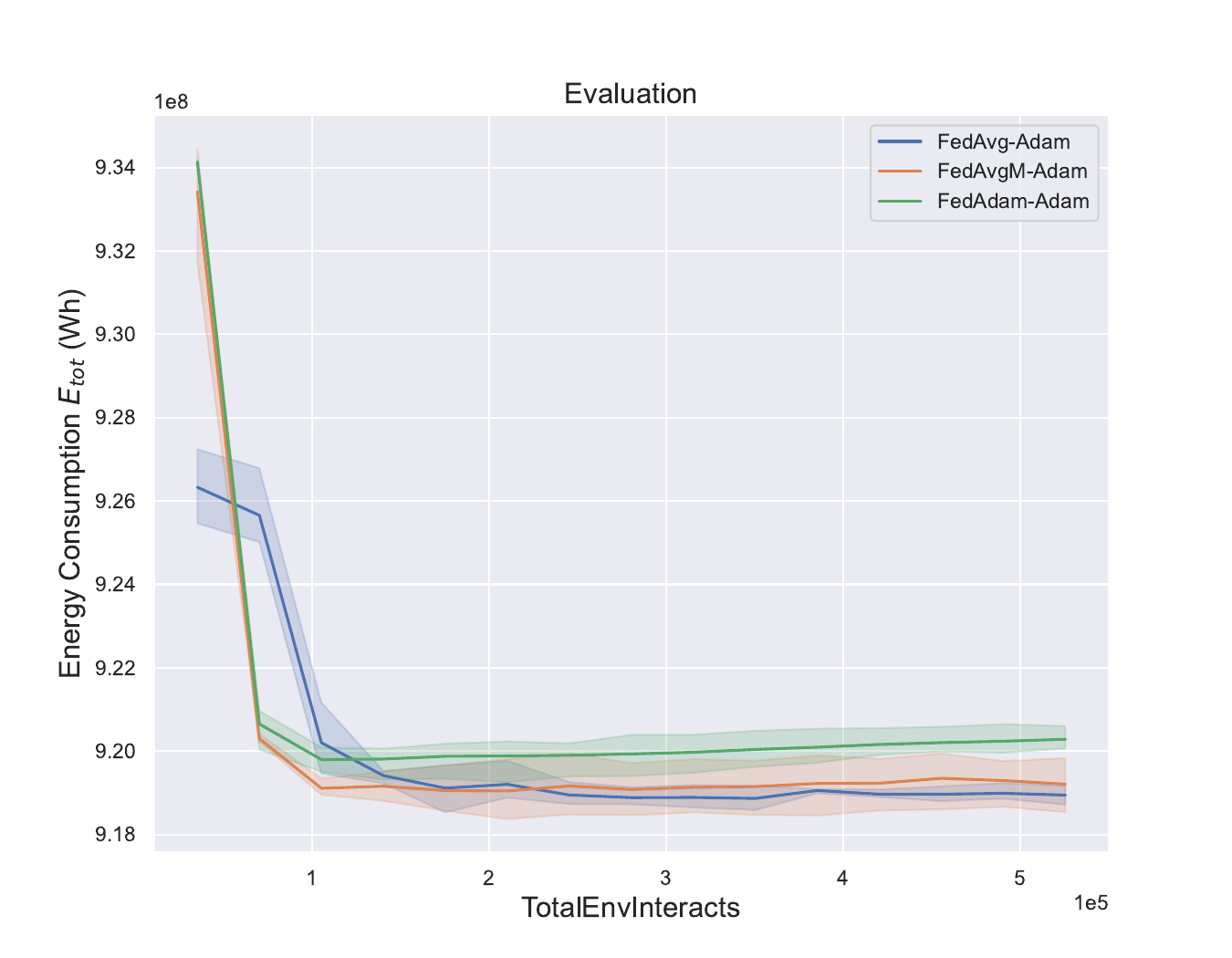}
         \caption{Energy consumption}
         \label{eval:server:all:energy}
    \end{subfigure}
    \hfill
    \begin{subfigure}[b]{0.49\textwidth}
         \centering
         \includegraphics[width=\textwidth]{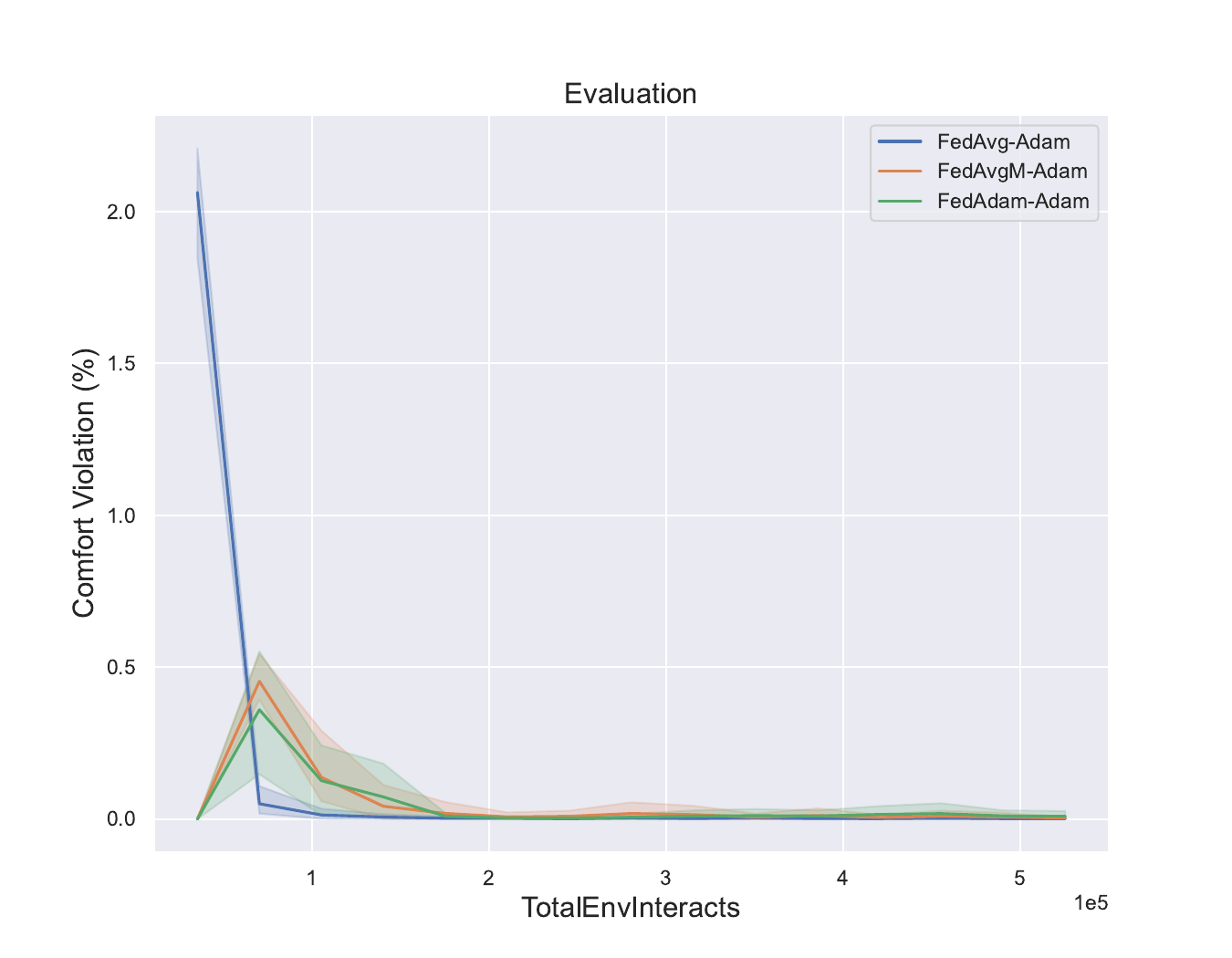}
         \caption{Comfort violation}
         \label{eval:server:all:comfort}
    \end{subfigure}
    \caption{Progression of the energy consumption and comfort violation on the Helsinki evaluation environment of FedAvg, FedAvgM and FedAdam.}
    \label{eval:server:all}
\end{figure}

Besides FedAvg, we also evaluate two alternative server optimizers: FedAvgM and FedAdam. Both use Adam as the client optimizer, with $\eta_l = 0.001$. The evaluation performances of the best-performing configurations of each optimizer at the end of training are presented in table \ref{tab:servopt:performance}. The progression plots of the energy consumption and comfort violations on the evaluation environment are presented in figure \ref{eval:server:all}. From table \ref{tab:servopt:performance}, we see that, although all perform similarly, FedAvg slightly outperforms the others in terms of both energy consumption and comfort violation. Looking at figures \ref{eval:server:all:energy} and \ref{eval:server:all:comfort}, the most striking difference is the early performance of the optimizers. FedAvg has significantly worse comfort violations than FedAvgM and FedAdam in the first episode, but performs better in the second. The opposite is true for the energy consumption. We also notice that FedAvg has tighter confidence intervals, and so offers better learning stability than FedAvgM and FedAdam.
\begin{table}
    \caption{Performance of the federated agent for different server optimizers on the evaluation environment (Helsinki) after 15 episodes of training. We choose the configuration that yields the highest return for reporting the performance of the federated agent, which are $\eta_g = 0.1$, $U = 24$, $\mu = 0.9$ for FedAvgM, and $\eta_g = 0.001$, $U = 24$, $\beta_1 = 0.8$, $\beta_2 = 0.9$ for FedAdam. The reported values are the means over three episodes of evaluation. $E_{tot}$ is the cumulative power consumption of the data center over one year, and Viol. is the comfort violation rate.} \label{tab:servopt:performance}
    \centering
    \begin{tabular}{clrr} \\ \toprule  
        &&$E_{tot}$ (GWh)&Viol. (\%)\\ \midrule
        & FedAvg  &  $\mathbf{0.9189}$ & $\mathbf{0.0016}$ \\
        & FedAvgM  & 0.9192 & 0.0035 \\
        & FedAdam & 0.9203 & 0.0092 \\
        \bottomrule
    \end{tabular}
\end{table} 
%


\section{Discussion}
\label{sec:Discussion}
Through our experiments, we have identified three key improvements from applying federated optimization to training reinforcement learning HVAC controllers. Firstly, by learning from experience collected from multiple heterogeneous environments, the agent gains access, albeit indirectly, to a larger amount of training data, which generally encompasses more variability than that available to any independent agent. In other words, there is an increase in exploration, which leads to a more informed global agent that can \textit{generalize better} to different environments. Secondly, the amount of total experience increases at a faster rate, which leads to an \textit{increase in learning speed}. Finally, when aggregating over the local agents, the dominant direction in the pseudo-gradient will have the most impact on the global update. This seems to have a regularizing effect, making it more difficult for the agent to branch off into sub-optimal regions of the policy space, \textit{increasing the learning stability}.

Our experiments show that the choice of client optimizer can have a significant impact on performance. The federated model can benefit from adaptivity on the local optimizer, as we found Adam to perform considerably better than both SGD and SGDM in terms of generalization, learning speed and learning stability. Meanwhile, the choice of server optimizer seems less critical. 

Regarding the server optimizer, both FedAvgM and FedAdam display comparable performance to FedAvg, although FedAvg slightly outperforms them. Moreover, FedAvg also has the advantage of having fewer hyperparameters to be tuned. As shown in \ref{appendix:sensAnal:Server}, these algorithms can be considerably sensitive to the choice of said hyperparameters. In a real-world scenario, we cannot evaluate multiple different hyperparameters and, as such, it is desirable to use an algorithm with as few adjustable parameters as possible, with minimal sensitivity to said hyperparameters. Considering that FedAvg has fewer adjustable components, combined with the observation that neither FedAvgM nor FedAdam seems to offer any significant improvement in terms of either generalization, learning speed or learning stability, we believe that FedAvg provides a more defendable choice for future efforts related to deployment in real-world settings.

A few limitations of our experiments are worth highlighting. The training of the federated agent was revealed to be considerably sensitive to the clients' learning rate. With lower learning rates, federated optimization offers stable and fast learning, but it is not suitable as-is for a real-world building environment due to the high degree of comfort violations at the beginning of training. By increasing the learning rate, it is possible to significantly reduce the comfort violations early on, but this comes at a trade-off for significantly worse final performance. Second, although significant progress can be achieved by using federated learning in this particular context, challenges remain in bridging the gap between simulation and real-world deployment, which can be noted from the time taken, between 3 and 24 weeks depending on the hyperparameter configuration, for the HVAC control agent to reach satisfactory energy consumption and comfort violation performance.

From a practical perspective in the context of HVAC control, federated learning offers significant benefits, being this a setting in which the underlying tasks across different buildings are largely similar yet subject to local variations. By allowing individual controllers to learn from their own operational data while sharing only aggregated model updates, the federated approach leverages commonalities across similar systems while preserving the confidentiality of sensitive information, such as occupancy patterns.

In our experiments, the primary focus was on energy consumption and comfort violations. However, practical real-world deployments would need to additionally account for communication or computational overhead. Nonetheless, the federated learning framework inherently reduces communication requirements by transmitting only aggregated model updates instead of raw data, while distributing the computational load across local nodes. This design suggests that, in a real-world HVAC system, the overhead from model aggregation is likely to be modest compared to the substantial benefits in terms of generalization, learning speed, and associated data privacy benefits.

\section{Conclusion}
\label{sec:Conclusion}
In this paper, we have experimentally evaluated the effects of training reinforcement learning HVAC control agents via federated optimization. We have trained Soft Actor-Critic (SAC) agents using Federated Averaging (FedAvg) with gradient masking, evaluating and comparing the performance of three different client optimizers: stochastic gradient descent (SGD), stochastic gradient descent with momentum (SGDM), and Adam. We have also compared the performance of federated agents to that of individual agents, trained on each respective client environment used in the federated learning scenario, both in terms of their performance in an unseen test environment and their performance in the training environments themselves. Furthermore, two alternative server optimizers, Federated Averaging with server momentum (FedAvgM) and FedAdam were compared to the FedAvg algorithm.

Our results have demonstrated that federated learning can improve generalization and the learning speed and stability of reinforcement learning-based HVAC controllers, which are critical bottlenecks for their adoption in real-world settings. However, there are still important challenges that must be addressed in that direction, mainly related to the time required for the learning-based controllers to learn policies that perform satisfactorily. 

Moreover, while our numerical experiments demonstrate clear benefits in terms of learning speed, generalization, and stability when employing federated optimization for reinforcement learning-based HVAC control, we acknowledge that these outcomes constitute only a first, albeit critical,  step towards their wider deployment. As such, real-world pilot deployments remain essential to conclusively verify practical benefits and applicability in realistic building settings and thus warrant further research efforts.

Future research could be dedicated to bridging trade-offs between learning rates and comfort violation at the early stages of training through, e.g., the use of learning rate schedules, starting with a high learning rate and gradually decreasing it as training progresses. The great generalization of the federated agent provides another promising direction for future research. Practical implementations could benefit from integrating additional techniques in a hybrid manner, such as rule-based controllers (including the PID tested as baseline) or model-based approaches (if feasible) for improving early sample efficiency. Another promising direction is to focus on transfer learning from simulated to real environments, where a pre-trained agent is deployed and tuned on real buildings. Alternatively, the federated agent could also be pre-trained on historical data. 

Finally, our choice of federated learning is driven by its pragmatic benefits — improving generalization, learning speed, and stability through the aggregation of local updates while preserving data privacy. In contrast, meta-reinforcement learning (meta-RL), though promising for showing rapid adaptation between unseen tasks, still faces practical challenges, such as the need for meticulously curated task distributions and increased computational complexity. Nonetheless, as it develops further, meta-RL represents an interesting avenue for future research on autonomous HVAC control.

\section*{Acknowledgements}
We acknowledge the computational resources provided by the Aalto Science-IT project. This research received financial support from the Research Council of Finland (decision number 348092).

\bibliographystyle{elsarticle-harv} 
\bibliography{main.bib}

\begin{thebibliography}{48}
\expandafter\ifx\csname natexlab\endcsname\relax\def\natexlab#1{#1}\fi
\providecommand{\url}[1]{\texttt{#1}}
\providecommand{\href}[2]{#2}
\providecommand{\path}[1]{#1}
\providecommand{\DOIprefix}{doi:}
\providecommand{\ArXivprefix}{arXiv:}
\providecommand{\URLprefix}{URL: }
\providecommand{\Pubmedprefix}{pmid:}
\providecommand{\doi}[1]{\href{http://dx.doi.org/#1}{\path{#1}}}
\providecommand{\Pubmed}[1]{\href{pmid:#1}{\path{#1}}}
\providecommand{\bibinfo}[2]{#2}
\ifx\xfnm\relax \def\xfnm[#1]{\unskip,\space#1}\fi
\bibitem[{Anderson et~al.(1997)Anderson, Hittle, Katz and Kretchmar}]{anderson1997synthesis}
\bibinfo{author}{Anderson, C.W.}, \bibinfo{author}{Hittle, D.C.}, \bibinfo{author}{Katz, A.D.}, \bibinfo{author}{Kretchmar, R.M.}, \bibinfo{year}{1997}.
\newblock \bibinfo{title}{Synthesis of reinforcement learning, neural networks and pi control applied to a simulated heating coil}.
\newblock \bibinfo{journal}{Artificial Intelligence in Engineering} \bibinfo{volume}{11}, \bibinfo{pages}{421--429}.
\bibitem[{Barrett and Linder(2015)}]{barrett2015autonomous}
\bibinfo{author}{Barrett, E.}, \bibinfo{author}{Linder, S.}, \bibinfo{year}{2015}.
\newblock \bibinfo{title}{Autonomous hvac control, a reinforcement learning approach}, in: \bibinfo{booktitle}{Machine Learning and Knowledge Discovery in Databases: European Conference, ECML PKDD 2015, Porto, Portugal, September 7-11, 2015, Proceedings, Part III 15}, \bibinfo{organization}{Springer}. pp. \bibinfo{pages}{3--19}.
\bibitem[{Biemann et~al.(2021)Biemann, Scheller, Liu and Huang}]{biemann2021experimental}
\bibinfo{author}{Biemann, M.}, \bibinfo{author}{Scheller, F.}, \bibinfo{author}{Liu, X.}, \bibinfo{author}{Huang, L.}, \bibinfo{year}{2021}.
\newblock \bibinfo{title}{Experimental evaluation of model-free reinforcement learning algorithms for continuous hvac control}.
\newblock \bibinfo{journal}{Applied Energy} \bibinfo{volume}{298}, \bibinfo{pages}{117164}.
\bibitem[{Brockman et~al.(2016)Brockman, Cheung, Pettersson, Schneider, Schulman, Tang and Zaremba}]{brockman2016openai}
\bibinfo{author}{Brockman, G.}, \bibinfo{author}{Cheung, V.}, \bibinfo{author}{Pettersson, L.}, \bibinfo{author}{Schneider, J.}, \bibinfo{author}{Schulman, J.}, \bibinfo{author}{Tang, J.}, \bibinfo{author}{Zaremba, W.}, \bibinfo{year}{2016}.
\newblock \bibinfo{title}{Openai gym}.
\newblock \bibinfo{journal}{arXiv preprint arXiv:1606.01540} .
\bibitem[{Chen et~al.(2020)Chen, Cai and Berg{\'e}s}]{chen2020gnu}
\bibinfo{author}{Chen, B.}, \bibinfo{author}{Cai, Z.}, \bibinfo{author}{Berg{\'e}s, M.}, \bibinfo{year}{2020}.
\newblock \bibinfo{title}{Gnu-rl: A practical and scalable reinforcement learning solution for building hvac control using a differentiable mpc policy}.
\newblock \bibinfo{journal}{Frontiers in Built Environment} \bibinfo{volume}{6}, \bibinfo{pages}{562239}.
\bibitem[{Costanzo et~al.(2016)Costanzo, Iacovella, Ruelens, Leurs and Claessens}]{costanzo2016experimental}
\bibinfo{author}{Costanzo, G.T.}, \bibinfo{author}{Iacovella, S.}, \bibinfo{author}{Ruelens, F.}, \bibinfo{author}{Leurs, T.}, \bibinfo{author}{Claessens, B.J.}, \bibinfo{year}{2016}.
\newblock \bibinfo{title}{Experimental analysis of data-driven control for a building heating system}.
\newblock \bibinfo{journal}{Sustainable Energy, Grids and Networks} \bibinfo{volume}{6}, \bibinfo{pages}{81--90}.
\bibitem[{Du et~al.(2021)Du, Zandi, Kotevska, Kurte, Munk, Amasyali, Mckee and Li}]{du2021intelligent}
\bibinfo{author}{Du, Y.}, \bibinfo{author}{Zandi, H.}, \bibinfo{author}{Kotevska, O.}, \bibinfo{author}{Kurte, K.}, \bibinfo{author}{Munk, J.}, \bibinfo{author}{Amasyali, K.}, \bibinfo{author}{Mckee, E.}, \bibinfo{author}{Li, F.}, \bibinfo{year}{2021}.
\newblock \bibinfo{title}{Intelligent multi-zone residential hvac control strategy based on deep reinforcement learning}.
\newblock \bibinfo{journal}{Applied Energy} \bibinfo{volume}{281}, \bibinfo{pages}{116117}.
\bibitem[{Fawzy et~al.(2020)Fawzy, Osman, Doran and Rooney}]{fawzy2020strategies}
\bibinfo{author}{Fawzy, S.}, \bibinfo{author}{Osman, A.I.}, \bibinfo{author}{Doran, J.}, \bibinfo{author}{Rooney, D.W.}, \bibinfo{year}{2020}.
\newblock \bibinfo{title}{Strategies for mitigation of climate change: a review}.
\newblock \bibinfo{journal}{Environmental Chemistry Letters} \bibinfo{volume}{18}, \bibinfo{pages}{2069--2094}.
\bibitem[{Fujita et~al.(2022)Fujita, Fujimura, Sun, Esaki and Ochiai}]{fujita2022federated}
\bibinfo{author}{Fujita, K.}, \bibinfo{author}{Fujimura, S.}, \bibinfo{author}{Sun, Y.}, \bibinfo{author}{Esaki, H.}, \bibinfo{author}{Ochiai, H.}, \bibinfo{year}{2022}.
\newblock \bibinfo{title}{Federated reinforcement learning for the building facilities}, in: \bibinfo{booktitle}{2022 IEEE International Conference on Omni-layer Intelligent Systems (COINS)}, \bibinfo{organization}{IEEE}. pp. \bibinfo{pages}{1--6}.
\bibitem[{Gao and Wang(2023)}]{gao2023comparative}
\bibinfo{author}{Gao, C.}, \bibinfo{author}{Wang, D.}, \bibinfo{year}{2023}.
\newblock \bibinfo{title}{Comparative study of model-based and model-free reinforcement learning control performance in hvac systems}.
\newblock \bibinfo{journal}{Journal of Building Engineering} \bibinfo{volume}{74}, \bibinfo{pages}{106852}.
\bibitem[{Gao et~al.(2019)Gao, Li and Wen}]{gao2019energy}
\bibinfo{author}{Gao, G.}, \bibinfo{author}{Li, J.}, \bibinfo{author}{Wen, Y.}, \bibinfo{year}{2019}.
\newblock \bibinfo{title}{Energy-efficient thermal comfort control in smart buildings via deep reinforcement learning}.
\newblock \bibinfo{journal}{arXiv preprint arXiv:1901.04693} .
\bibitem[{Gao et~al.(2021)Gao, Wang, Liu, Billah and Campbell}]{gao2021decentralized}
\bibinfo{author}{Gao, J.}, \bibinfo{author}{Wang, W.}, \bibinfo{author}{Liu, Z.}, \bibinfo{author}{Billah, M.F.R.M.}, \bibinfo{author}{Campbell, B.}, \bibinfo{year}{2021}.
\newblock \bibinfo{title}{Decentralized federated learning framework for the neighborhood: a case study on residential building load forecasting}, in: \bibinfo{booktitle}{Proceedings of the 19th ACM Conference on Embedded Networked Sensor Systems}, pp. \bibinfo{pages}{453--459}.
\bibitem[{Guo et~al.(2020)Guo, Wang, Vishwanath, Xu and Li}]{guo2020towards}
\bibinfo{author}{Guo, Y.}, \bibinfo{author}{Wang, D.}, \bibinfo{author}{Vishwanath, A.}, \bibinfo{author}{Xu, C.}, \bibinfo{author}{Li, Q.}, \bibinfo{year}{2020}.
\newblock \bibinfo{title}{Towards federated learning for hvac analytics: A measurement study}, in: \bibinfo{booktitle}{Proceedings of the Eleventh ACM International Conference on Future Energy Systems}, pp. \bibinfo{pages}{68--73}.
\bibitem[{Haarnoja et~al.(2018c)Haarnoja, Ha, Zhou, Tan, Tucker and Levine}]{haarnoja2018softlearning}
\bibinfo{author}{Haarnoja, T.}, \bibinfo{author}{Ha, S.}, \bibinfo{author}{Zhou, A.}, \bibinfo{author}{Tan, J.}, \bibinfo{author}{Tucker, G.}, \bibinfo{author}{Levine, S.}, \bibinfo{year}{2018c}.
\newblock \bibinfo{title}{Learning to walk via deep reinforcement learning}.
\newblock \bibinfo{journal}{arXiv preprint arXiv:1812.11103} .
\bibitem[{Haarnoja et~al.(2018a)Haarnoja, Zhou, Abbeel and Levine}]{haarnoja2018softfirst}
\bibinfo{author}{Haarnoja, T.}, \bibinfo{author}{Zhou, A.}, \bibinfo{author}{Abbeel, P.}, \bibinfo{author}{Levine, S.}, \bibinfo{year}{2018a}.
\newblock \bibinfo{title}{Soft actor-critic: Off-policy maximum entropy deep reinforcement learning with a stochastic actor}, in: \bibinfo{booktitle}{International conference on machine learning}, \bibinfo{organization}{PMLR}. pp. \bibinfo{pages}{1861--1870}.
\bibitem[{Haarnoja et~al.(2018b)Haarnoja, Zhou, Hartikainen, Tucker, Ha, Tan, Kumar, Zhu, Gupta, Abbeel et~al.}]{haarnoja2018softapp}
\bibinfo{author}{Haarnoja, T.}, \bibinfo{author}{Zhou, A.}, \bibinfo{author}{Hartikainen, K.}, \bibinfo{author}{Tucker, G.}, \bibinfo{author}{Ha, S.}, \bibinfo{author}{Tan, J.}, \bibinfo{author}{Kumar, V.}, \bibinfo{author}{Zhu, H.}, \bibinfo{author}{Gupta, A.}, \bibinfo{author}{Abbeel, P.}, et~al., \bibinfo{year}{2018b}.
\newblock \bibinfo{title}{Soft actor-critic algorithms and applications}.
\newblock \bibinfo{journal}{arXiv preprint arXiv:1812.05905} .
\bibitem[{Hagstr{\"o}m(2023)}]{hagstrom2023using}
\bibinfo{author}{Hagstr{\"o}m, F.}, \bibinfo{year}{2023}.
\newblock \bibinfo{title}{Using federated learning techniques to train deep reinforcement learning agents for hvac control} .
\bibitem[{Henze and Schoenmann(2003)}]{henze2003evaluation}
\bibinfo{author}{Henze, G.P.}, \bibinfo{author}{Schoenmann, J.}, \bibinfo{year}{2003}.
\newblock \bibinfo{title}{Evaluation of reinforcement learning control for thermal energy storage systems}.
\newblock \bibinfo{journal}{HVAC\&R Research} \bibinfo{volume}{9}, \bibinfo{pages}{259--275}.
\bibitem[{Hsu et~al.(2019)Hsu, Qi and Brown}]{hsu2019measuring}
\bibinfo{author}{Hsu, T.M.H.}, \bibinfo{author}{Qi, H.}, \bibinfo{author}{Brown, M.}, \bibinfo{year}{2019}.
\newblock \bibinfo{title}{Measuring the effects of non-identical data distribution for federated visual classification}.
\newblock \bibinfo{journal}{arXiv preprint arXiv:1909.06335} .
\bibitem[{Jiménez-Raboso et~al.(2021)Jiménez-Raboso, Campoy-Nieves, Manjavacas-Lucas, Gómez-Romero and Molina-Solana}]{2021sinergym}
\bibinfo{author}{Jiménez-Raboso, J.}, \bibinfo{author}{Campoy-Nieves, A.}, \bibinfo{author}{Manjavacas-Lucas, A.}, \bibinfo{author}{Gómez-Romero, J.}, \bibinfo{author}{Molina-Solana, M.}, \bibinfo{year}{2021}.
\newblock \bibinfo{title}{Sinergym: A building simulation and control framework for training reinforcement learning agents}, in: \bibinfo{booktitle}{Proceedings of the 8th ACM International Conference on Systems for Energy-Efficient Buildings, Cities, and Transportation}, \bibinfo{publisher}{Association for Computing Machinery}, \bibinfo{address}{New York, NY, USA}. p. \bibinfo{pages}{319–323}.
\newblock \URLprefix \url{https://doi.org/10.1145/3486611.3488729}, \DOIprefix\doi{10.1145/3486611.3488729}.
\bibitem[{Khalil et~al.(2021)Khalil, Esseghir and Merghem-Boulahia}]{khalil2021federated}
\bibinfo{author}{Khalil, M.}, \bibinfo{author}{Esseghir, M.}, \bibinfo{author}{Merghem-Boulahia, L.}, \bibinfo{year}{2021}.
\newblock \bibinfo{title}{Federated learning for energy-efficient thermal comfort control service in smart buildings}, in: \bibinfo{booktitle}{2021 IEEE Global Communications Conference (GLOBECOM)}, \bibinfo{organization}{IEEE}. pp. \bibinfo{pages}{01--06}.
\bibitem[{Khalil et~al.(2022)Khalil, Esseghir and Merghem-Boulahia}]{khalil2022federated}
\bibinfo{author}{Khalil, M.}, \bibinfo{author}{Esseghir, M.}, \bibinfo{author}{Merghem-Boulahia, L.}, \bibinfo{year}{2022}.
\newblock \bibinfo{title}{A federated learning approach for thermal comfort management}.
\newblock \bibinfo{journal}{Advanced Engineering Informatics} \bibinfo{volume}{52}, \bibinfo{pages}{101526}.
\bibitem[{Kingma and Ba(2014)}]{kingma2014adam}
\bibinfo{author}{Kingma, D.P.}, \bibinfo{author}{Ba, J.}, \bibinfo{year}{2014}.
\newblock \bibinfo{title}{Adam: A method for stochastic optimization}.
\newblock \bibinfo{journal}{arXiv preprint arXiv:1412.6980} .
\bibitem[{Kochenderfer et~al.(2022)Kochenderfer, Wheeler and Wray}]{kochenderfer2022algorithms}
\bibinfo{author}{Kochenderfer, M.J.}, \bibinfo{author}{Wheeler, T.A.}, \bibinfo{author}{Wray, K.H.}, \bibinfo{year}{2022}.
\newblock \bibinfo{title}{Algorithms for decision making}.
\newblock \bibinfo{publisher}{MIT press}.
\bibitem[{Lee et~al.(2021)Lee, Xie and Choi}]{lee2021privacy}
\bibinfo{author}{Lee, S.}, \bibinfo{author}{Xie, L.}, \bibinfo{author}{Choi, D.H.}, \bibinfo{year}{2021}.
\newblock \bibinfo{title}{Privacy-preserving energy management of a shared energy storage system for smart buildings: A federated deep reinforcement learning approach}.
\newblock \bibinfo{journal}{Sensors} \bibinfo{volume}{21}, \bibinfo{pages}{4898}.
\bibitem[{Li and Xia(2015)}]{li2015multi}
\bibinfo{author}{Li, B.}, \bibinfo{author}{Xia, L.}, \bibinfo{year}{2015}.
\newblock \bibinfo{title}{A multi-grid reinforcement learning method for energy conservation and comfort of hvac in buildings}, in: \bibinfo{booktitle}{2015 IEEE International Conference on Automation Science and Engineering (CASE)}, \bibinfo{organization}{IEEE}. pp. \bibinfo{pages}{444--449}.
\bibitem[{Liu and Henze(2005)}]{liu2005evaluation}
\bibinfo{author}{Liu, S.}, \bibinfo{author}{Henze, G.P.}, \bibinfo{year}{2005}.
\newblock \bibinfo{title}{Evaluation of reinforcement learning for optimal control of building active and passive thermal storage inventory}, in: \bibinfo{booktitle}{International Solar Energy Conference}, pp. \bibinfo{pages}{301--311}.
\bibitem[{Liu and Henze(2006a)}]{liu2006experimental-1}
\bibinfo{author}{Liu, S.}, \bibinfo{author}{Henze, G.P.}, \bibinfo{year}{2006a}.
\newblock \bibinfo{title}{Experimental analysis of simulated reinforcement learning control for active and passive building thermal storage inventory. part 1. theoretical foundation}.
\newblock \bibinfo{journal}{Energy and buildings} \bibinfo{volume}{38}, \bibinfo{pages}{142--147}.
\bibitem[{Liu and Henze(2006b)}]{liu2006experimental-2}
\bibinfo{author}{Liu, S.}, \bibinfo{author}{Henze, G.P.}, \bibinfo{year}{2006b}.
\newblock \bibinfo{title}{Experimental analysis of simulated reinforcement learning control for active and passive building thermal storage inventory: Part 2: Results and analysis}.
\newblock \bibinfo{journal}{Energy and buildings} \bibinfo{volume}{38}, \bibinfo{pages}{148--161}.
\bibitem[{Lu et~al.(2023)Lu, Cui, Wang, Sun and Liu}]{lu2023residential}
\bibinfo{author}{Lu, Y.}, \bibinfo{author}{Cui, L.}, \bibinfo{author}{Wang, Y.}, \bibinfo{author}{Sun, J.}, \bibinfo{author}{Liu, L.}, \bibinfo{year}{2023}.
\newblock \bibinfo{title}{Residential energy consumption forecasting based on federated reinforcement learning with data privacy protection}.
\newblock \bibinfo{journal}{CMES-Computer Modeling in Engineering \& Sciences} \bibinfo{volume}{137}.
\bibitem[{McMahan et~al.(2017)McMahan, Moore, Ramage, Hampson and y~Arcas}]{mcmahan2017communication}
\bibinfo{author}{McMahan, B.}, \bibinfo{author}{Moore, E.}, \bibinfo{author}{Ramage, D.}, \bibinfo{author}{Hampson, S.}, \bibinfo{author}{y~Arcas, B.A.}, \bibinfo{year}{2017}.
\newblock \bibinfo{title}{Communication-efficient learning of deep networks from decentralized data}, in: \bibinfo{booktitle}{Artificial intelligence and statistics}, \bibinfo{organization}{PMLR}. pp. \bibinfo{pages}{1273--1282}.
\bibitem[{Mozer(1998)}]{mozer1998neural}
\bibinfo{author}{Mozer, M.C.}, \bibinfo{year}{1998}.
\newblock \bibinfo{title}{The neural network house: An environment that adapts to its inhabitants}, in: \bibinfo{booktitle}{Proc. AAAI Spring Symp. Intelligent Environments}.
\bibitem[{Nagy et~al.(2018)Nagy, Kazmi, Cheaib and Driesen}]{nagy2018deep}
\bibinfo{author}{Nagy, A.}, \bibinfo{author}{Kazmi, H.}, \bibinfo{author}{Cheaib, F.}, \bibinfo{author}{Driesen, J.}, \bibinfo{year}{2018}.
\newblock \bibinfo{title}{Deep reinforcement learning for optimal control of space heating}.
\newblock \bibinfo{journal}{arXiv preprint arXiv:1805.03777} .
\bibitem[{Paszke et~al.(2017)Paszke, Gross, Chintala, Chanan, Yang, DeVito, Lin, Desmaison, Antiga and Lerer}]{paszke2017automatic}
\bibinfo{author}{Paszke, A.}, \bibinfo{author}{Gross, S.}, \bibinfo{author}{Chintala, S.}, \bibinfo{author}{Chanan, G.}, \bibinfo{author}{Yang, E.}, \bibinfo{author}{DeVito, Z.}, \bibinfo{author}{Lin, Z.}, \bibinfo{author}{Desmaison, A.}, \bibinfo{author}{Antiga, L.}, \bibinfo{author}{Lerer, A.}, \bibinfo{year}{2017}.
\newblock \bibinfo{title}{Automatic differentiation in pytorch} .
\bibitem[{Perera and Kamalaruban(2021)}]{perera2021applications}
\bibinfo{author}{Perera, A.}, \bibinfo{author}{Kamalaruban, P.}, \bibinfo{year}{2021}.
\newblock \bibinfo{title}{Applications of reinforcement learning in energy systems}.
\newblock \bibinfo{journal}{Renewable and Sustainable Energy Reviews} \bibinfo{volume}{137}, \bibinfo{pages}{110618}.
\bibitem[{Raffin et~al.(2021)Raffin, Hill, Gleave, Kanervisto, Ernestus and Dormann}]{stable-baselines3}
\bibinfo{author}{Raffin, A.}, \bibinfo{author}{Hill, A.}, \bibinfo{author}{Gleave, A.}, \bibinfo{author}{Kanervisto, A.}, \bibinfo{author}{Ernestus, M.}, \bibinfo{author}{Dormann, N.}, \bibinfo{year}{2021}.
\newblock \bibinfo{title}{Stable-baselines3: Reliable reinforcement learning implementations}.
\newblock \bibinfo{journal}{Journal of Machine Learning Research} \bibinfo{volume}{22}, \bibinfo{pages}{1--8}.
\newblock \URLprefix \url{http://jmlr.org/papers/v22/20-1364.html}.
\bibitem[{Reddi et~al.(2020)Reddi, Charles, Zaheer, Garrett, Rush, Kone{\v{c}}n{\`y}, Kumar and McMahan}]{reddi2020adaptive}
\bibinfo{author}{Reddi, S.}, \bibinfo{author}{Charles, Z.}, \bibinfo{author}{Zaheer, M.}, \bibinfo{author}{Garrett, Z.}, \bibinfo{author}{Rush, K.}, \bibinfo{author}{Kone{\v{c}}n{\`y}, J.}, \bibinfo{author}{Kumar, S.}, \bibinfo{author}{McMahan, H.B.}, \bibinfo{year}{2020}.
\newblock \bibinfo{title}{Adaptive federated optimization}.
\newblock \bibinfo{journal}{arXiv preprint arXiv:2003.00295} .
\bibitem[{Ruelens et~al.(2016)Ruelens, Claessens, Vandael, De~Schutter, Babu{\v{s}}ka and Belmans}]{ruelens2016residential}
\bibinfo{author}{Ruelens, F.}, \bibinfo{author}{Claessens, B.J.}, \bibinfo{author}{Vandael, S.}, \bibinfo{author}{De~Schutter, B.}, \bibinfo{author}{Babu{\v{s}}ka, R.}, \bibinfo{author}{Belmans, R.}, \bibinfo{year}{2016}.
\newblock \bibinfo{title}{Residential demand response of thermostatically controlled loads using batch reinforcement learning}.
\newblock \bibinfo{journal}{IEEE Transactions on Smart Grid} \bibinfo{volume}{8}, \bibinfo{pages}{2149--2159}.
\bibitem[{Sun et~al.(2013)Sun, Luh, Jia and Yan}]{sun2013event}
\bibinfo{author}{Sun, B.}, \bibinfo{author}{Luh, P.B.}, \bibinfo{author}{Jia, Q.S.}, \bibinfo{author}{Yan, B.}, \bibinfo{year}{2013}.
\newblock \bibinfo{title}{Event-based optimization with non-stationary uncertainties to save energy costs of hvac systems in buildings}, in: \bibinfo{booktitle}{2013 IEEE International Conference on Automation Science and Engineering (CASE)}, \bibinfo{organization}{IEEE}. pp. \bibinfo{pages}{436--441}.
\bibitem[{Sutton and Barto(2018)}]{sutton2018reinforcement}
\bibinfo{author}{Sutton, R.S.}, \bibinfo{author}{Barto, A.G.}, \bibinfo{year}{2018}.
\newblock \bibinfo{title}{Reinforcement learning: An introduction}.
\newblock \bibinfo{publisher}{MIT press}.
\bibitem[{TC et~al.(2016)}]{tc2016data}
\bibinfo{author}{TC, A.}, et~al., \bibinfo{year}{2016}.
\newblock \bibinfo{title}{Data center power equipment thermal guidelines and best practices}.
\newblock \bibinfo{journal}{ASHRAE TC 9.9, ASHRAE, USA} .
\bibitem[{Tenison et~al.(2022)Tenison, Sreeramadas, Mugunthan, Oyallon, Belilovsky and Rish}]{tenison2022gradient}
\bibinfo{author}{Tenison, I.}, \bibinfo{author}{Sreeramadas, S.A.}, \bibinfo{author}{Mugunthan, V.}, \bibinfo{author}{Oyallon, E.}, \bibinfo{author}{Belilovsky, E.}, \bibinfo{author}{Rish, I.}, \bibinfo{year}{2022}.
\newblock \bibinfo{title}{Gradient masked averaging for federated learning}.
\newblock \bibinfo{journal}{arXiv preprint arXiv:2201.11986} .
\bibitem[{V{\'a}zquez-Canteli and Nagy(2019)}]{vazquez2019reinforcement}
\bibinfo{author}{V{\'a}zquez-Canteli, J.R.}, \bibinfo{author}{Nagy, Z.}, \bibinfo{year}{2019}.
\newblock \bibinfo{title}{Reinforcement learning for demand response: A review of algorithms and modeling techniques}.
\newblock \bibinfo{journal}{Applied energy} \bibinfo{volume}{235}, \bibinfo{pages}{1072--1089}.
\bibitem[{Wang and Hong(2020)}]{wang2020reinforcement}
\bibinfo{author}{Wang, Z.}, \bibinfo{author}{Hong, T.}, \bibinfo{year}{2020}.
\newblock \bibinfo{title}{Reinforcement learning for building controls: The opportunities and challenges}.
\newblock \bibinfo{journal}{Applied Energy} \bibinfo{volume}{269}, \bibinfo{pages}{115036}.
\bibitem[{Wang et~al.(2022)Wang, Yu and Zhang}]{wang2022privacy}
\bibinfo{author}{Wang, Z.}, \bibinfo{author}{Yu, P.}, \bibinfo{author}{Zhang, H.}, \bibinfo{year}{2022}.
\newblock \bibinfo{title}{Privacy-preserving regulation capacity evaluation for hvac systems in heterogeneous buildings based on federated learning and transfer learning}.
\newblock \bibinfo{journal}{IEEE Transactions on Smart Grid} .
\bibitem[{Wei et~al.(2017)Wei, Wang and Zhu}]{wei2017deep}
\bibinfo{author}{Wei, T.}, \bibinfo{author}{Wang, Y.}, \bibinfo{author}{Zhu, Q.}, \bibinfo{year}{2017}.
\newblock \bibinfo{title}{Deep reinforcement learning for building hvac control}, in: \bibinfo{booktitle}{Proceedings of the 54th annual design automation conference 2017}, pp. \bibinfo{pages}{1--6}.
\bibitem[{Weinberg et~al.(2022)Weinberg, Wang, Timoudas and Fischione}]{weinberg2022review}
\bibinfo{author}{Weinberg, D.}, \bibinfo{author}{Wang, Q.}, \bibinfo{author}{Timoudas, T.O.}, \bibinfo{author}{Fischione, C.}, \bibinfo{year}{2022}.
\newblock \bibinfo{title}{A review of reinforcement learning for controlling building energy systems from a computer science perspective}.
\newblock \bibinfo{journal}{Sustainable cities and society} , \bibinfo{pages}{104351}.
\bibitem[{Wiering and Van~Otterlo(2012)}]{wiering2012reinforcement}
\bibinfo{author}{Wiering, M.A.}, \bibinfo{author}{Van~Otterlo, M.}, \bibinfo{year}{2012}.
\newblock \bibinfo{title}{Reinforcement learning}.
\newblock \bibinfo{journal}{Adaptation, learning, and optimization} \bibinfo{volume}{12}, \bibinfo{pages}{729}.

\end{thebibliography}






\newpage 
\onecolumn
\appendix

\section{Soft Actor-Critic (SAC) pseudo-code}
\label{appendix:SAC}

The final practical SAC algorithm used in our experiments includes a few additional features when compared to the algorithm presented in section \ref{sec:SAC}. The final algorithm learns two concurrent soft Q-functions, parameterized by $\theta_i, i \in \{1, 2\}$, which are trained independently to minimize $\mathcal{L}(\theta_i)$ in equation \eqref{eq:loss}. They both have their respective target networks $\bar{\theta}_i, i \in \{1, 2\}$. The equation \eqref{eq:SACtarget} for the target $y$ is modified to utilize the minimum of the two Q-functions
\begin{align}
    y = r + \gamma \big(\min_{i=1,2} Q_{\bar{\theta}_i} (s', \Tilde{a}') - \alpha \log \pi (\Tilde{a}' | s') \big), \quad \Tilde{a}' \sim \pi_\phi(\cdot | s') \label{eq:SACTargetMin}
\end{align}
and similarly for the performance $J(\phi)$ in equation \eqref{eq:ReParamPolicyUpdate}
\begin{align}
     J(\phi) = \mathbb{E}_{s \sim \mathcal{D}, \epsilon \sim \mathcal{N}} \Big[  \alpha \log \pi(f_\phi(\epsilon, s) | s) - \min_{i=1,2} Q_{\theta_i} (s, f_\phi(\epsilon, s)) \Big]. \label{eq:ReParamPolicyUpdateMin}
\end{align}
This double Q-learning trick is used to mitigate positive bias in the policy improvement step, which can degrade performance \citep{haarnoja2018softapp}.

The SAC algorithm is particularly sensitive to the temperature $\alpha$, which has to be fine-tuned to the task at hand in order to achieve appropriate performance. \cite{haarnoja2018softapp} develop a method for automatically adjusting its value during training to stabilise learning across different tasks. The temperature is updated at each gradient step by minimizing the following objective
\begin{align}
    J(\alpha) = \mathbb{E}_{s \sim \mathcal{D}, a \sim \pi_\phi} \Big[ -\alpha \log \pi_\phi(a | s) - \alpha \Tilde{\mathcal{H}} \Big], \label{tempAdjust}
\end{align}
where $\Tilde{\mathcal{H}}$ is the minimum desired entropy. \cite{haarnoja2018softlearning} find that the algorithm is quite robust with respect to the minimum entropy, and generally setting it to $-1$ times the action dimension yields good results.

\begin{algorithm}
\caption{SAC}
\label{alg:SAC}
\begin{algorithmic}[1]
    \STATE \textbf{Initialize:} \\
    Critic networks $Q_{\theta_1}$, $Q_{\theta_2}$ and actor network $\pi_{\phi'}$ with random parameters $\theta_1$, $\theta_2$, $\phi$. \\
    Target networks $\theta_1' \leftarrow \theta_1$, $\theta_2' \leftarrow \theta_2$. \\
    Replay buffer $\mathcal{B}.$
    \FOR{each iteration}
        \FOR {each environment step}
            \STATE Sample action $a_t \sim \pi_\phi (\cdot | s_t)$ and observe reward $r_t$ and next state $s_{t+1}$.
            \STATE Store transition tuple $(s, a, r, s')$ in replay buffer $\mathcal{B}$.
        \ENDFOR
        \FOR {each gradient step}
            \STATE Sample mini-batch $\mathcal{D}$ from replay buffer $\mathcal{B}$.
            \STATE Compute targets $y$ for all $(s, a, r, s') \in \mathcal{D}$, equation \eqref{eq:SACTargetMin}
            \STATE
            \STATE Update critics: $\theta_i \leftarrow \theta_i - \lambda_{Q} \nabla_{\theta_i} \mathcal{L}(\theta_i)$, equation \eqref{eq:loss}.
            \STATE Update actor: $\phi \leftarrow \phi - \lambda_\pi \nabla_{\phi} J(\phi)$, equation \eqref{eq:ReParamPolicyUpdateMin}.
            \STATE Update temperature: $\alpha \leftarrow \alpha - \lambda_\alpha \nabla_{\alpha} J(\alpha)$, equation \eqref{tempAdjust}.
            \STATE
            \STATE Update target networks: $\Bar{\theta}_i \leftarrow \rho \theta_i + (1 - \rho) \Bar{\theta}_i$
        \ENDFOR
    \ENDFOR
\end{algorithmic}
\end{algorithm}

\section{State space}
\label{appendix:ObsSpace}

Table \ref{tab:obsVar} shows the complete list of observed state features in the data center environment.
\begin{table}[htb]
    \caption{Description of the state space.} \label{tab:obsVar}
    \centering
    \begin{tabular}{l c} \\ \toprule
         Feature & Unit  \\ \midrule 
         Site Outdoor air drybulb temperature & $\degree C$ \\ 
         Site Outdoor Air Relative Humidity &  \% \\
         Site Wind Speed & $m / s$ \\
         Site Wind Direction & degree \\
         Site Diffuse Solar Radiation Rate per Area & $W / m^2$ \\
         Site Direct Solar Radiation Rate per Area & $W / m^2$ \\
         Zone Air Temperature(West Zone) & $\degree C$ \\
         Zone Air Relative Humidity(West Zone) & \% \\
         Zone Air Temperature(East Zone) & $\degree C$ \\
         Zone Air Relative Humidity(East Zone) & \% \\
         Facility Total HVAC Electricity Demand Rate &  $W$ \\
         Facility Total Building Electricity Demand Rate & $W$ \\
         Forecasted Outdoor Air Drybulb Temp (+1h) & $\degree C$ \\
         Forecasted Outdoor Air Relative Humidity (+1h) & \% \\
         Forecasted Outdoor Air Drybulb Temp (+3h) &  $\degree C$ \\
         Forecasted Outdoor Air Relative Humidity (+3h) & \% \\
         Forecasted Outdoor Air Drybulb Temp (+6h) & $\degree C$ \\
         Forecasted Outdoor Air Relative Humidity (+6h) & \% \\
         \bottomrule
    \end{tabular}
\end{table}

\section{Weather Profiles} \label{appendix:weather}

To add stochasticity to the weather from year to year, Sinergym modifies the base weather profiles via the Ornstein-Uhlenbeck process at the beginning of each year. The Ornstein-Uhlenbeck process $X_t$ is defined by the stochastic differential equation
\begin{align}
    dX_t = \tau (\mu - X_t) dt + \sigma dW_t,
\end{align}
where $W_t$ is Brownian motion with unit variance, and $\tau, \sigma \geq 0$ and $\mu$ are parameters affecting the evolution of the process. In our experiments, we set $\tau = 0.001$, $\sigma = 2.0$ and $\mu = 0$.

As mentioned in section \ref{sec:MDPExperiment}, we include ``forecasted'' outside temperature and relative humidity in our observations. These forecasted values are retrieved from the base weather profile. Since the weather over each year is stochastically modified from the base weather profile, the base weather profile provides us with values that are close to the ``true'' observed values, much like a typical weather forecast. Hence the base profile gives us a good proxy for a real weather forecast.

\begin{table}[htb]
    \caption{The base weather files available in Sinergym. M.T is the mean temperature and M.H is the mean relative humidity of the file.} \label{tab:weatherFiles}
    \centering
    \begin{tabular}{l c c} \\ \toprule
         Location & M.T ($\degree C$) & M.H (\%) \\  \midrule
         Sydney, Australia & 17.9 & 68.83 \\ 
         Bogota, Colombia & 13.2 & 80.3 \\
         Granada, Spain & 14.84 & 59.83 \\
         Helsinki, Finland & 5.1 & 79.25 \\
         Tokyo, Japan  & 8.9 & 78.6 \\
         Antananarivo, Madagascar & 18.35 &  75.91 \\
         Arizona, USA  & 21.7 & 34.9 \\
         Colorado, USA & 9.95 & 55.25 \\
         Illinois, USA & 9.92 & 70.3 \\
         New York, USA & 12.6 & 68.5 \\
         Pennsylvania, USA  & 10.5 & 66.41 \\
         Washington, USA & 9.3 & 81.1 \\
         \bottomrule
    \end{tabular}
\end{table}

\newpage

\section{Implementation details} \label{appendix:ImpDetails}

We use the implementation of the SAC algorithm provided by the Stable Baselines3 framework \citep{stable-baselines3}, which offers reliable implementations of reinforcement learning algorithms in PyTorch \citep{paszke2017automatic}. The Q-value functions and policy are approximated using simple feed-forward neural networks with an input layer, two hidden layers, and an output layer.  As argued by \cite{biemann2021experimental}, in a real-world application, tuning all the hyperparameters of the algorithms becomes infeasible, hence the algorithms should perform well out-of-the-box. We therefore use the default hyperparameters of the Stable Baselines3 implementation. An exception is the rate at which the policy and Q-networks are updated. We set the training frequency to once every hour, i.e., after every 4 environment steps. At every update, the model takes a number of gradient steps equal to the number of environment steps taken between updates. See table \ref{tab:hyperparam:sac} for a list of the exact hyperparameter values used for SAC.

In deep reinforcement learning and when training neural networks in general, it is often useful to ensure that all the features of the input vectors are on the same scale. This prevents very large features from dominating the calculated gradient, as well as maintains a more consistent range of values for the gradient, which often leads to more stable and faster learning. Hence we normalize the observations. The reward also affects the scale of the gradient, and as such, normalizing the rewards can also have a stabilizing effect. We therefore normalize the rewards as well. We use the VecNormalize wrapper in Stable Baselines3 with default values to normalize using a moving average.

SAC learns a stochastic policy. However, \cite{haarnoja2018softapp} find that making the final policy deterministic often results in better performance than choosing actions stochastically, and so we set the SAC policy to be deterministic as well during evaluation. This is done by choosing the mean $\mu_\phi(s)$ of the policy distribution as the action.
\begin{table}[htb]
    \caption{SAC hyperparameters.} \label{tab:hyperparam:sac}
    \centering
    \begin{tabular}{l c l} \\ \toprule
         Critic networks && $24 \rightarrow 256 \rightarrow 256 \rightarrow 1$ \\
         Actor networks && $24 \rightarrow 256 \rightarrow 256 \rightarrow (2 \times 4)$ \\
         Activation function && ReLU \\
         Discount factor $\gamma$ && 0.99 \\
         Batch size && 256 \\
         Polyak averaging $\rho$ && 0.005 \\
         Buffer size && $10^6$ \\
         Temperature $\alpha$ && auto \\
         Target entropy && auto \\
         Train frequency && 4 \\
         Gradient steps && -1 (match train frequency)\\
         Learning starts && 100 \\
         Exploration (action) noise $\xi$ && None \\
         \bottomrule
    \end{tabular}
\end{table}

For the client optimizers, we only vary the learning rate, and use the default values for other hyperparameters. The default values are presented in table \ref{tab:cliopt:hyperparams}. For SGDM we set the momentum $\mu = 0.9$. 
\begin{table}
    \caption{Hyperparameters of the client optimizers. These are held constant throughout all experiments.} \label{tab:cliopt:hyperparams}
    \centering
    \begin{tabular}{clrrrrrr} \\ \toprule  
        && $\beta_1$ & $\beta_2$ & $\epsilon$ & $\lambda$ & $\mu$ & $\tau$\\ \midrule
        & Adam &  0.9 & 0.999 & $10^{-8}$ & 0 & - & - \\
        & SGD  & - & - & - & 0 & 0 & 0 \\
        & SGDM & - & - & - & 0 & 0.9 & 0 \\
        \bottomrule
    \end{tabular}
\end{table} 

\section{Sensitivity analysis} \label{appendix:sensAnal}

In this section, we perform a sensitivity analysis of the hyperparameters of both sets of experiments. We present the analysis of the client optimizers in section \ref{appendix:sensAnal:Client}, and analyse the server optimizers in section \ref{appendix:sensAnal:Server}.

\subsection{Client optimizers} \label{appendix:sensAnal:Client}

\begin{figure}[h]
    \centering
    \begin{subfigure}[b]{0.49\textwidth}
         \centering
         \includegraphics[width=\textwidth]{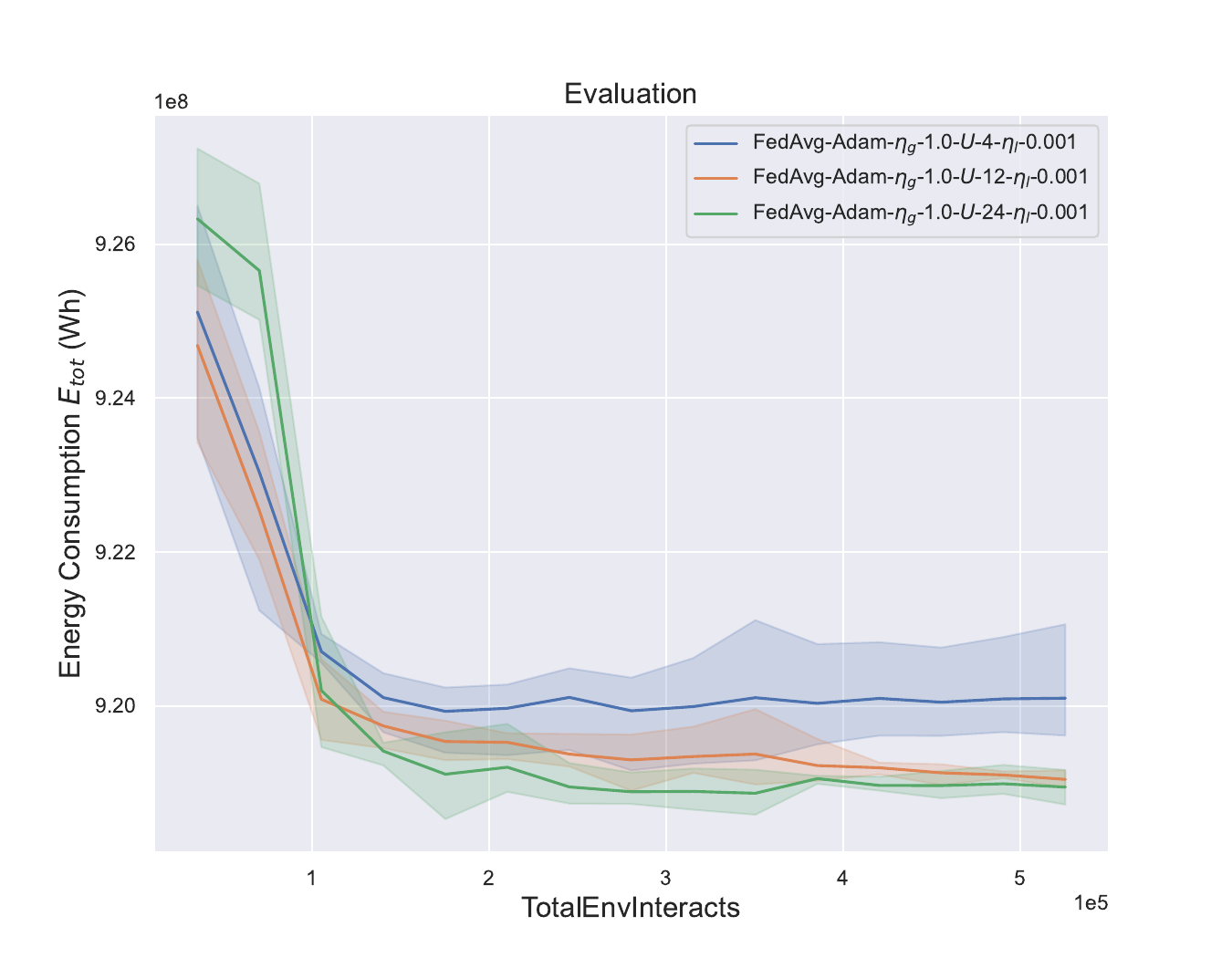}
         \caption{Energy consumption}
         \label{sens:Adam:lupr:energy}
    \end{subfigure}
    \hfill
    \begin{subfigure}[b]{0.49\textwidth}
         \centering
         \includegraphics[width=\textwidth]{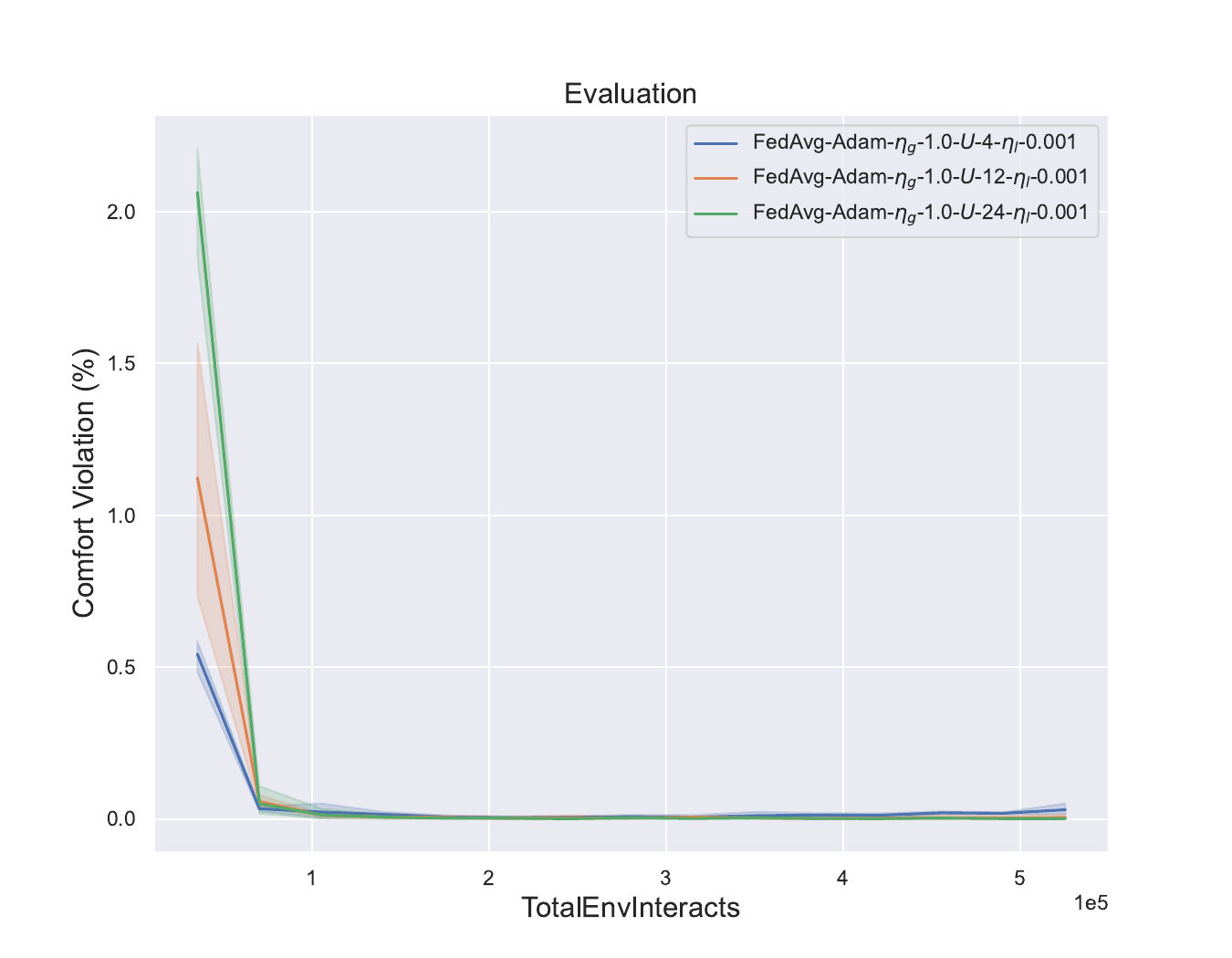}
         \caption{Comfort violation}
         \label{sens:Adam:lupr:comfort}
    \end{subfigure}
    \caption{Comparing the performance of FedAvg with Adam as client optimizer on the evaluation environment for different local updates per round $U$. We fix $\eta_l = 0.001$.}
    \label{sens:Adam:lupr}
\end{figure}

\begin{figure}[h]
    \centering
    \begin{subfigure}[b]{0.49\textwidth}
         \centering
         \includegraphics[width=\textwidth]{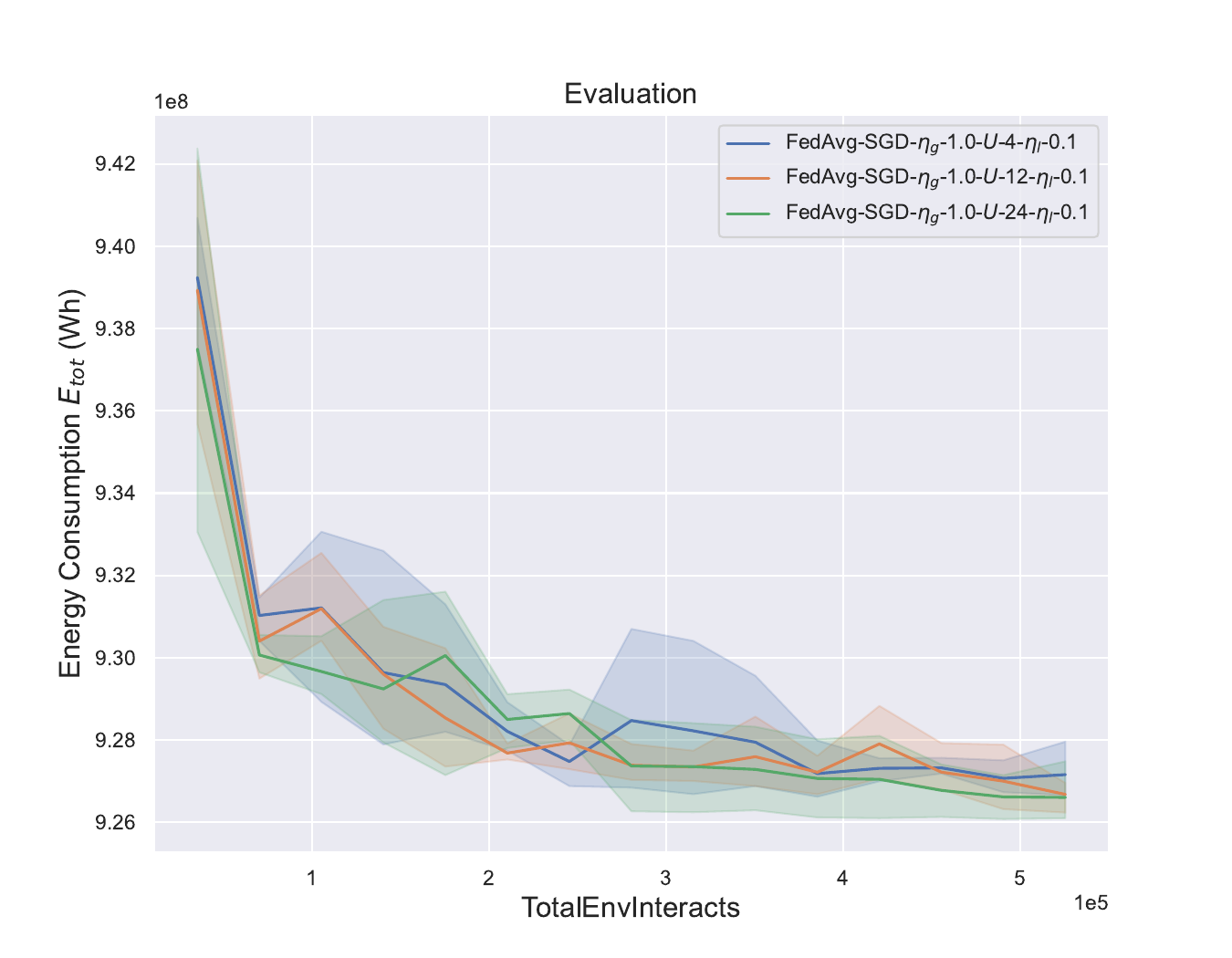}
         \caption{Energy consumption}
         \label{sens:SGD:lupr:energy}
    \end{subfigure}
    \hfill
    \begin{subfigure}[b]{0.49\textwidth}
         \centering
         \includegraphics[width=\textwidth]{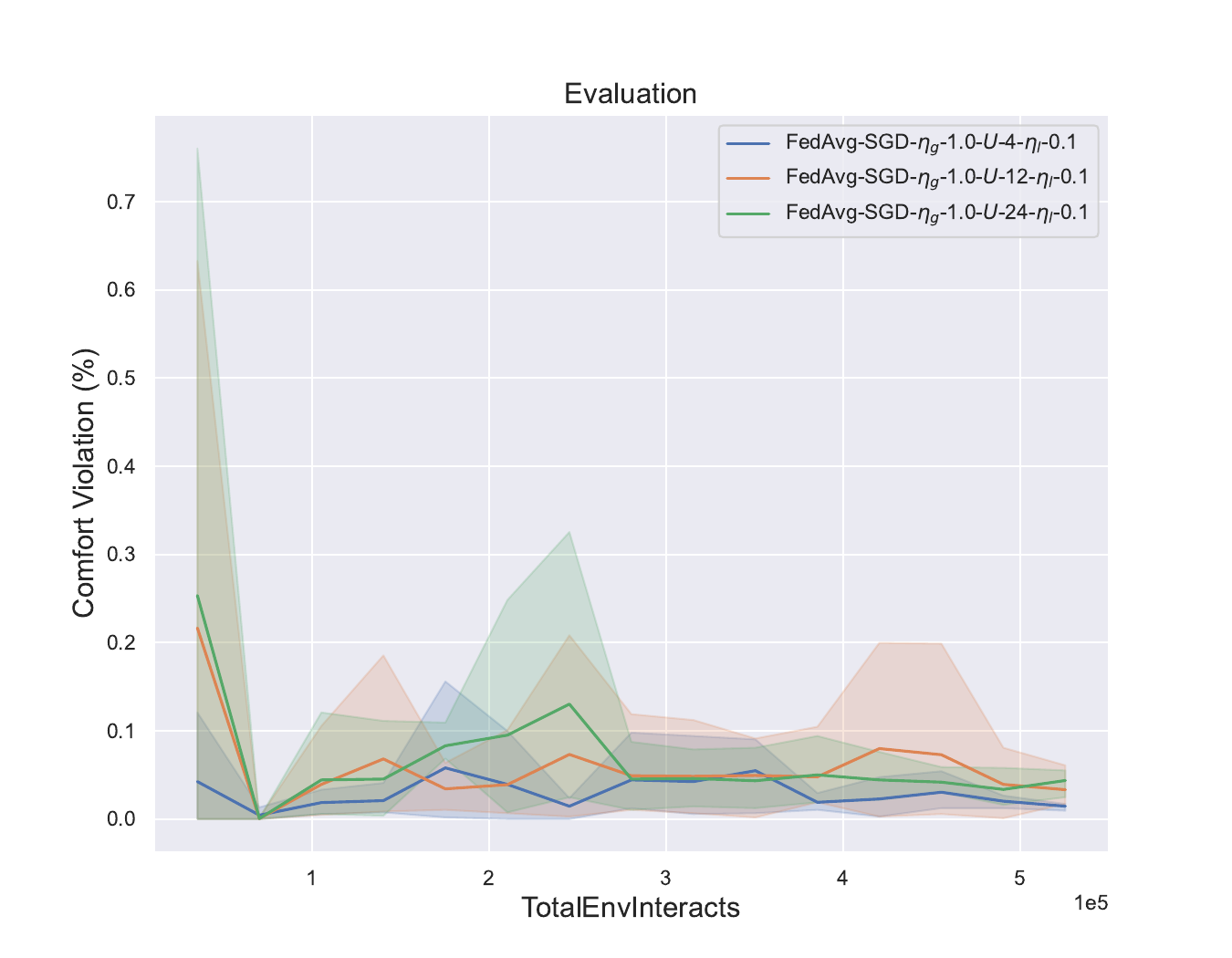}
         \caption{Comfort violation}
         \label{sens:SGD:lupr:comfort}
    \end{subfigure}
    \caption{Comparing the performance of FedAvg with SGD as client optimizer on the evaluation environment for different local updates per round $U$. We fix $\eta_l = 0.1$.}
    \label{sens:SGD:lupr}
\end{figure}

\begin{figure}[h]
    \centering
    \begin{subfigure}[b]{0.49\textwidth}
         \centering
         \includegraphics[width=\textwidth]{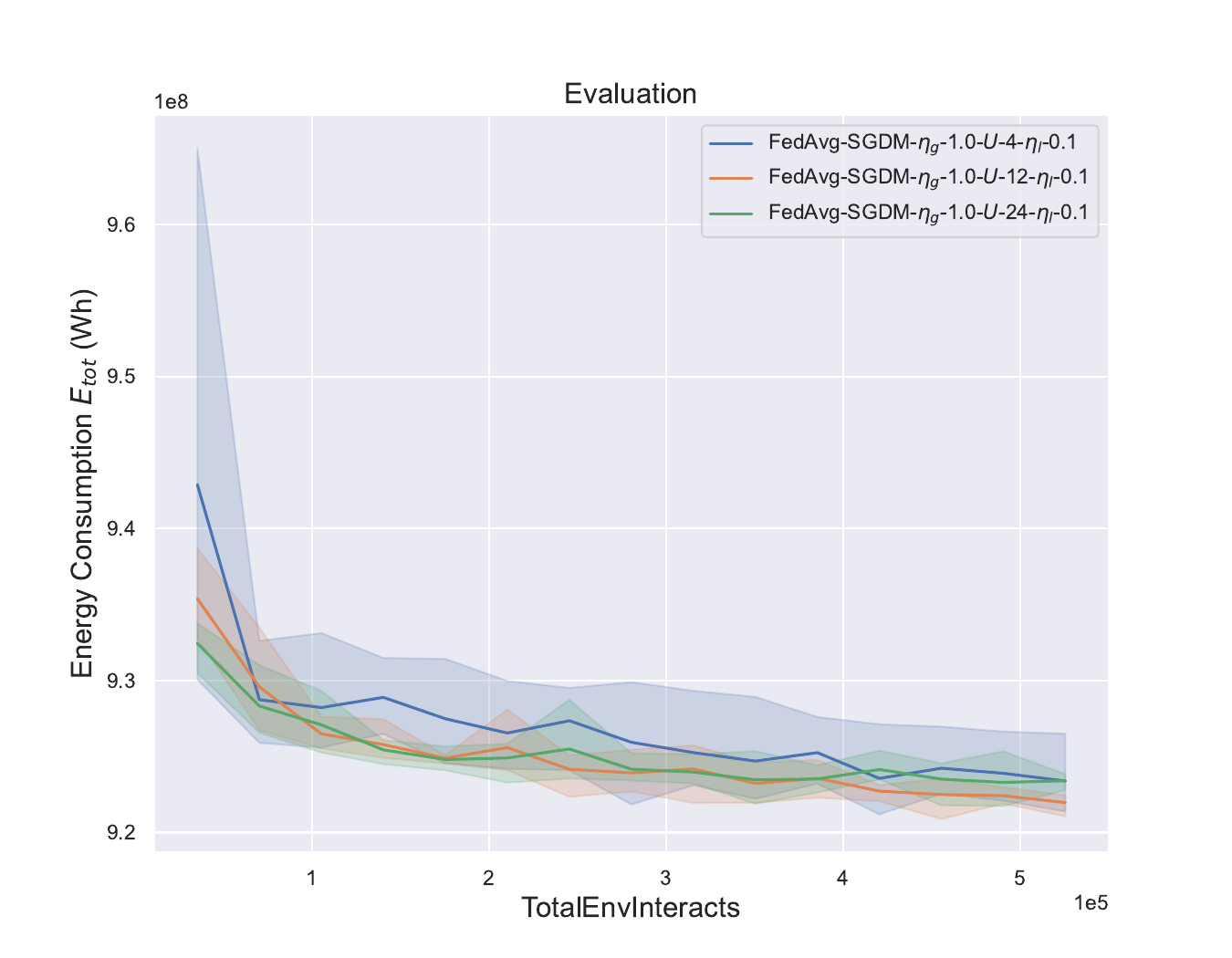}
         \caption{Energy consumption}
         \label{sens:SGDM:lupr:energy}
    \end{subfigure}
    \hfill
    \begin{subfigure}[b]{0.49\textwidth}
         \centering
         \includegraphics[width=\textwidth]{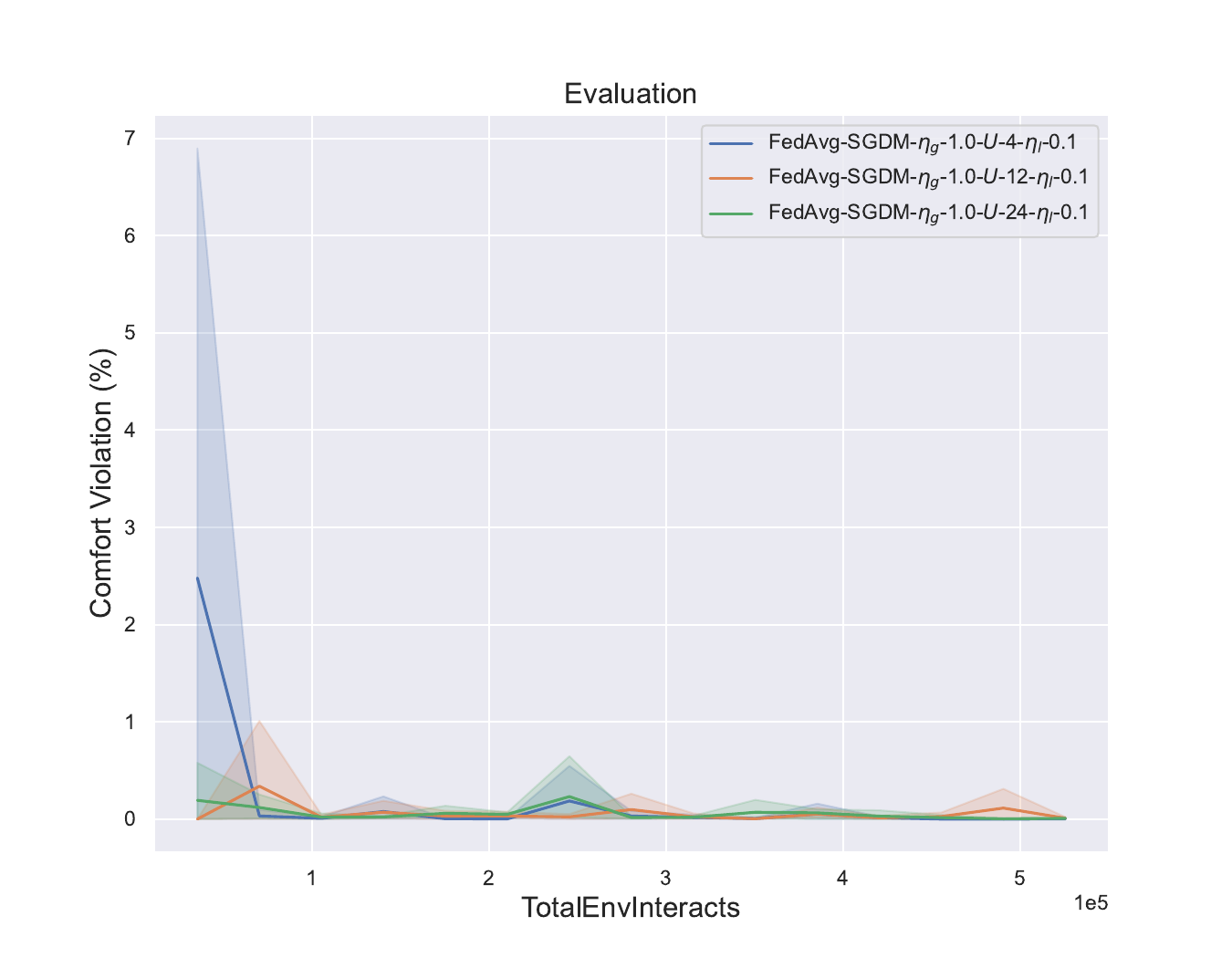}
         \caption{Comfort violation}
         \label{sens:SGDM:lupr:comfort}
    \end{subfigure}
    \caption{Comparing the performance of FedAvg with SGDM as client optimizer on the evaluation environment for different local updates per round $U$. We fix $\eta_l = 0.1$.}
    \label{sens:SGDM:lupr}
\end{figure}

In our experiments comparing different client optimizers, we had two tunable hyperparameters: the local updates per round $U$ and client learning rate $\eta_l$. First, we look at how the choice of $U$ affects the performance of the federated agent. We present the progression of energy consumption and comfort violation on the evaluation environment for different values of $U$ in figures \ref{sens:Adam:lupr}, \ref{sens:SGD:lupr} and \ref{sens:SGDM:lupr}, for Adam, SGD and SGDM, respectively. The performance for different values of $U$ tends to be quite comparable, for every tested client optimizer. We do not observe any one value of $U$ that consistently outperforms the others, though $U=4$ tends to fall short of the others in terms of energy consumption. We notice that $U=4$ also has slightly worse stability than other values of $U$, both with respect to energy consumption and comfort violation. This is in line with the previous experiments reported in \cite{hagstrom2023using}, where also a centralized agent trained on data pooled from all environments (i.e., having $U=1$) was shown to underperform against the federated agents. 

Considering these results, conclude that the federated agent is robust to the choice of $U$. It is, however, advisable to use larger values, not only because of the worse stability when performing global aggregation after every local update but also because larger values of $U$ mean fewer communication rounds, reducing the communication costs of the federated algorithm. 

\begin{figure}[h]
    \centering
    \begin{subfigure}[b]{0.49\textwidth}
         \centering
         \includegraphics[width=\textwidth]{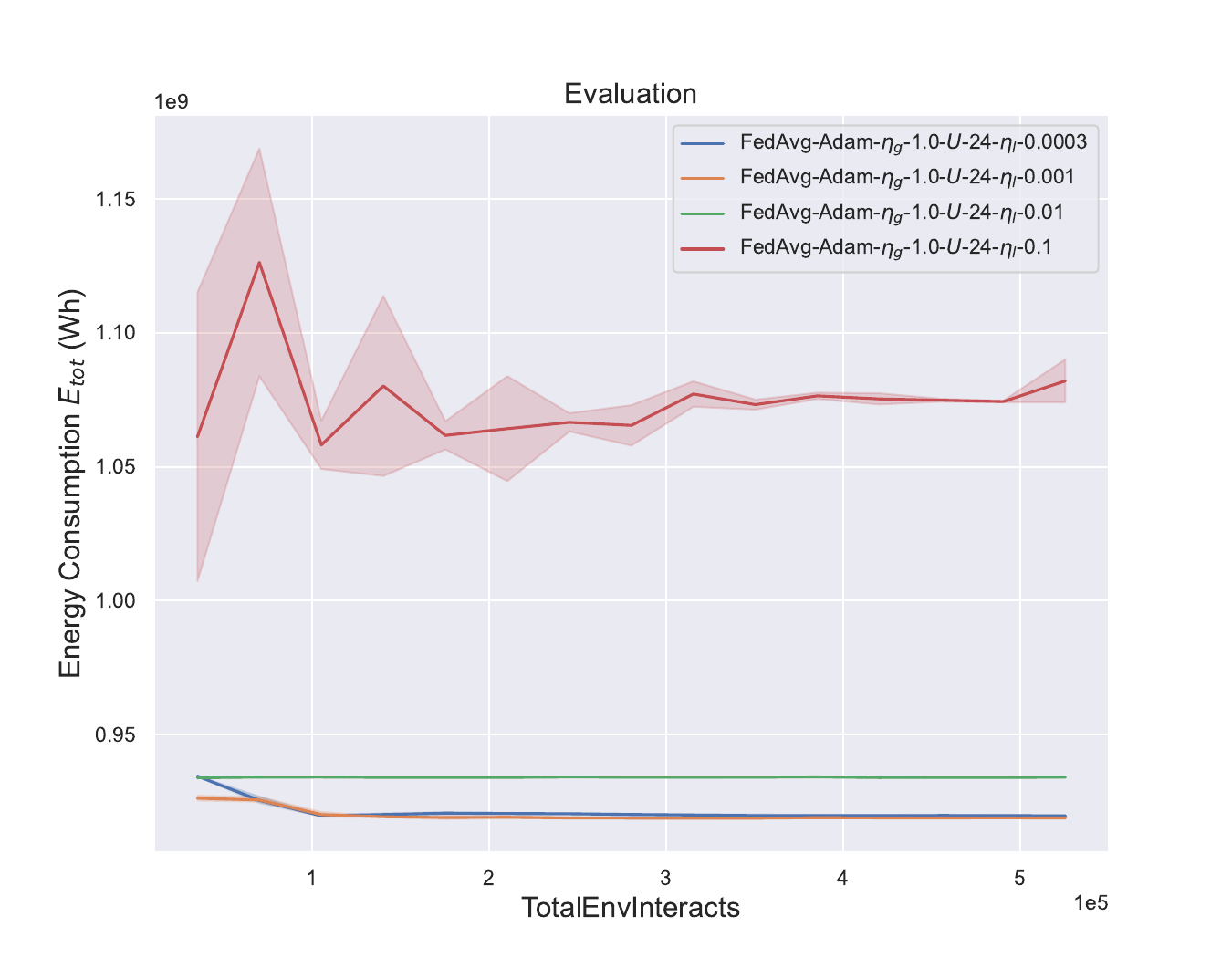}
         \caption{Energy consumption}
         \label{sens:Adam:clr:energy}
    \end{subfigure}
    \hfill
    \begin{subfigure}[b]{0.49\textwidth}
         \centering
         \includegraphics[width=\textwidth]{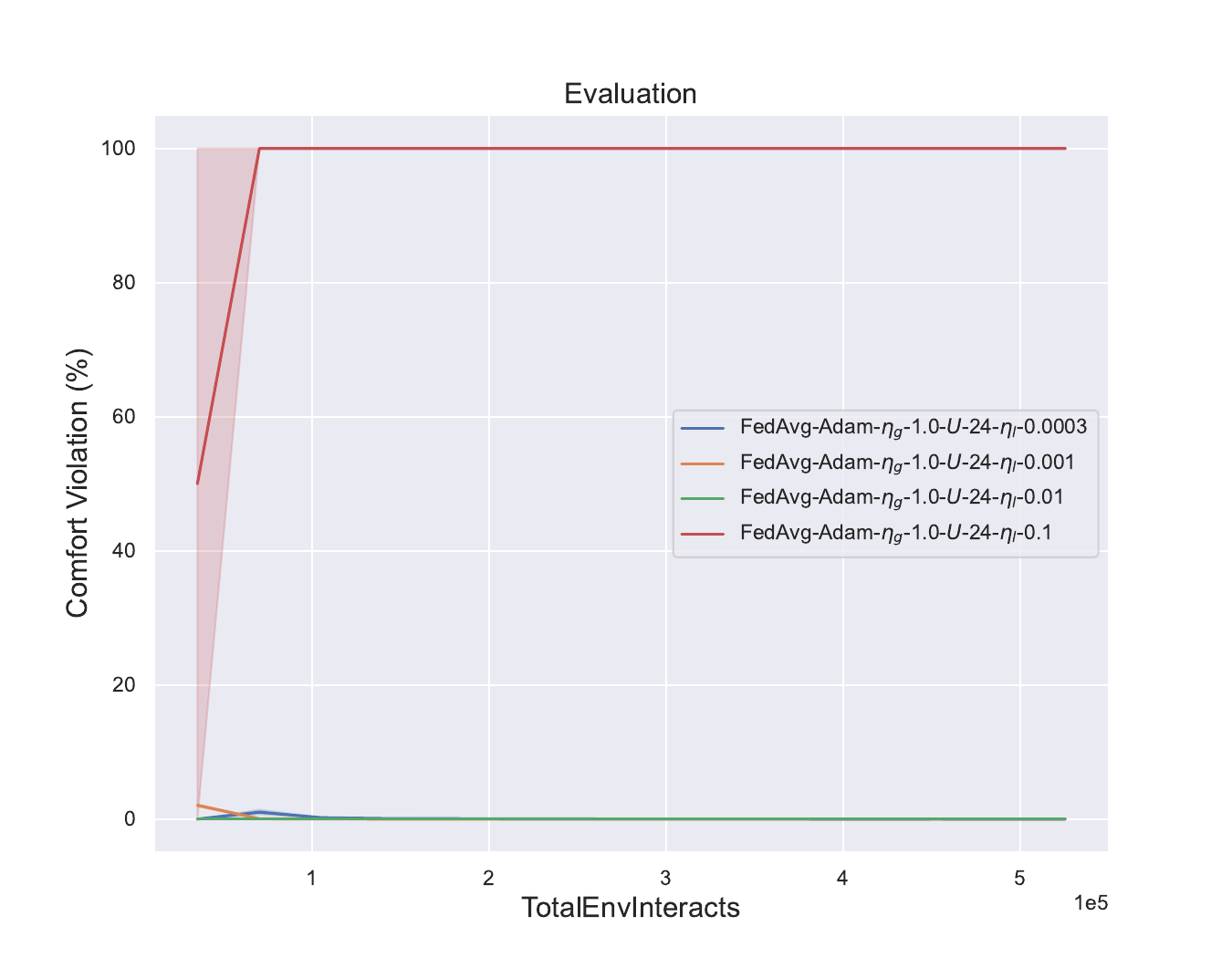}
         \caption{Comfort violation}
         \label{sens:Adam:clr:comfort}
    \end{subfigure}
    \caption{Comparing the performance of FedAvg with Adam as client optimizer on the evaluation environment for different client learning rates $\eta_l$. We fix $U = 24$.}
    \label{sens:Adam:clr}
\end{figure}

\begin{figure}[h]
    \centering
    \begin{subfigure}[b]{0.49\textwidth}
         \centering
         \includegraphics[width=\textwidth]{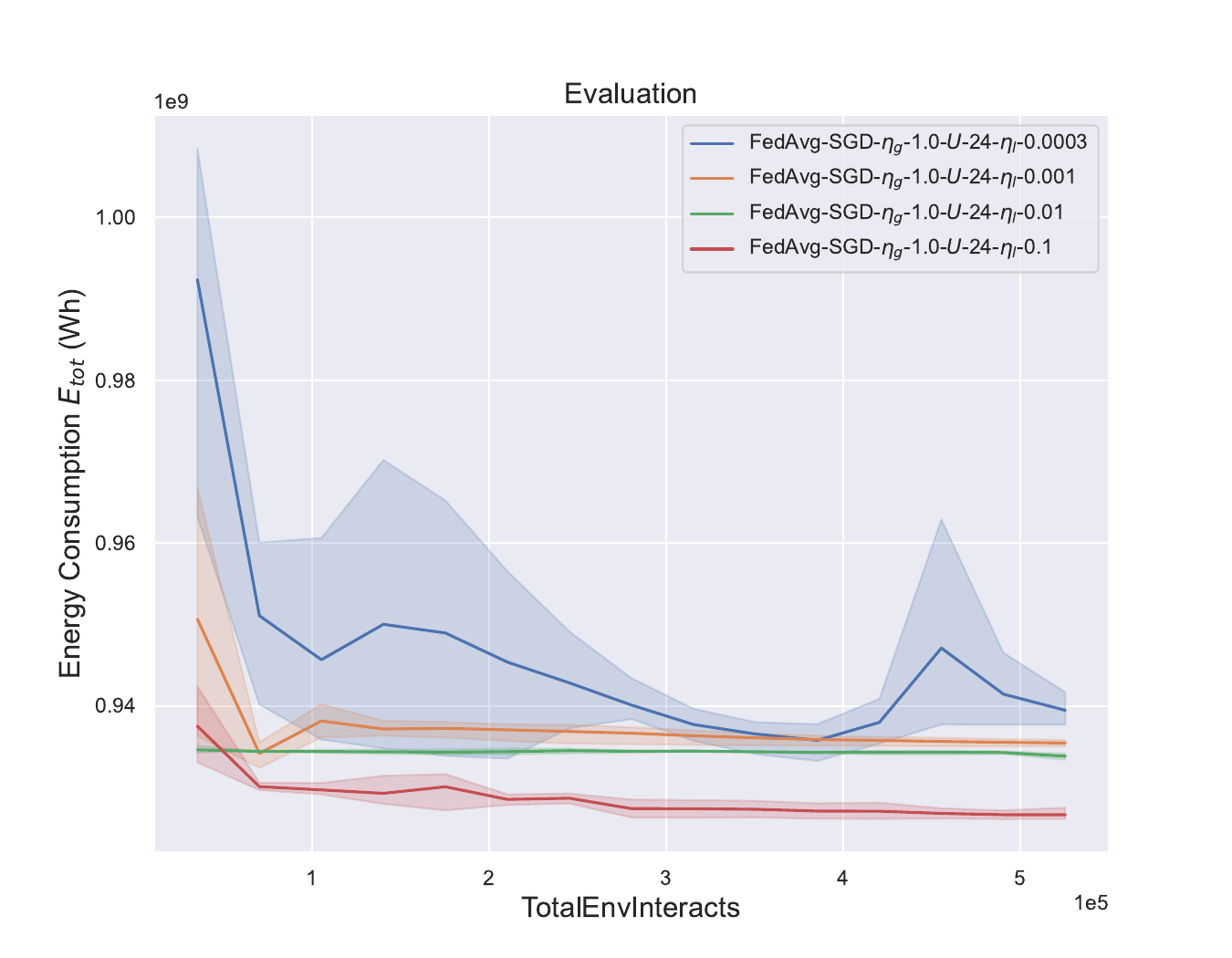}
         \caption{Energy consumption}
         \label{sens:SGD:clr:energy}
    \end{subfigure}
    \hfill
    \begin{subfigure}[b]{0.49\textwidth}
         \centering
         \includegraphics[width=\textwidth]{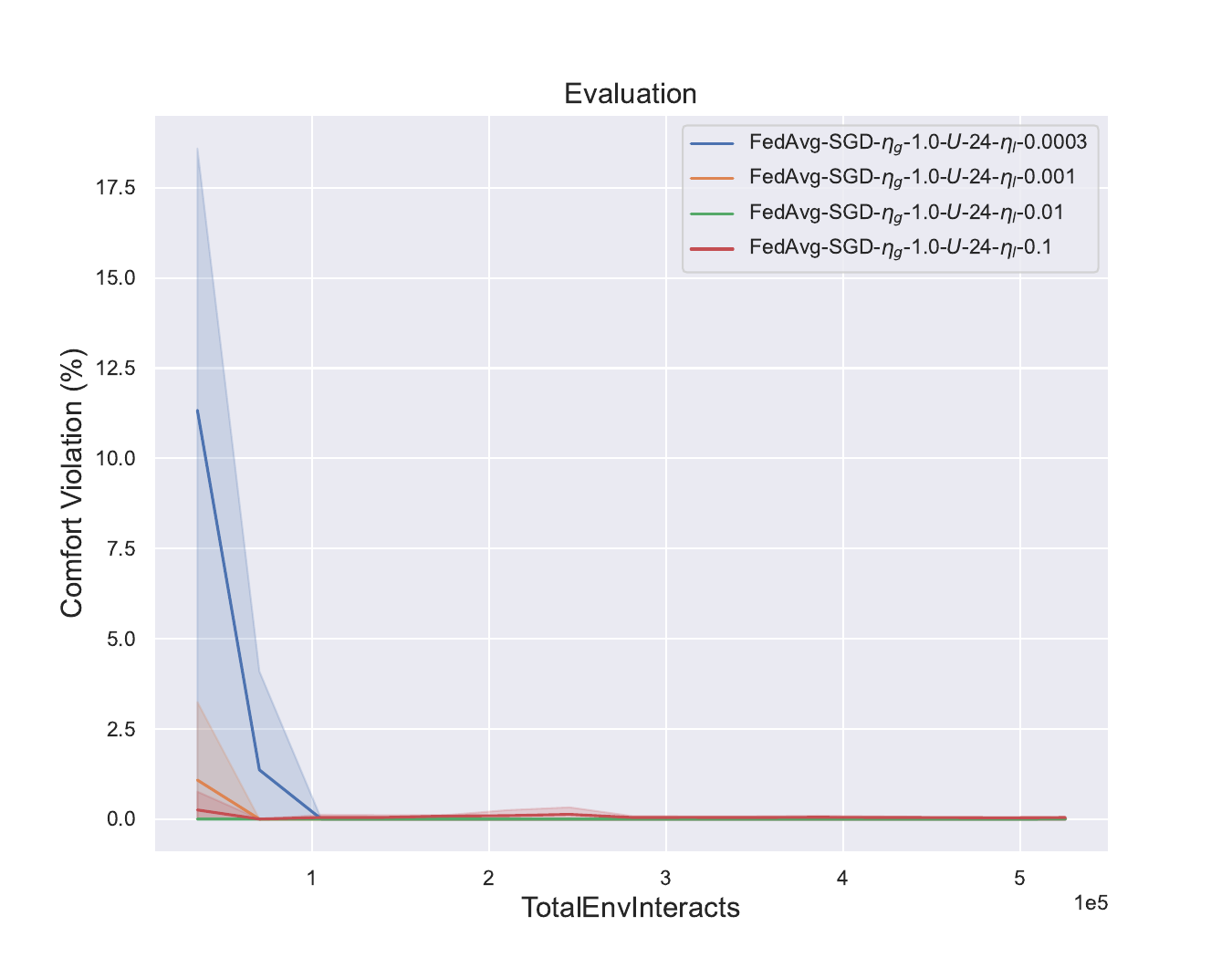}
         \caption{Comfort violation}
         \label{sens:SGD:clr:comfort}
    \end{subfigure}
    \caption{Comparing the performance of FedAvg with SGD as client optimizer on the evaluation environment for different client learning rates $\eta_l$. We fix $U = 24$.}
    \label{sens:SGD:clr}
\end{figure}

\begin{figure}[h]
    \centering
    \begin{subfigure}[b]{0.49\textwidth}
         \centering
         \includegraphics[width=\textwidth]{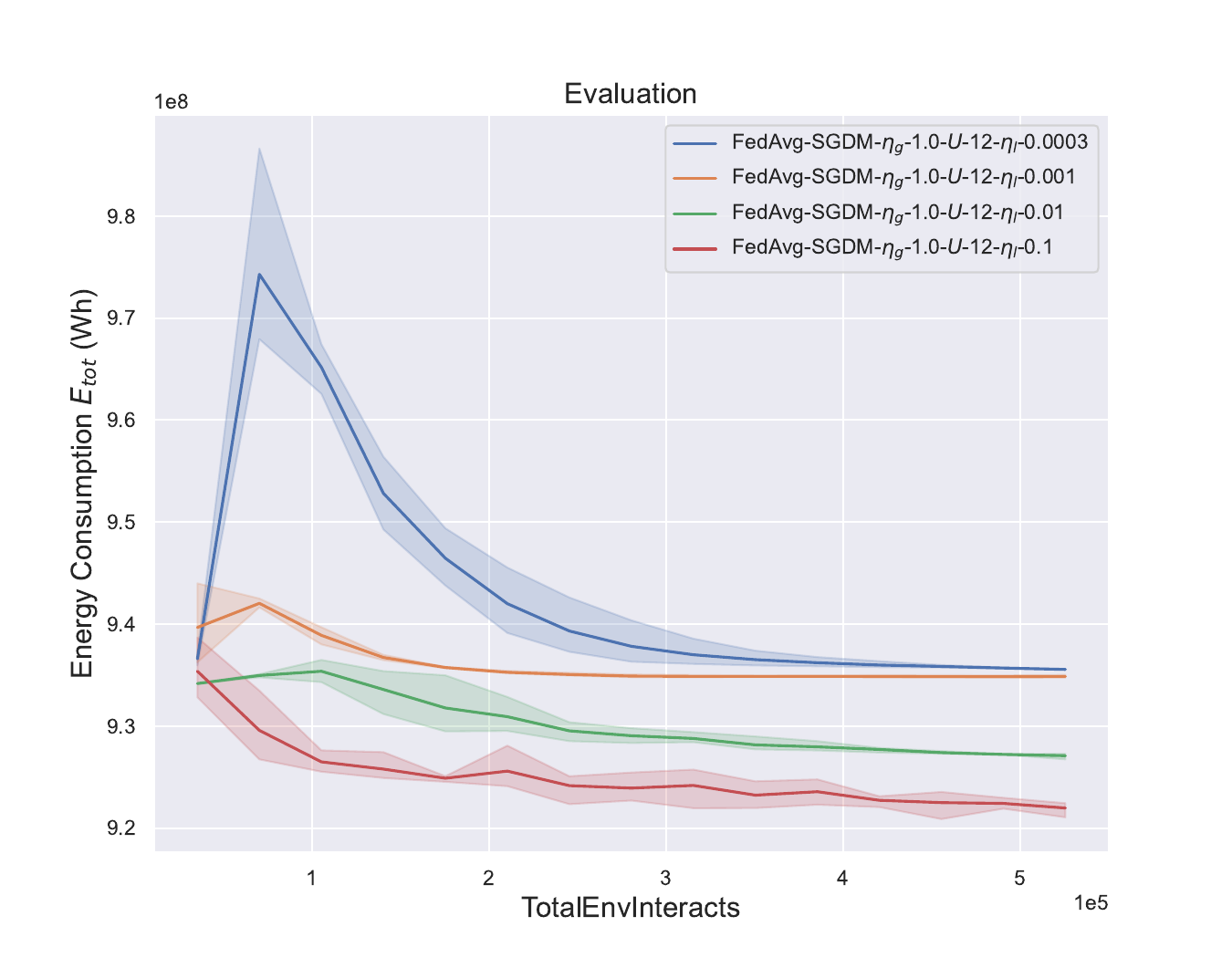}
         \caption{Energy consumption}
         \label{sens:SGDM:clr:energy}
    \end{subfigure}
    \hfill
    \begin{subfigure}[b]{0.49\textwidth}
         \centering
         \includegraphics[width=\textwidth]{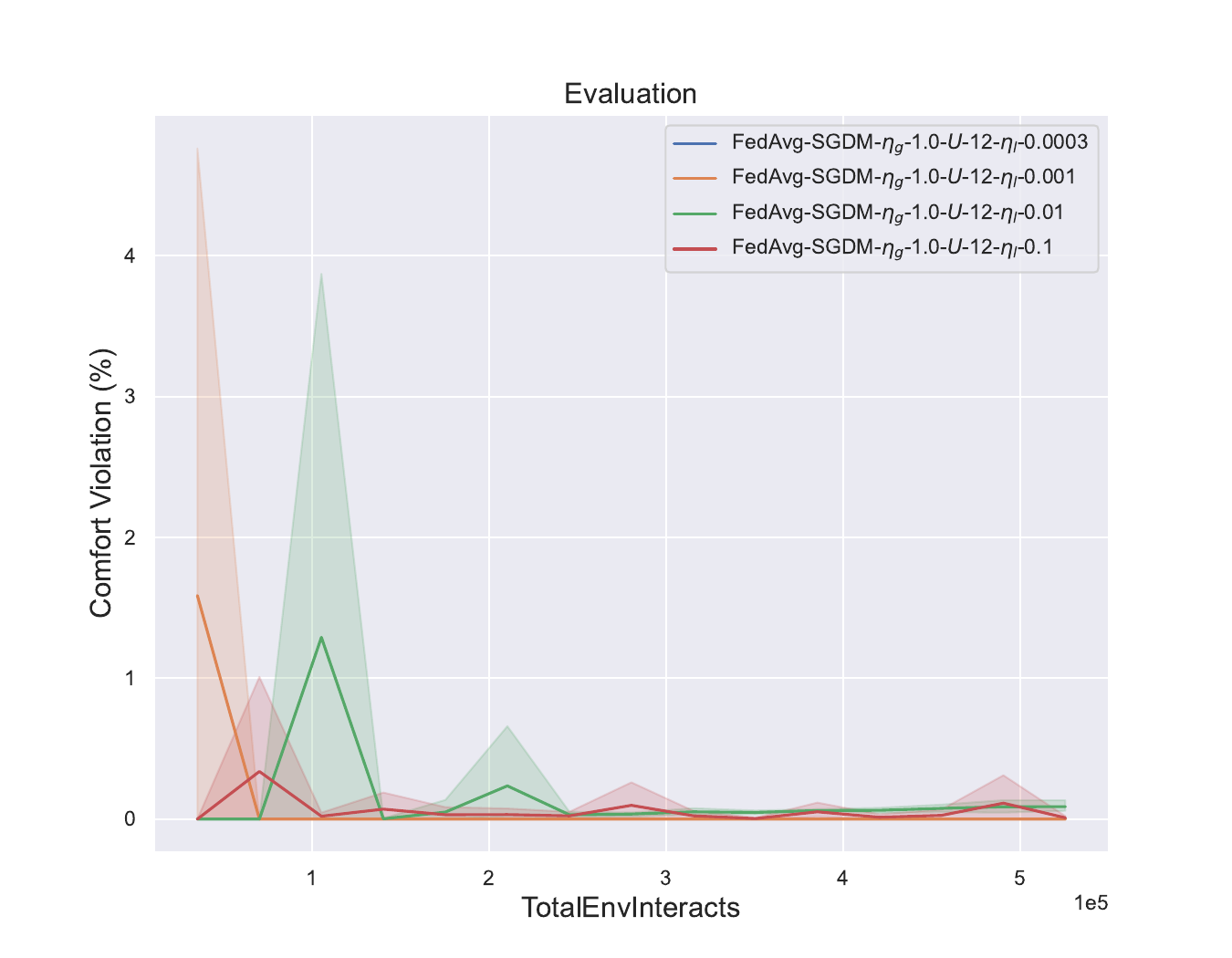}
         \caption{Comfort violation}
         \label{sens:SGDM:clr:comfort}
    \end{subfigure}
    \caption{Comparing the performance of FedAvg with SGDM as client optimizer on the evaluation environment for different client learning rates $\eta_l$. We fix $U = 12$.}
    \label{sens:SGDM:clr}
\end{figure}

Next, we focus on the client learning rate $\eta_l$. We present the progression of energy consumption and comfort violation on the evaluation environment for different values of $\eta_l$ in figures \ref{sens:Adam:clr}, \ref{sens:SGD:clr} and \ref{sens:SGDM:clr}, for Adam, SGD and SGDM, respectively. The performance of the federated agent is sensitive to the client learning rate. In figure \ref{sens:Adam:clr}, we see that higher learning rates can significantly increase the energy consumption of the agent. It can also lead to complete failure in learning a comfortable policy, with $\eta_l = 0.1$ having 100 \% comfort violation. The inverse relationship is true for SGD and SGDM, as can be seen in figures \ref{sens:SGD:clr} and \ref{sens:SGDM:clr}. They tend to achieve lower energy consumption with higher $\eta_l$ and the choice of $\eta_l$ seems to have less of an impact on the comfort violation. 

In section \ref{sec:Results:eval}, we concluded Adam to be the best choice of client optimizer. Based on the observed sensitivity to the client learning rate, it is advisable to use values of $\eta_l \leq 0.001$ for safe performance.

\subsection{Server optimizers} \label{appendix:sensAnal:Server}

\begin{figure}[h]
    \centering
    \begin{subfigure}[b]{0.49\textwidth}
         \centering
         \includegraphics[width=\textwidth]{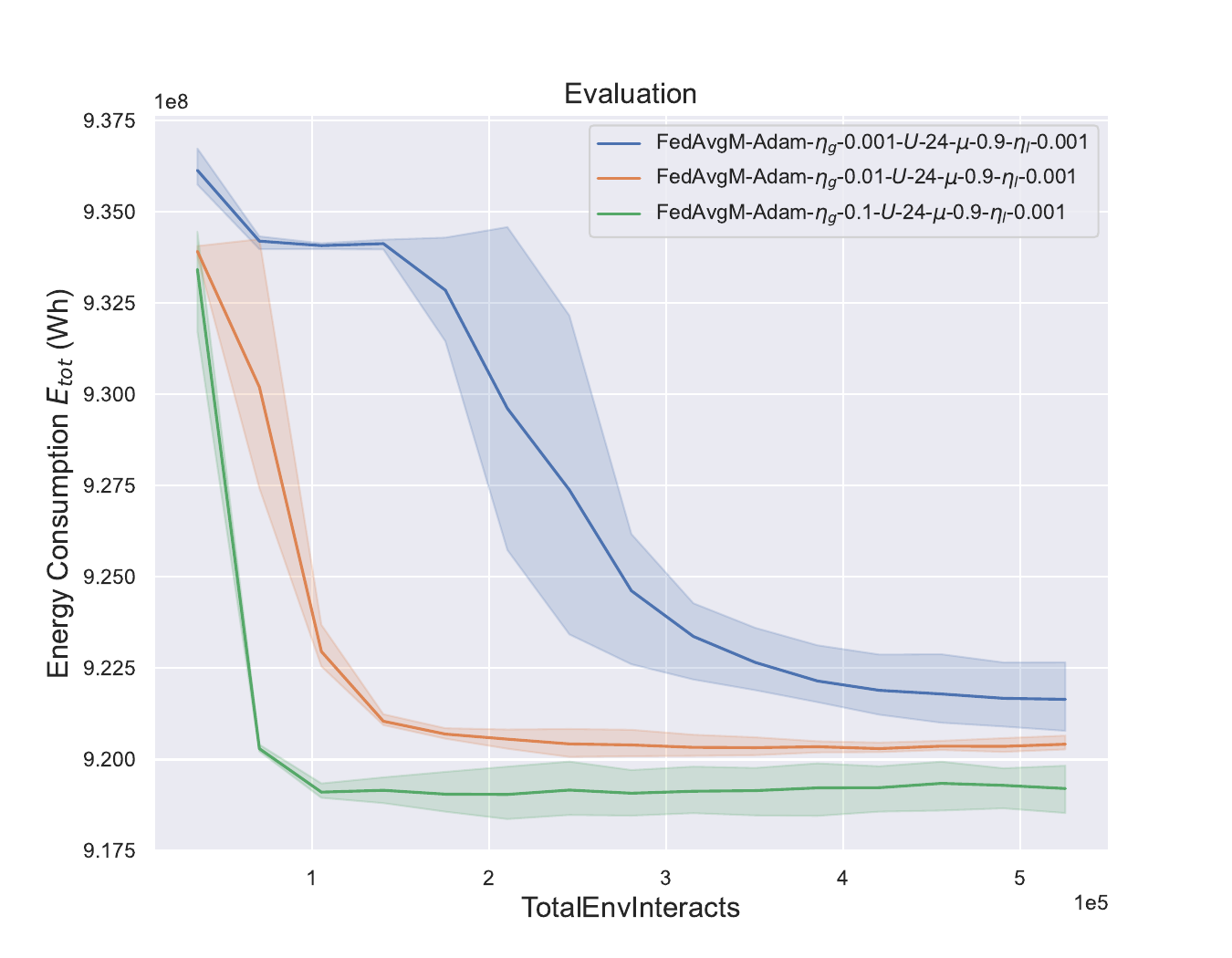}
         \caption{Energy consumption}
         \label{sens:FedAvgM:glr:energy}
    \end{subfigure}
    \hfill
    \begin{subfigure}[b]{0.49\textwidth}
         \centering
         \includegraphics[width=\textwidth]{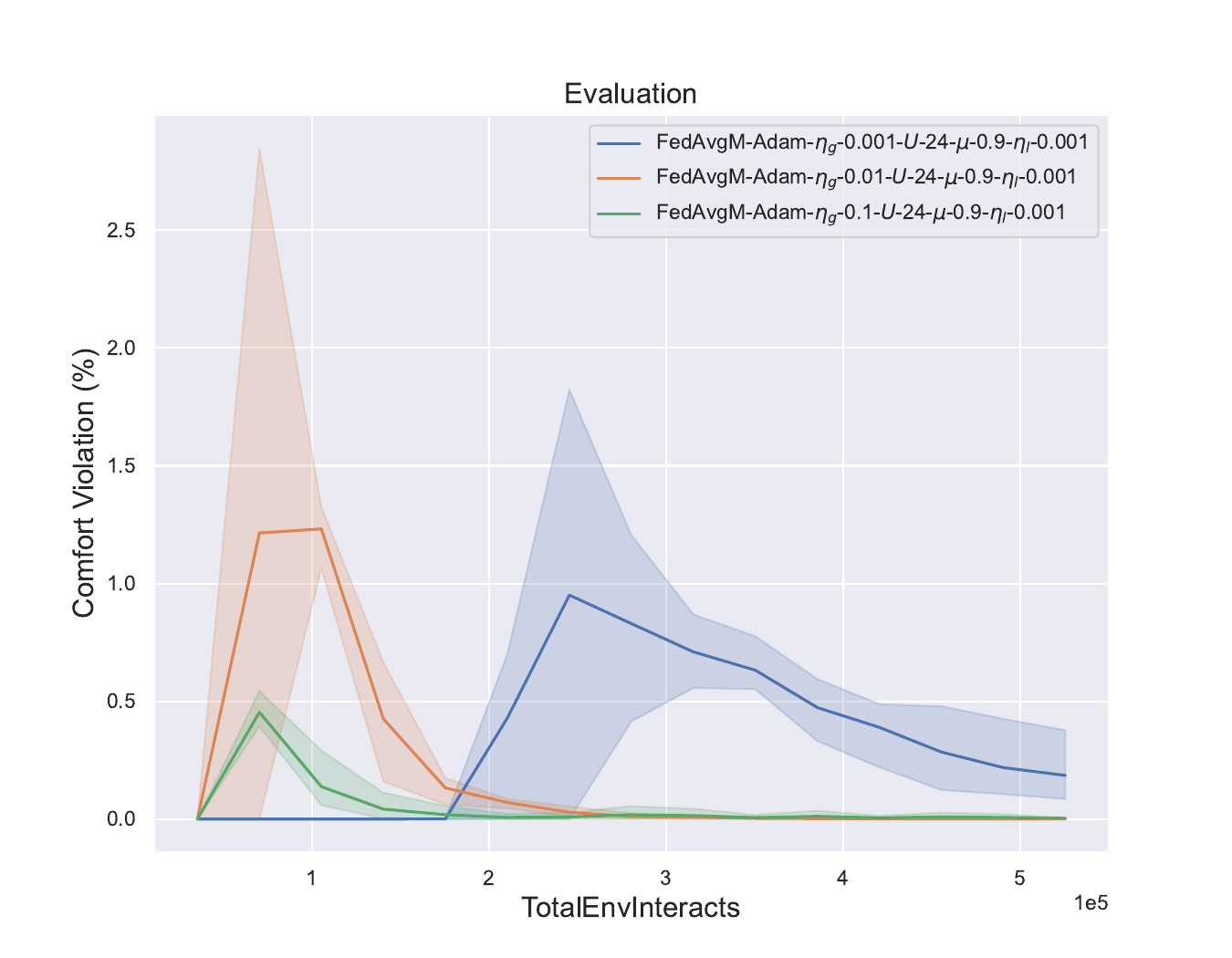}
         \caption{Comfort violation}
         \label{sens:FedAvgM:glr:comfort}
    \end{subfigure}
    \caption{Comparing the performance of FedAvgM on the evaluation environment for different global learning rates $\eta_g$. We fix $U = 24$ and $\beta = 0.9$.}
    \label{sens:FedAvgM:glr}
\end{figure}

\begin{figure}[h]
    \centering
    \begin{subfigure}[b]{0.49\textwidth}
         \centering
         \includegraphics[width=\textwidth]{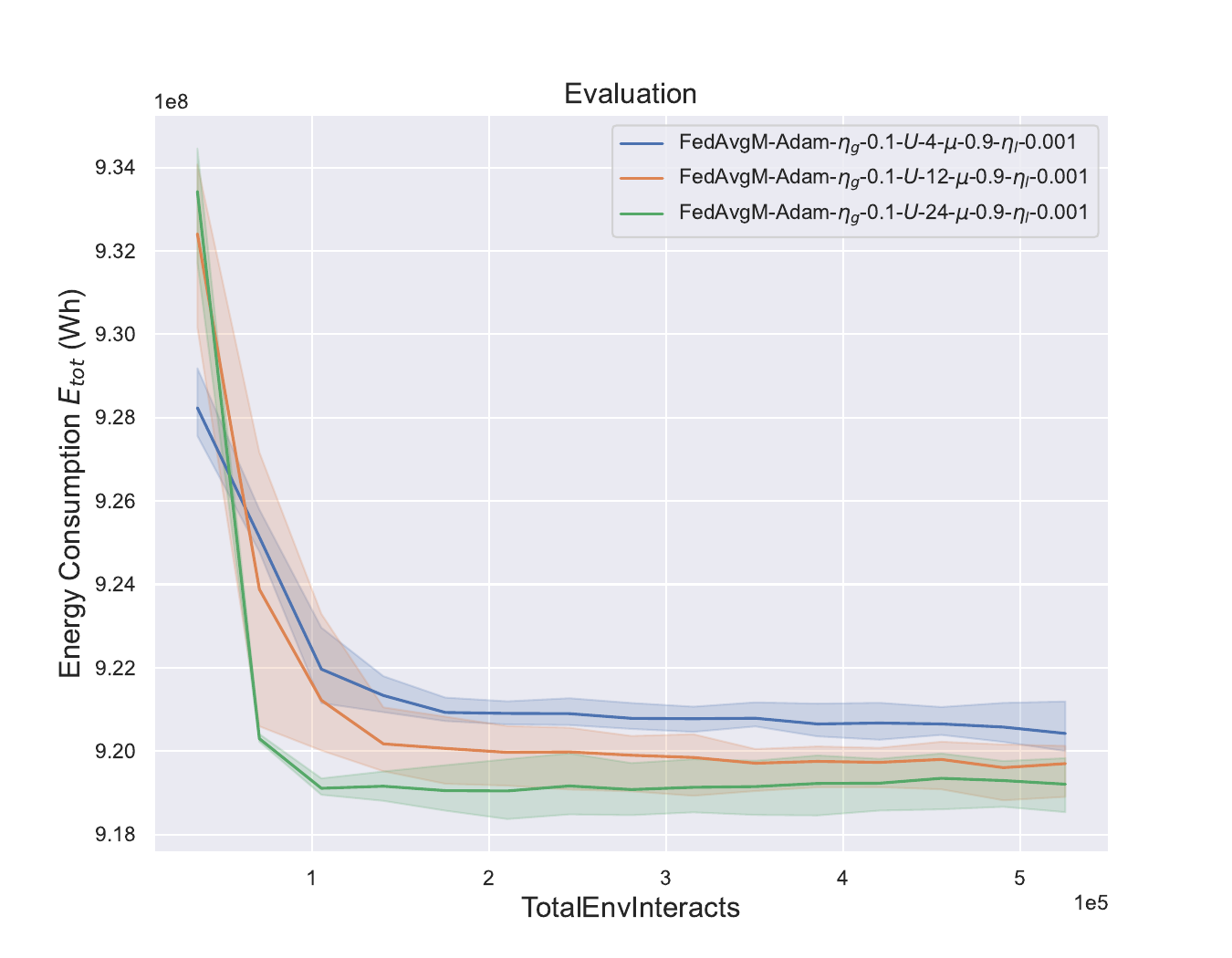}
         \caption{Energy consumption}
         \label{sens:FedAvgM:lupr:energy}
    \end{subfigure}
    \hfill
    \begin{subfigure}[b]{0.49\textwidth}
         \centering
         \includegraphics[width=\textwidth]{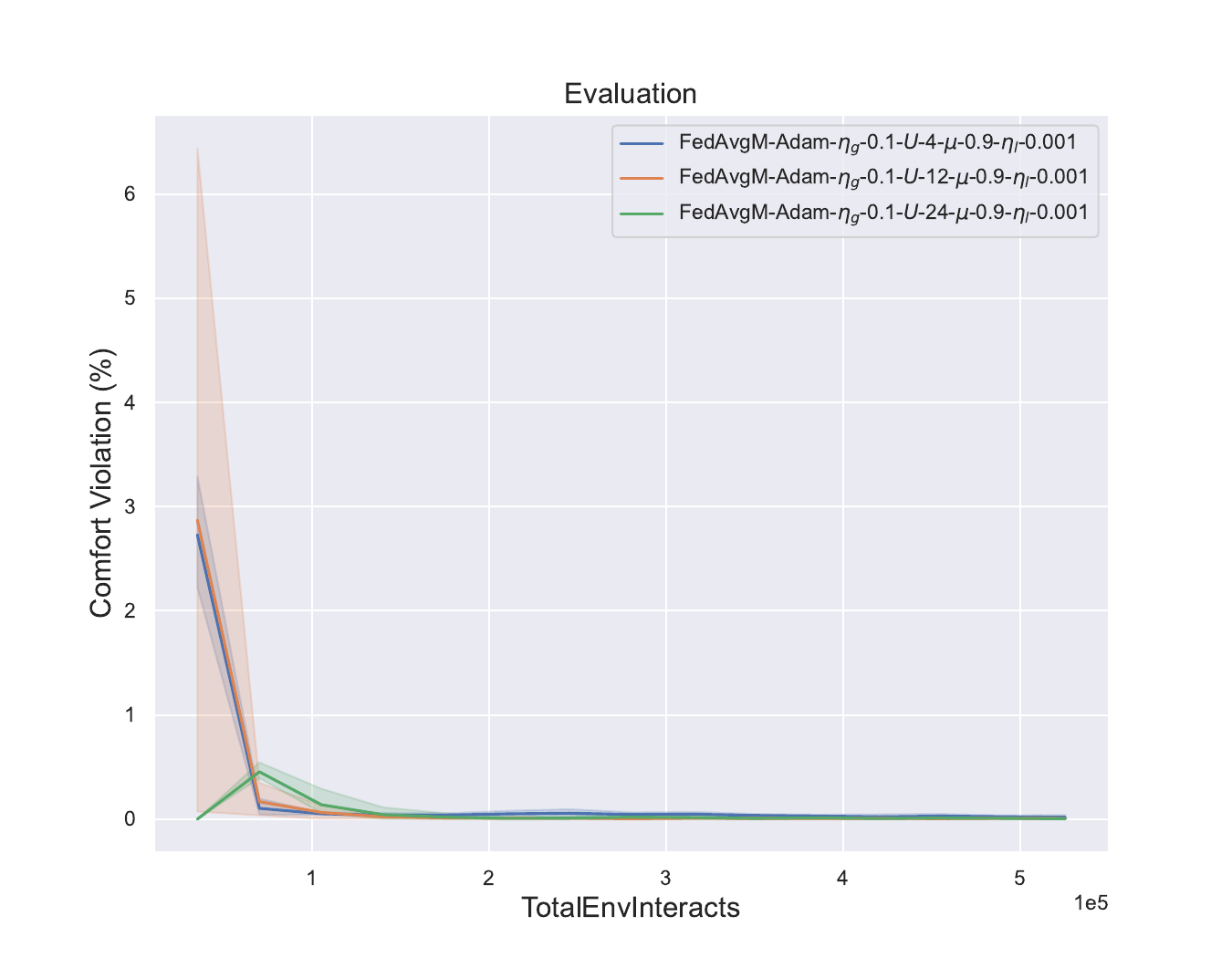}
         \caption{Comfort violation}
         \label{sens:FedAvgM:lupr:comfort}
    \end{subfigure}
    \caption{Comparing the performance of FedAvgM on the evaluation environment for different local updates per round $U$. We fix $\eta_g = 0.1$ and $\beta = 0.9$.}
    \label{sens:FedAvgM:lupr}
\end{figure}

We first consider the sensitivity of the FedAvgM algorithm to its hyperparameters. In figure \ref{sens:FedAvgM:glr}, we present the progression of energy consumption and comfort violation in the evaluation environment for different values of the global learning rate $\eta_g$. We notice a trend of improved performance for larger values of $\eta_g$, both in terms of energy consumption and comfort violation. The global learning rate also affects the learning speed and, to some extent, the learning stability. While increasing the learning rate tends to improve the learning speed and, thus, the performance of FedAvgM, one cannot use arbitrarily large values. We observed in our experiments that setting $\eta_g = 1.0$ tends to lead to exploding gradients, thus leading to an unusable policy. The FedAvgM algorithm is sensitive to the choice of global learning rate, and it needs to be chosen carefully for optimal performance. From figure \ref{sens:FedAvgM:lupr}, we see that FedAvgM is less sensitive to the choice of $U$. FedAvgM displays similar performance, learning speed and learning stability for different values of $U$, though larger values perform slightly better in terms of energy consumption.

\begin{figure}[h]
    \centering
    \begin{subfigure}[b]{0.49\textwidth}
         \centering
         \includegraphics[width=\textwidth]{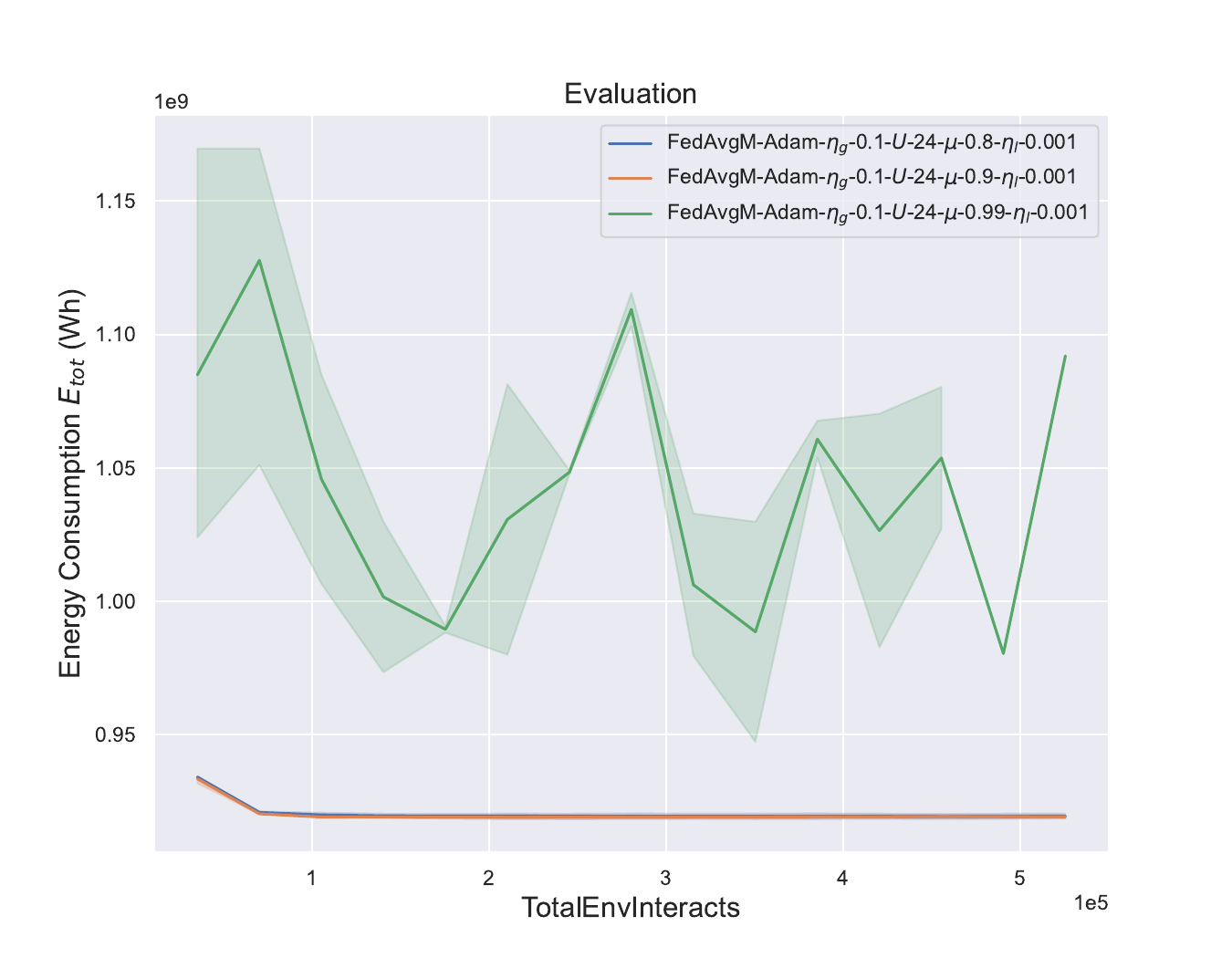}
         \caption{Energy consumption}
         \label{sens:FedAvgM:smom:energy}
    \end{subfigure}
    \hfill
    \begin{subfigure}[b]{0.49\textwidth}
         \centering
         \includegraphics[width=\textwidth]{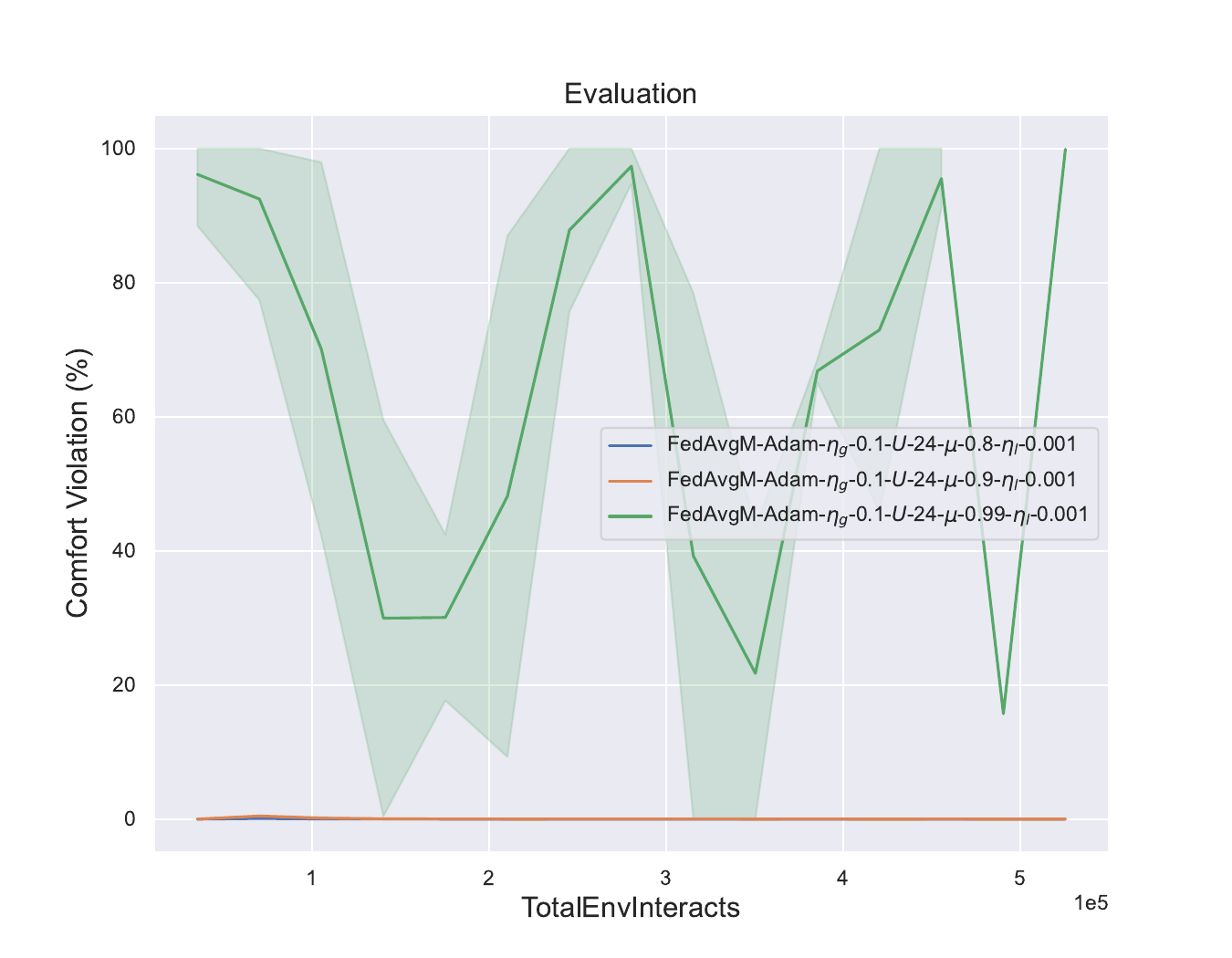}
         \caption{Comfort violation}
         \label{sens:FedAvgM:smom:comfort}
    \end{subfigure}
    \caption{Comparing the performance of FedAvgM on the evaluation environment for different momentums $\beta$. We fix $\eta_g = 0.1$ and $U = 24$.}
    \label{sens:FedAvgM:smom}
\end{figure}

FedAvgM introduces the server momentum parameter $\mu$. The progression of energy consumption and comfort violation for different values of $\mu$ are presented in figure \ref{sens:FedAvgM:smom}. The choice of $\mu$ has a considerable effect on the learning of the model. Too large a value leads to a significant increase in both energy consumption and comfort violation. The learning never converges, and the learning stability is significantly worsened, showing that the FedAvgM is also sensitive to the choice of momentum.

\begin{figure}[h]
    \centering
    \begin{subfigure}[b]{0.49\textwidth}
         \centering
         \includegraphics[width=\textwidth]{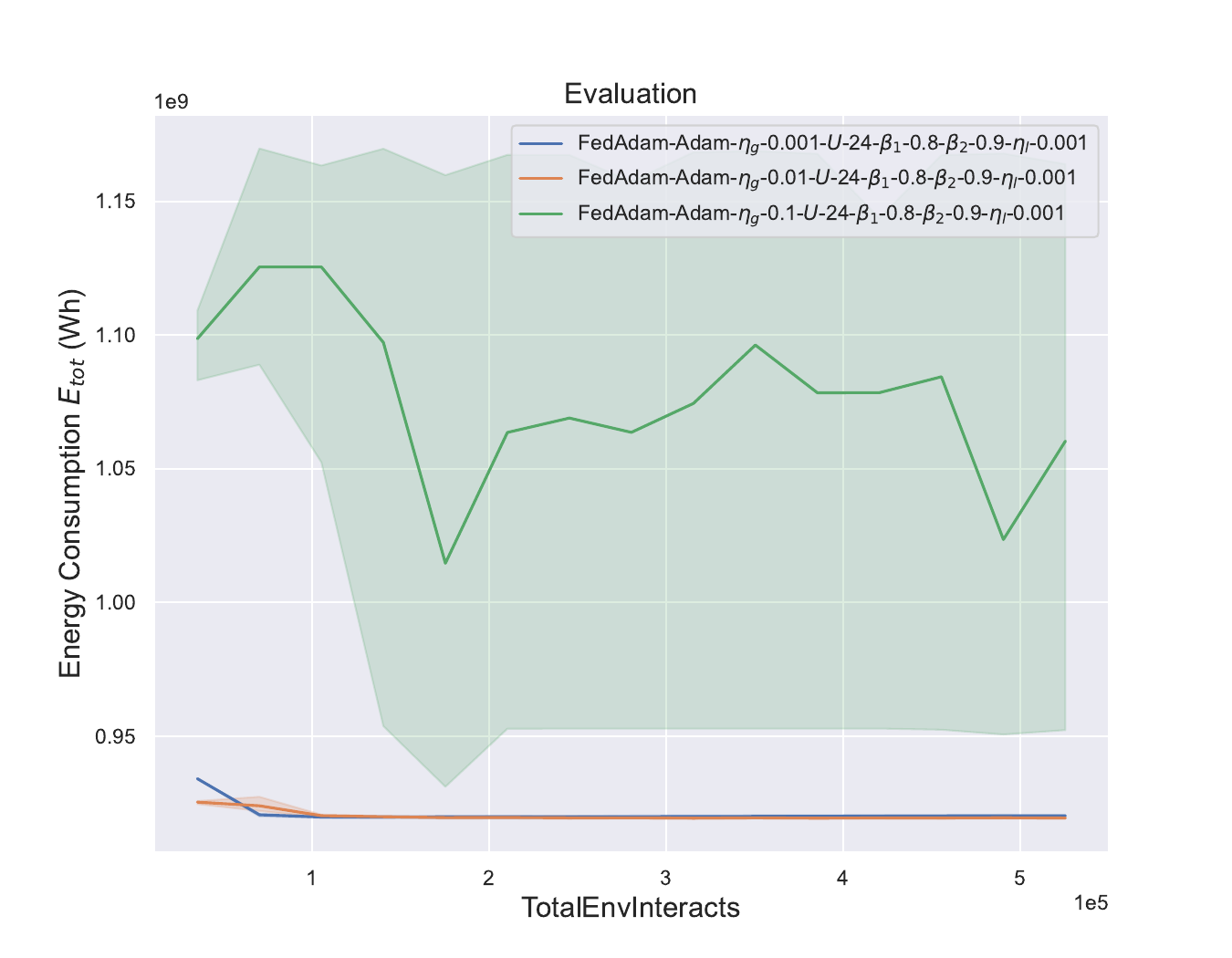}
         \caption{Energy consumption}
         \label{sens:FedAdam:glr:energy}
    \end{subfigure}
    \hfill
    \begin{subfigure}[b]{0.49\textwidth}
         \centering
         \includegraphics[width=\textwidth]{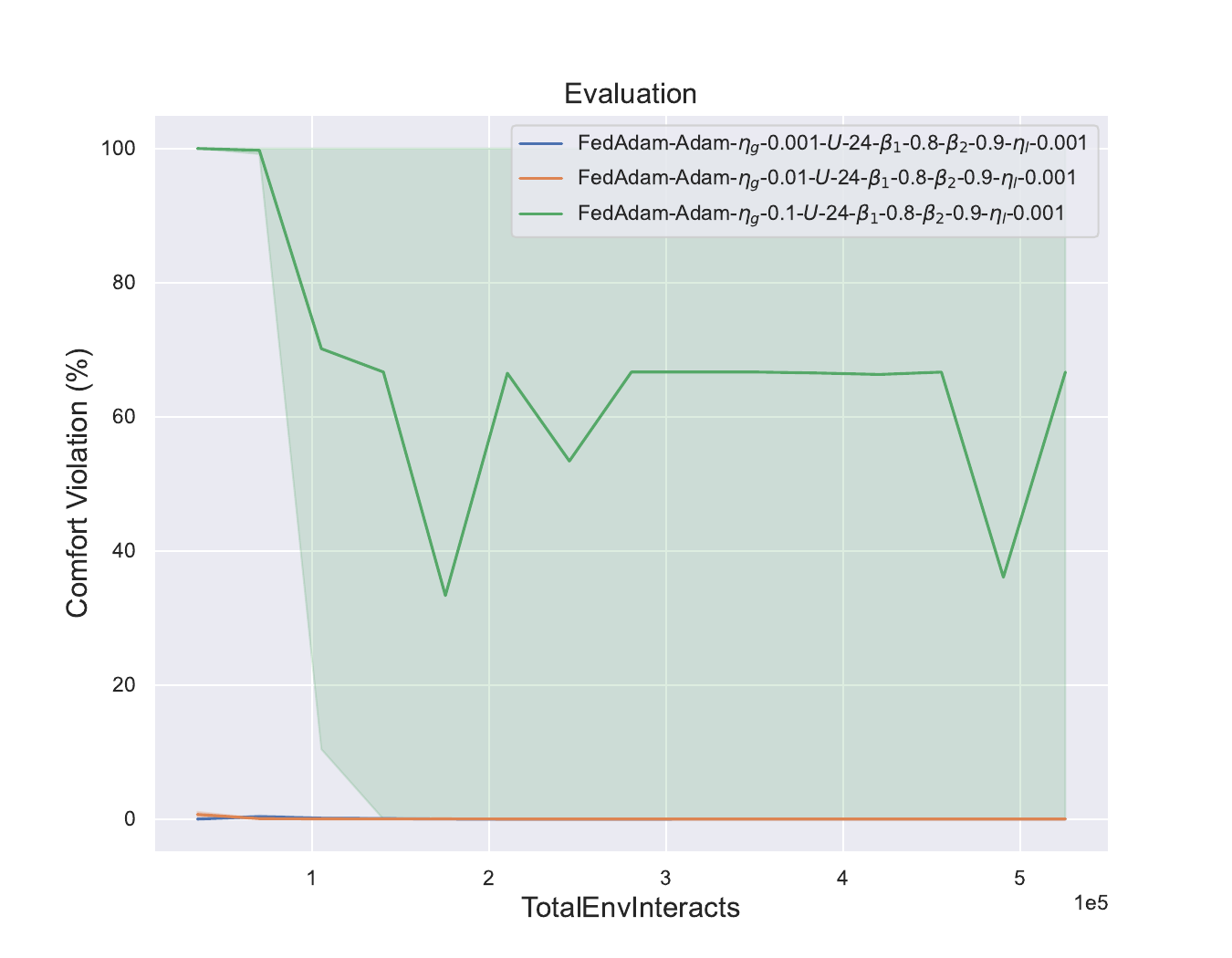}
         \caption{Comfort violation}
         \label{sens:FedAdam:glr:comfort}
    \end{subfigure}
    \caption{Comparing the performance of FedAdam on the evaluation environment for different global learning rates $\eta_g$. We fix $U = 24$, $\beta_1 = 0.8$ and $\beta_2 = 0.9$.}
    \label{sens:FedAdam:glr}
\end{figure}

Similarly to FedAvgM, FedAdam is sensitive to the choice of global learning rate $\eta_g$, as can be seen in figure \ref{sens:FedAdam:glr}. FedAdam, however, performs better with smaller learning rates. Larger learning rates lead to a significant reduction in performance and learning stability, both in terms of energy consumption and comfort violation. Too large a global learning rate can also lead to failure to learn, as we observed that setting $\eta_g = 1.0$ to result in exploding gradients during training.

\begin{figure}[h]
    \centering
    \begin{subfigure}[b]{0.49\textwidth}
         \centering
         \includegraphics[width=\textwidth]{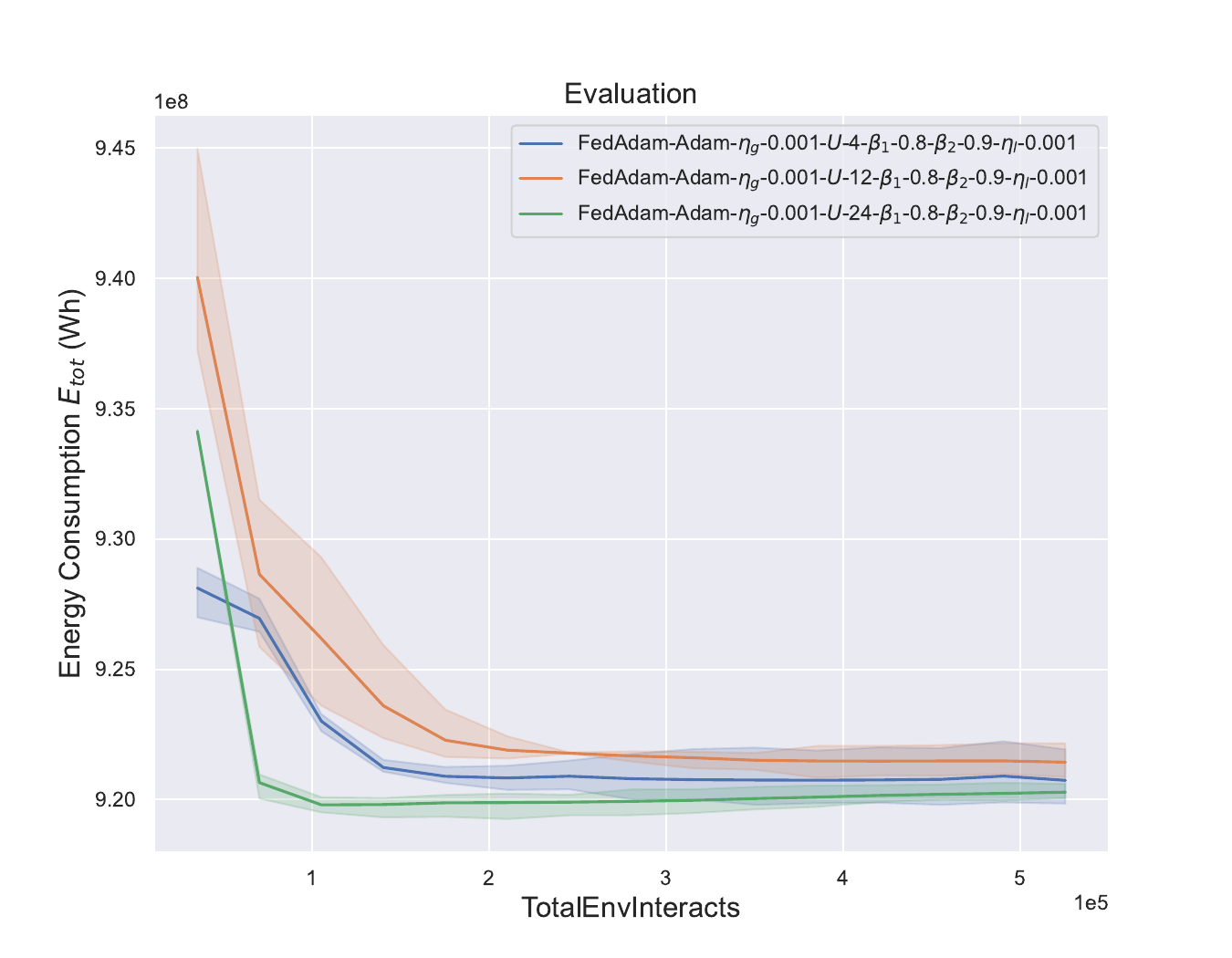}
         \caption{Energy consumption}
         \label{sens:FedAdam:lupr:energy}
    \end{subfigure}
    \hfill
    \begin{subfigure}[b]{0.49\textwidth}
         \centering
         \includegraphics[width=\textwidth]{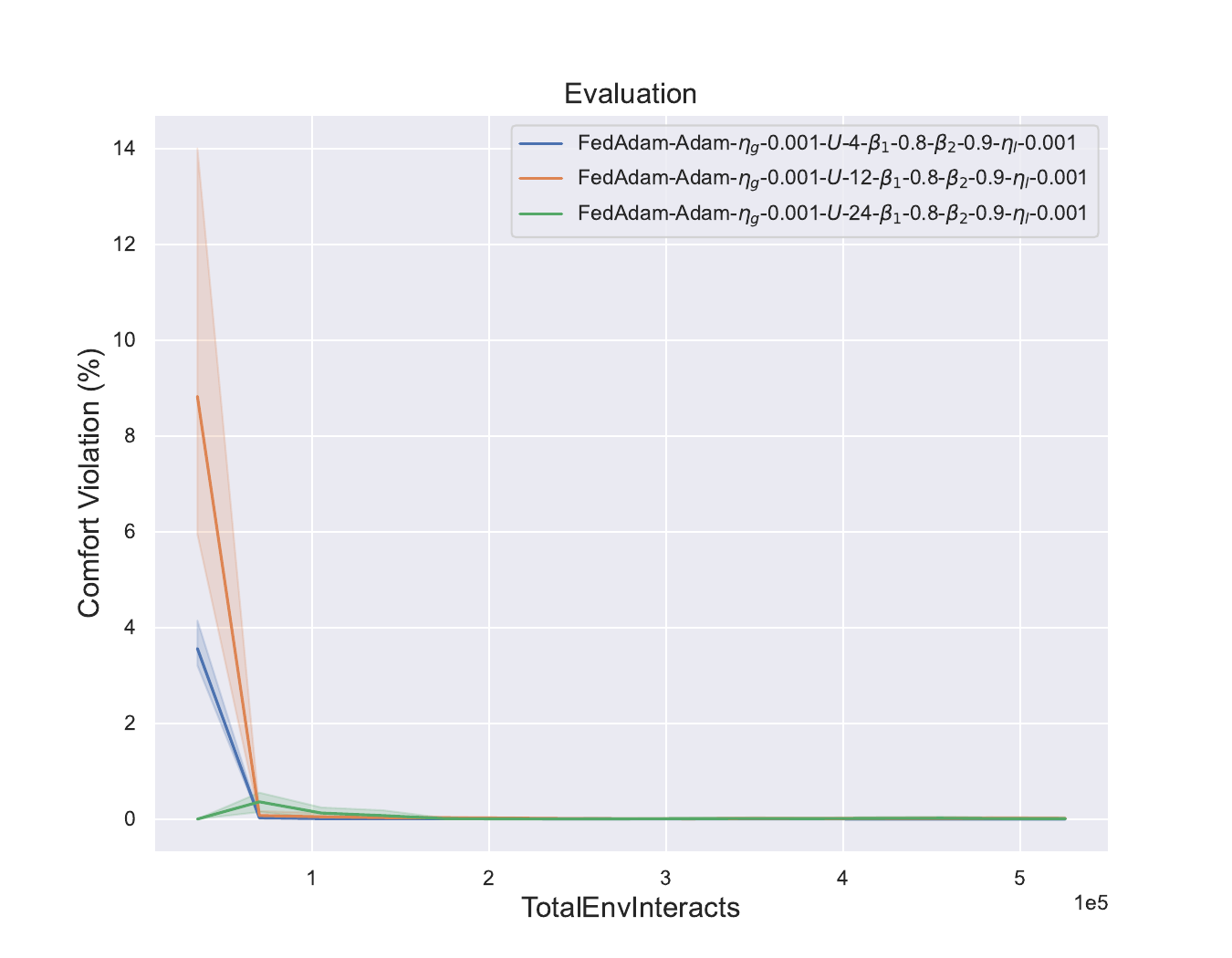}
         \caption{Comfort violation}
         \label{sens:FedAdam:lupr:comfort}
    \end{subfigure}
    \caption{Comparing the performance of FedAdam on the evaluation environment for different local updates per round $U$. We fix $\eta_g = 0.001$, $\beta_1 = 0.8$ and $\beta_2 = 0.9$.}
    \label{sens:FedAdam:lupr}
\end{figure}

In figure \ref{sens:FedAdam:lupr}, we present the learning curves for different values of $U$. As with both FedAvg and FedAvgM, the performance, learning speed and stability are comparable for all tested values of $U$, and larger values display slightly improved energy consumption.

\begin{figure}[h]
    \centering
    \begin{subfigure}[b]{0.49\textwidth}
         \centering
         \includegraphics[width=\textwidth]{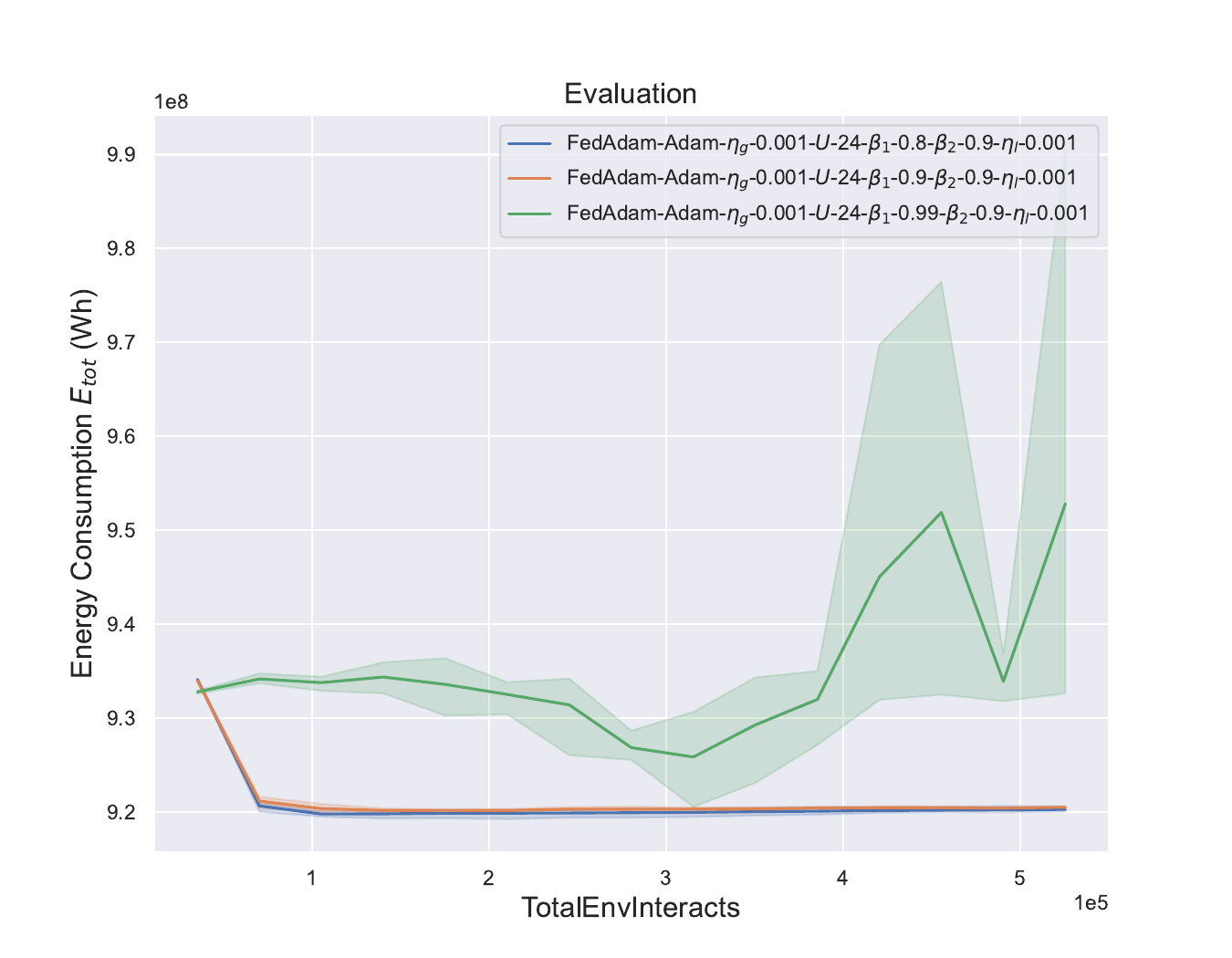}
         \caption{Energy consumption}
         \label{sens:FedAdam:beta_1:energy}
    \end{subfigure}
    \hfill
    \begin{subfigure}[b]{0.49\textwidth}
         \centering
         \includegraphics[width=\textwidth]{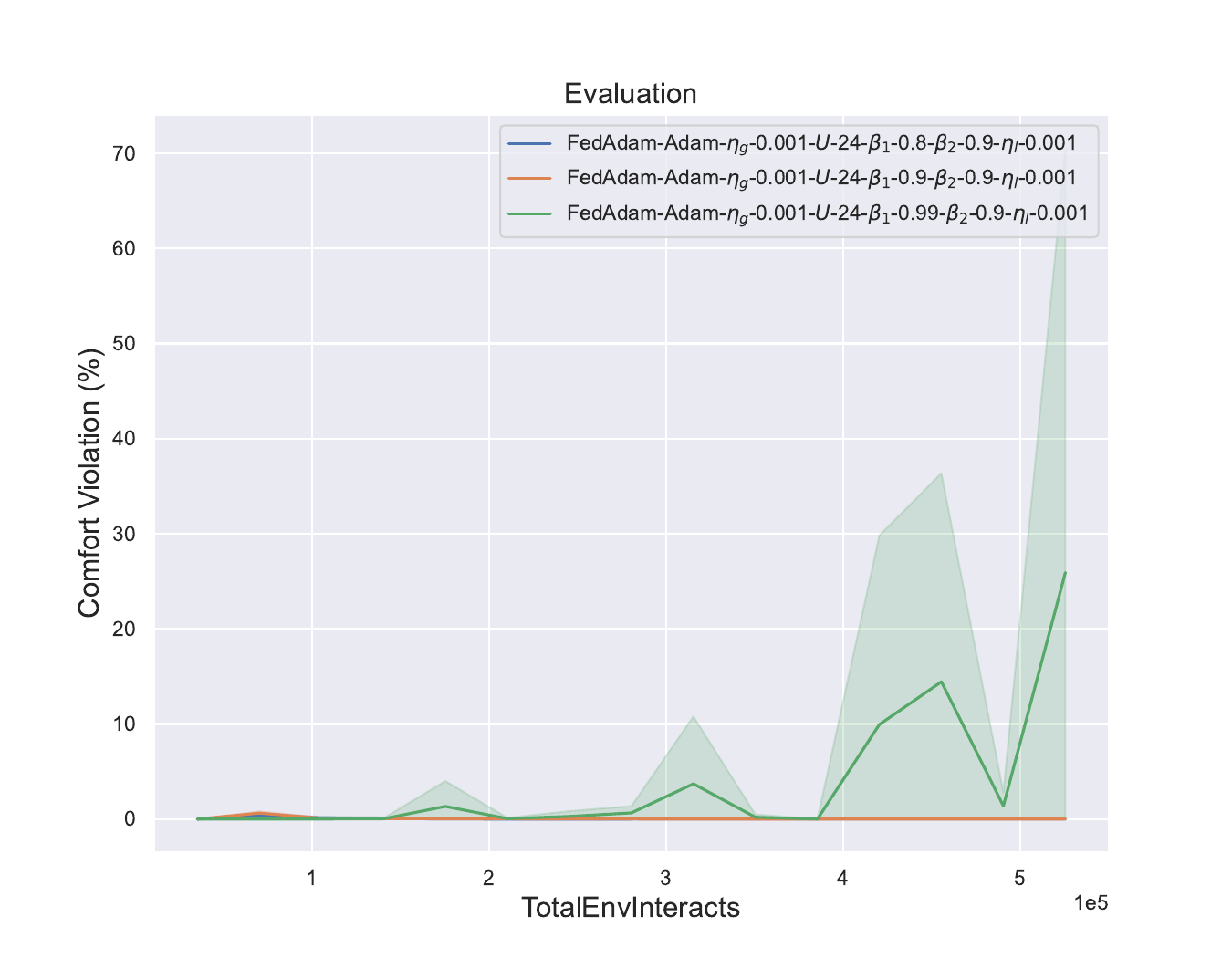}
         \caption{Comfort violation}
         \label{sens:FedAdam:beta_1:comfort}
    \end{subfigure}
    \caption{Comparing the performance of FedAdam on the evaluation environment for different $\beta_1$. We fix $\eta_g = 0.001$, $U = 24$ and $\beta_2 = 0.9$.}
    \label{sens:FedAdam:beta_1}
\end{figure}

\begin{figure}[h]
    \centering
    \begin{subfigure}[b]{0.49\textwidth}
         \centering
         \includegraphics[width=\textwidth]{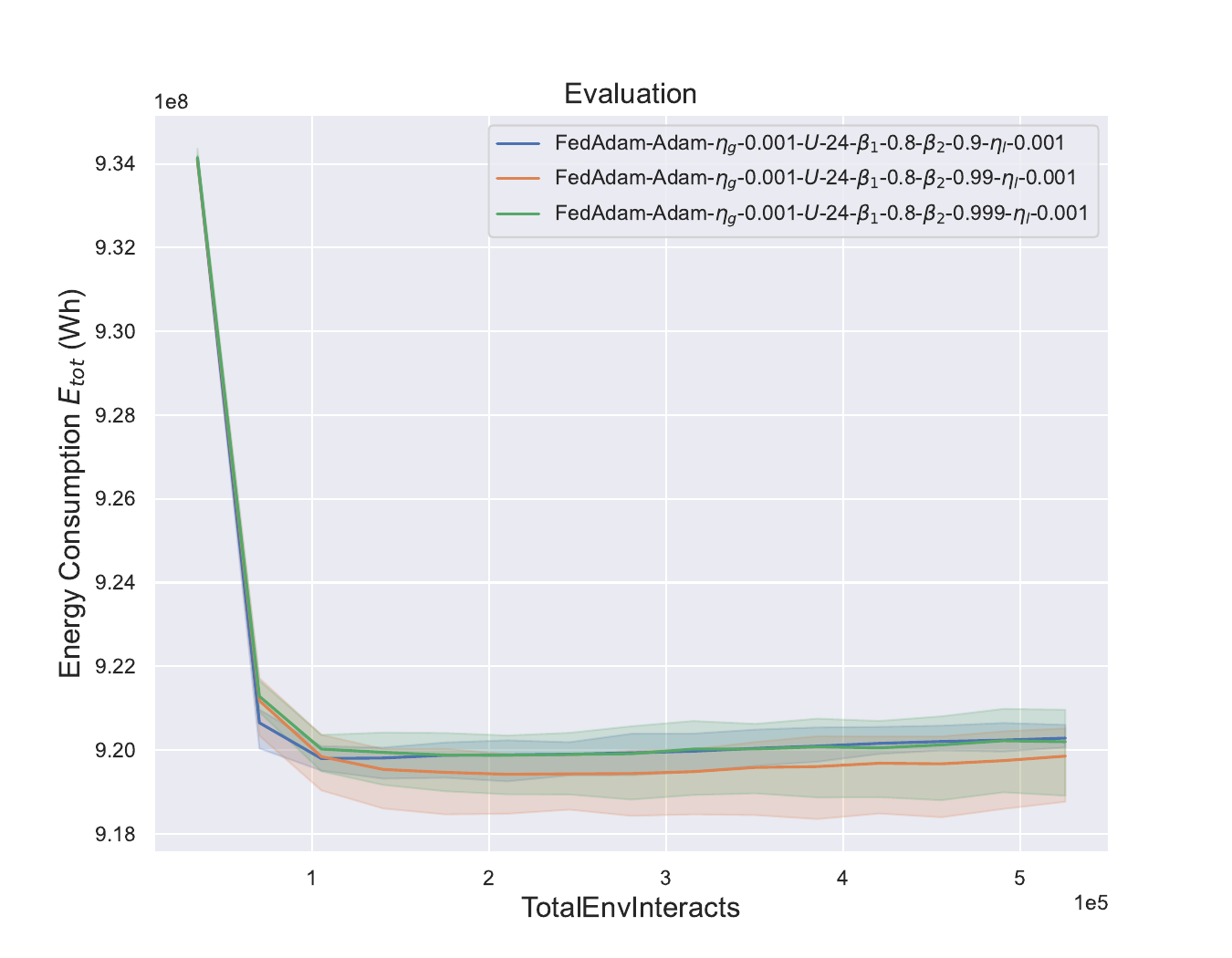}
         \caption{Energy consumption}
         \label{sens:Adam:beta_2:energy}
    \end{subfigure}
    \hfill
    \begin{subfigure}[b]{0.49\textwidth}
         \centering
         \includegraphics[width=\textwidth]{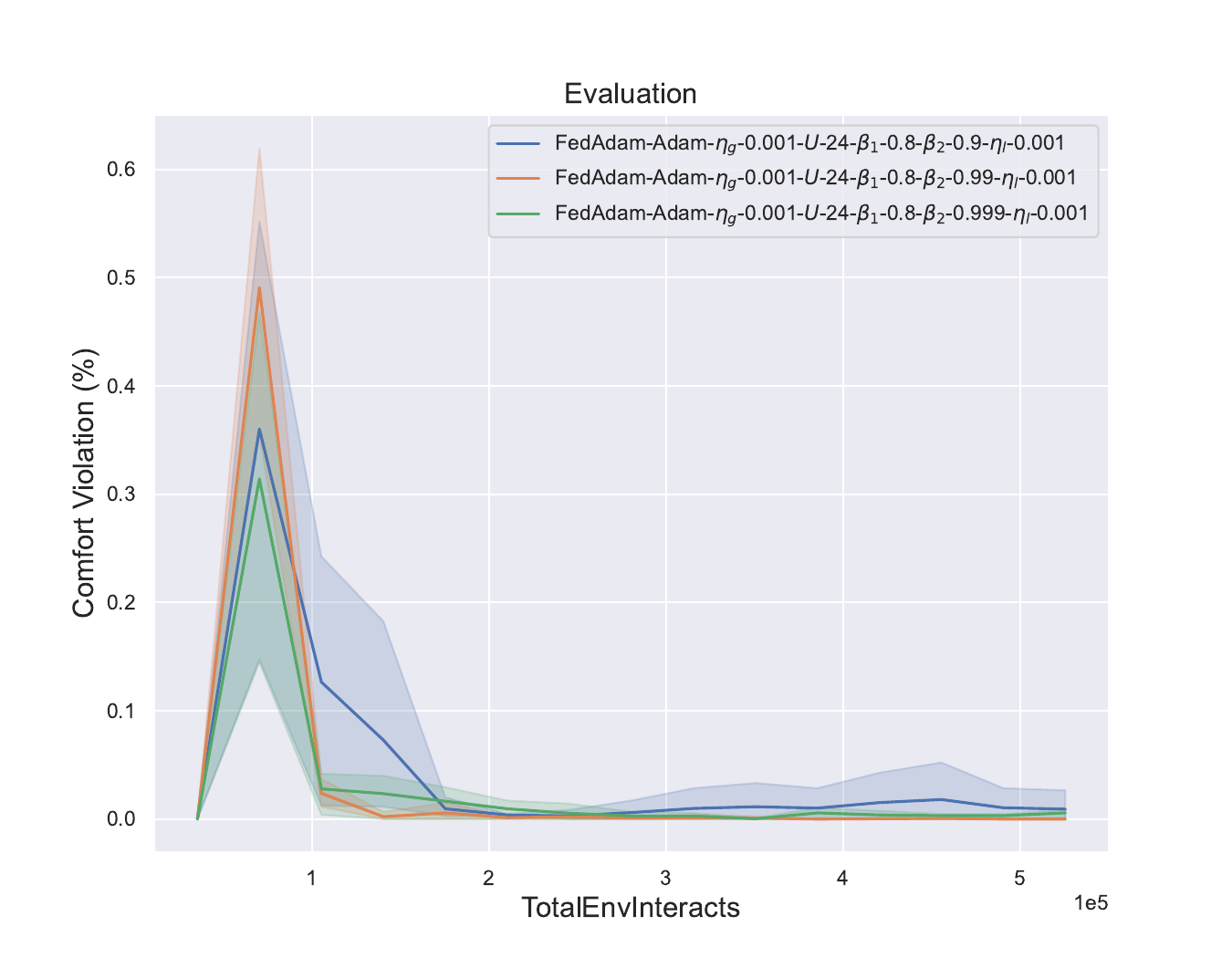}
         \caption{Comfort violation}
         \label{sens:FedAdam:beta_2:comfort}
    \end{subfigure}
    \caption{Comparing the performance of FedAdam on the evaluation environment for different $\beta_2$. We fix $\eta_g = 0.001$, $U = 24$ and $\beta_1 = 0.8$.}
    \label{sens:FedAdam:beta_2}
\end{figure}

FedAdam has two adjustable moment parameters $\beta_1$ and $\beta_2$. FedAdam seems to be more sensitive to the choice of $\beta_1$ than the choice of $\beta_2$. In figure \ref{sens:FedAdam:beta_1}, we see that too large a value of $\beta_1$ leads to a significant degradation in performance and learning stability. $\beta_2$ seems significantly more robust, with all tested values having comparable performance, learning speed and stability in terms of both energy consumption and comfort violation, as can be seen in figure \ref{sens:FedAdam:beta_2}.

\section{Additional plots}
\label{appendix:AddPlots}

\begin{figure}[h]
    \centering
    \begin{subfigure}[b]{0.49\textwidth}
         \centering
         \includegraphics[width=\textwidth]{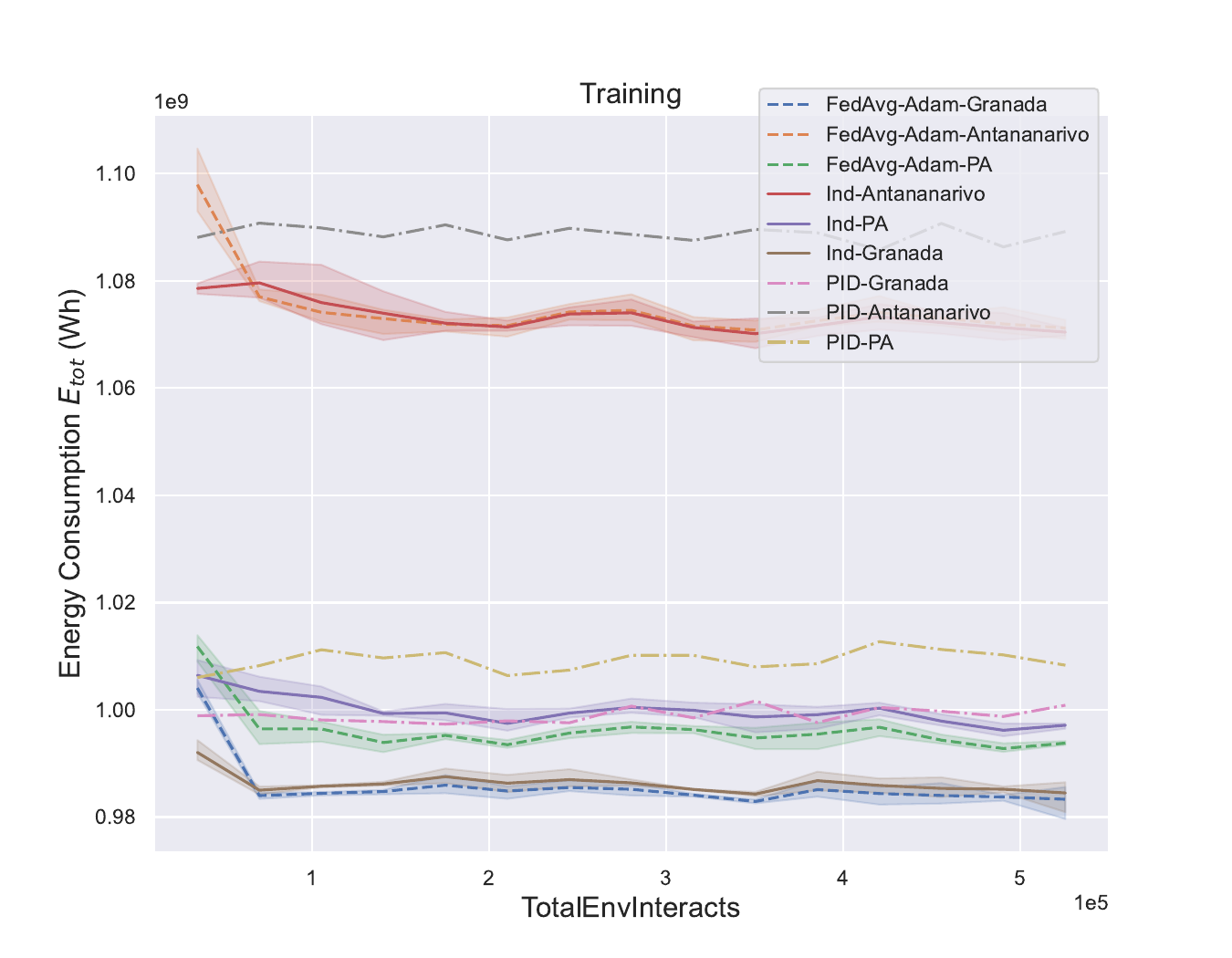}
         \caption{Energy consumption}
         \label{train:fedavg:ind:GraAntPA:energy}
    \end{subfigure}
    \hfill
    \begin{subfigure}[b]{0.49\textwidth}
         \centering
         \includegraphics[width=\textwidth]{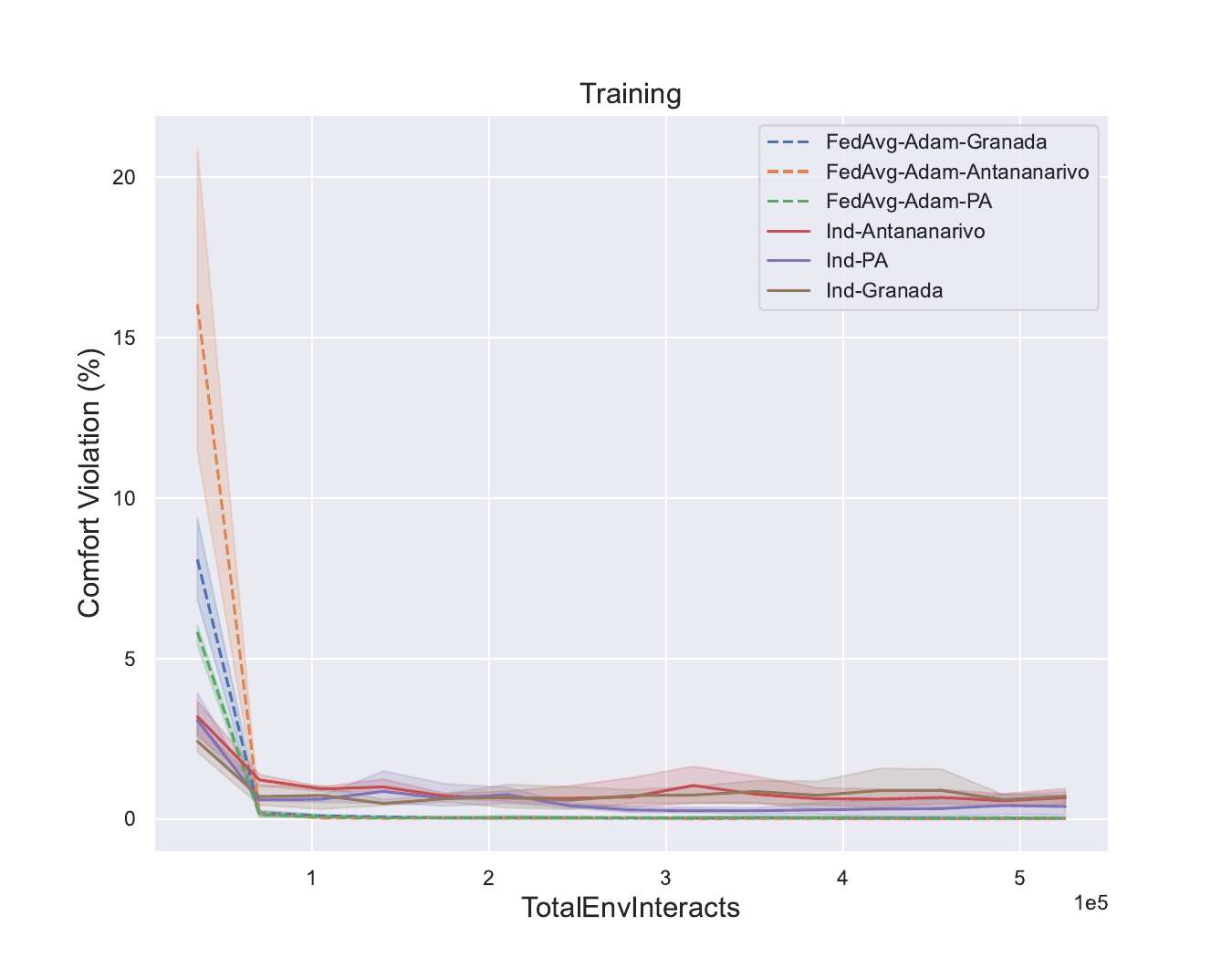}
         \caption{Comfort violation}
         \label{train:fedavg:ind:GraAntPA:comfort}
    \end{subfigure}
    \caption{Progression of the energy consumption and comfort violation of FedAvg and independent agents on training environments Granada, Antananarivo and PA.}
    \label{train:fedavg:ind:GraAntPA}
\end{figure}

\begin{figure}[h]
    \centering
    \begin{subfigure}[b]{0.49\textwidth}
         \centering
         \includegraphics[width=\textwidth]{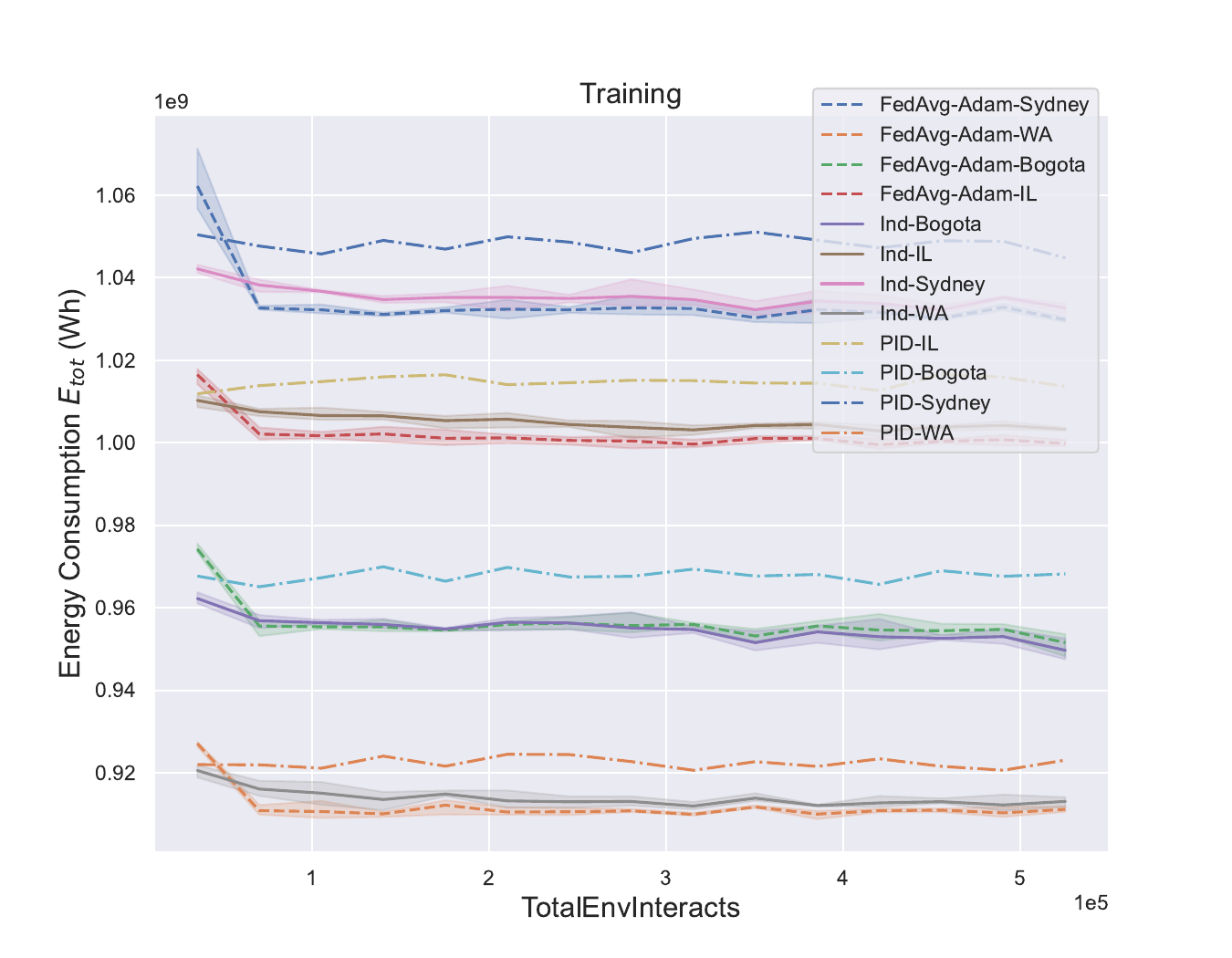}
         \caption{Energy consumption}
         \label{train:fedavg:ind:SydWABogIL:energy}
    \end{subfigure}
    \hfill
    \begin{subfigure}[b]{0.49\textwidth}
         \centering
         \includegraphics[width=\textwidth]{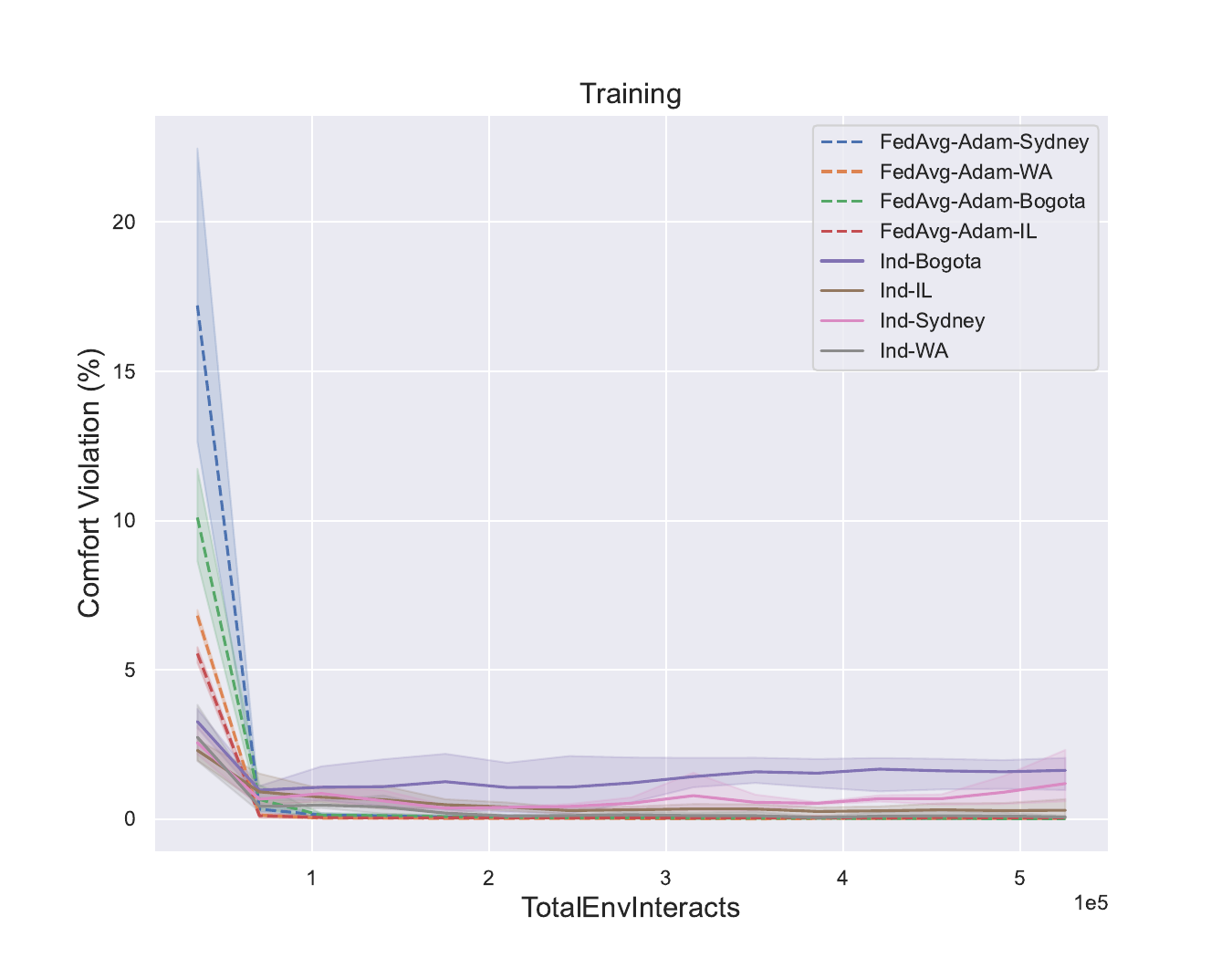}
         \caption{Comfort violation}
         \label{train:fedavg:ind:SydWABogIL:comfort}
    \end{subfigure}
    \caption{Progression of the energy consumption and comfort violation of FedAvg and independent agents on training environments Sydney, Bogota, WA and IL.}
    \label{train:fedavg:ind:SydWABogIL}
\end{figure}

In figures \ref{train:fedavg:ind:GraAntPA} and \ref{train:fedavg:ind:SydWABogIL} we show the training energy consumption and comfort violation curves for the environments omitted in section \ref{sec:TrainPerf}.

\end{document}